\documentclass[12pt]{article}
\usepackage{latexsym}
\usepackage{amssymb}
\usepackage[OT2,OT1]{fontenc}
\usepackage{graphicx}
\DeclareGraphicsExtensions{.eps, .jpg}
\def\cyr{%
\renewcommand\rmdefault{wncyr}%
\renewcommand\sfdefault{wncyss}%
\renewcommand\encodingdefault{OT2}%
\normalfont
\selectfont}
\DeclareTextFontCommand{\textcyr}{\cyr}
\def\p#1{\frac{\partial}{\partial #1}}

\newcommand{\Dye}{\mbox{\cyr   D}}

\newcommand{\C}{\mathbb{C}}
\newcommand{\CC}{\mathbb{C}}
\newcommand{\HH}{\mathbb{H}}
\newcommand{\CP}{\mathbb C\mathbb P}
\newcommand{\RP}{\mathbb R\mathbb P}
\newcommand{\PP}{\mathbb P}
\newcommand{\RR}{\mathbb R}

\newcommand{\Z}{\mathbb{Z}}
\newcommand{\ZZ}{\mathbb{Z}}

\renewcommand{\O}{\mathcal{O}}

\newcommand{\zz}{{\mathfrak z}}

\newtheorem{thm}{Theorem}[section]
\newtheorem{defn}[thm]{Definition}
\newtheorem{prop}[thm]{Proposition}
\newtheorem{cor}[thm]{Corollary}
\newtheorem{lem}[thm]{Lemma}
\newtheorem{conjecture}[thm]{Conjecture}
\newenvironment{remark}{\medskip \noindent {\bf Remark.}}{\hfill $\diamondsuit$\\}
\newenvironment{proof}{\medskip \noindent {\bf Proof.}}{\hfill $\blacksquare$\\}
\begin{document}
\title{Zoll Manifolds and Complex Surfaces}
\author{Claude LeBrun\thanks{Supported 
in part by  NSF grant DMS-0072591.}     ~and  
L.J. Mason}

\maketitle

\abstract{
We classify compact surfaces with 
 torsion-free 
affine connections
for which every geodesic is a simple closed curve.  In the process, 
we obtain  completely  new proofs of all the 
 major results \cite{beszoll}  concerning the Riemannian case. 
  In contrast to previous work, our 
approach is 
twistor-theoretic, and    depends  fundamentally on the 
fact that, up to biholomorphism, there is only one 
 complex structure on $\CP_2$.} 

\section{Introduction}

A  {\em Zoll metric} on a smooth manifold $M$  is a  Riemannian
metric $g$   whose geodesics are all simple closed curves of equal length. 
This terminology \cite{guillzoll} celebrates Otto Zoll's (now century-old)
 discovery \cite{zoll}
that $S^2$  admits many such metrics besides the obvious metrics of 
 constant curvature \cite{beszoll}. Indeed, in terms of  cylindrical coordinates
$(z,\theta)\in [-1,1]\times [0, 2\pi ]$, 
\begin{equation}
\label{ansatz}
g= \frac{\left[1+f(z)\right]^2}{1-z^2}dz^2 + (1-z^2)d\theta^2
\end{equation}
defines a Zoll metric on $S^2$
 for any smooth odd function 
$$
f: [-1 , 1] \to (-1, 1) , ~~ f(-z) = -f(z) 
$$
which vanishes at the end-points of the interval. A formal perturbation argument of 
Funk \cite{funk} later indicated that, modulo isometries and rescalings, 
 the {general}  Zoll metric
on $S^2$ depends on  one {odd} function 
$f: S^2 \to \RR$. This formal calculation was later turned
into a theorem by Guillemin \cite{guillzoll}, whose proof
depends, e.g.\ on an implicit function theorem 
of Nash-Moser type. Because the function $f$ is required to 
satisfy $f(-\vec{x})=-f(\vec{x})$, 
however,  these constructions never  give rise to  non-standard Zoll
metrics on $\RR\PP^2$. Indeed, the so-called Blaschke conjecture,  
proved by 
Leon  Green \cite{grezoll}, asserts that, up to isometries
and rescaling, the only  Zoll metric on $\RR\PP^2$
is the standard one. For an outstanding survey of these results, as well as 
 an exploration of their higher-dimensional Riemannian generalizations,
see \cite{beszoll}.

The aims of the present article are twofold. First of all, instead of
limiting ourselves to the study of  Riemannian metrics, we will more generally 
consider symmetric affine connections $\nabla$, and 
ask how many such connections on a given manifold $M$ have the property
that all of their geodesics are simple closed curves. In order
to make this  a sensible problem, however, one must first observe 
that for any   $1$-form $\beta$ on $M$, 
the  symmetric affine connection $\hat{\nabla}$ defined by 
$$\hat{\nabla}_{\bf u} {\bf v} = \nabla_{\bf u}{\bf v} + \beta ({\bf u}) {\bf v} + \beta ({\bf v}) {\bf u}$$ 
has exactly the same {\em unparameterized} geodesics
as the  connection $\nabla$; two connections
related in this manner are said to be {\em projectively equivalent},
and  obviously one should therefore only try to classify 
such connections  modulo projective equivalence.

Even in this rather general setting, our methods will allow us to 
 obtain results  very much like  
to the classical Riemannian results alluded to above. Indeed, in 
\S \ref{prelim}, we begin by showing that  the only
compact surfaces which admit Zoll projective  connections are $S^2$
and $\RR\PP^2$. In \S \ref{rip}, we then go on to show 
that,   modulo  diffeomorphisms, there is only one such
 projective class of  connections on $\RR\PP^2$.
Finally, in \S \ref{zoe}, we prove that 
there is a non-trivial moduli space of such 
projective classes on $S^2$, locally  parameterized by the  space of 
vector fields on $\RP^2$. 

But  even in the Riemannian case,   we seem to 
have something fundamentally 
new to contribute to the subject, as our proofs rest on foundations 
 completely
 different   
 from those used by  of our predecessors. 
Blaschke's unsuccessful approach to the problem of classifying
Zoll metrics on $\RP^2$ amounted to a direct attempt to identify the 
space of all geodesics with the standard dual projective plane
$\RP^{2*}$, the points of which which are by definition the   real projective lines 
$\RP^1$ in $\RP^2$.
The essence of  our method is to instead use {\em complex}, rather than real,
projective geometry to solve the problem. 
Indeed, we will 
 construct a  complex $2$-manifold
from any given   Zoll structure, modeled on the  dual complex projective plane $\CC\PP_2^*$.
The punch line of the proof is then that, up to biholomorphism, there is 
\cite{bpv,yau} 
only one complex structure on $\CP_2$. 
Our proof of the generalized Blaschke conjecture 
then proceeds by recognizing the points of $\RP^2$ as the set of those 
complex projective lines $\CP_1$ in this $\CP_2$
which are invariant under the action of a certain anti-holomorphic involution. 
By contrast, 
the flexibility of Zoll structure on $S^2$ arises because 
 the points in this case are instead represented 
by holomorphic disks with boundary on a totally real 
embedding of $\RP^2$ in $\CP_2$;  deformations of  this
embedding then correspond to deformations of
the Zoll structure. In this way, we are not only able to 
construct the general small deformation of the standard 
Zoll structure without recourse to Nash-Moser, 
 but, more importantly, we  are also
able to glean a significant amount of information 
regarding arbitrary Zoll structures, even when they are quite
far from the model case. 

Finally, by way of  an appendix, 
this article ends where it began,  with a discussion of the 
axisymmetric  case.  After all, since we have chosen to 
generalize Zoll's problem by focusing on projective structures, 
it is only fitting that we should also generalize
Zoll's  construction by  writing down all the 
axisymmetric Zoll projective structures on $S^2$ in closed form. 
In the process, we are able to show how the conceptual framework 
used in \S \ref{zoe}  can be implemented in concrete,
calculational terms. We  hope  that our discussion of this special case will not only
help clarify 
 our general approach, but also make it seem all the more compelling. 

\pagebreak 
 
\section{Zoll  Projective Structures} \label{prelim}

We  begin  by recalling the notion \cite{schouten} of projective equivalence of 
affine connections.
\begin{defn}
Two  torsion-free affine connections $\nabla$ and $\hat{\nabla}$ on a
 manifold $M$ are said to be {\em projectively equivalent}
if they have the same  geodesics, considered as {\em unparameterized} curves. 
\end{defn}
This condition may be re-expressed as
the requirement that 
$$
\hat{\nabla}_{\bf v}{\bf v}\propto {\bf v} ~~~ \Longleftrightarrow ~~~ {\nabla}_{\bf v}{\bf v}\propto {\bf v}.
$$
We therefore have  \cite{schouten}
\begin{prop}
Two $C^k$ symmetric affine connections $\nabla$ and $\hat{\nabla}$ 
are projectively equivalent iff
$$\hat{\nabla}_{\bf u} {\bf v} = \nabla_{\bf u}{\bf v} + \beta ({\bf u}) {\bf v} + \beta ({\bf v}) {\bf u}$$ 
for some $C^k$  $1$-form $\beta$. 
\end{prop}
Here a connection is said to be of differentiability class $C^k$ 
with respect to a fixed $C^{k+2}$ structure if 
the covariant derivative of any $C^{k+1}$ vector field is 
a $C^k$ tensor field; this is equivalent to requiring that 
the Christoffel
symbols
$$
\Gamma^j_{k\ell} = \left\langle dx^j , \nabla_{\p {x^k}}\p {x^\ell}\right\rangle
$$
are all $C^k$ functions in any admissible local coordinate system. 
We also note, in passing, that the symmetric (or torsion-free) condition 
employed here can been
imposed without any loss of generality; given  an {\em arbitrary} affine connection, one
can construct a unique torsion-free connection with precisely the same {\em parameterized}
geodesics by replacing the Christoffel symbols with their symmetrizations:
 $$
\Gamma^j_{k\ell} \rightsquigarrow \hat{\Gamma}^j_{k\ell}
=\frac{1}{2}\left(\Gamma^j_{k\ell} + \Gamma^j_{\ell k}\right)~.
$$
\begin{defn}
A $C^k$ {\em projective structure} on a smooth manifold is  the 
projective equivalence class $[\nabla]$ of some torsion-free $C^k$ affine
connection $\nabla$. 
\end{defn}
By definition, 
a  projective structure $[\nabla]$ on $M$  defines a certain family of geodesics;
these are to be thought of as abstract immersed curves in $M$,  without preferred parameterizations. Conversely, a projective structure is completely
specified once its geodesics are known.

In this paper, we will be interested in projective structures for
which every  geodesic is a simple closed curve. 

\begin{defn}
Let $\nabla$ be a $C^1$  torsion-free affine connection on a 
   smooth manifold $M$. We will say that the projective equivalence class
$[\nabla ]$ of 
$\nabla$   is a {\em  Zoll projective structure}
 if 
 the image ${\mathfrak C}$ of any maximal geodesic of $\nabla$
is an embedded circle $S^1\subset M$.
\end{defn}


%


If $c: (a,b) \looparrowright M$ is any immersed curve, its derivative 
$dc/dt$ is non-zero at every point, so that $[dc/dt ]$ is a well-defined
element of the {\em projectivized tangent bundle}
$$
\PP TM = (TM -0_M)/\RR^\times ;
$$
thus $t\mapsto [dc/dt]$ defines a  curve $\tilde{c} : \RR \to \PP TM$, called
the 
 {\em canonical lift}  of $c$. 
Given a  $C^k$ Zoll projective structure $[\nabla ]$ on $M$, the canonical lifts 
of its geodesics 
 give us a $C^k$ 
 foliation $\mathcal F$ of  $\PP TM$ by circles. Let  $N$ denote the
leaf space of this foliation. 

\begin{defn}
Let $(M,[\nabla ])$ be an $n$-manifold with $C^k$ Zoll projective structure. 
We will say that $[\nabla ]$ is {\em tame} if the corresponding 
foliation $\mathcal F$ of $\PP TM$
by lifted geodesics is 
 {\em locally trivial},  in the sense that   each  
leaf has a neighborhood which is $C^k$ diffeomorphic to 
$\RR^{2n-2} \times S^1$ in such a manner that every leaf 
 corresponds to a  circle of the form  $\{ pt \} \times S^1$. 
\end{defn}

These local trivializations give $N$ the structure
of a $C^k$  $(2n-2)$-manifold in a canonical manner, making the quotient map 
$\nu : \PP TM \to N$ into a $C^k$  submersion.  We will call the surface
$N$ the space of (undirected) geodesics of the tame  Zoll  projective structure
$[\nabla ]$. The situation is encapsulated by a diagram 
\setlength{\unitlength}{1ex}
\begin{center}\begin{picture}(20,17)(0,3)
\put(10,17){\makebox(0,0){$\PP TM$}}
\put(2,5){\makebox(0,0){$M$}}\put(18,5){\makebox(0,0){$N$}}
\put(15,12){\makebox(0,0){$\nu$}}
\put(5,12){\makebox(0,0){$\mu$}}
\put(11,15.5){\vector(2,-3){6}}
\put(9,15.5){\vector(-2,-3){6}}
\end{picture}\end{center}
which we shall refer to as the (real) {\em  double fibration} of $[\nabla ]$. 
Here $\mu : \PP TM \to M$ of course denotes the bundle projection. 
Notice that, by construction, the tangent spaces of the fibers of 
$\mu$ and $\nu$ are everywhere linearly independent:
$$(\ker \mu_*)\cap (\ker \nu_* ) =0.
$$
Moreover, the restriction of  $\nu$  to any fiber  of $\mu$
gives us an embedding  $\RR\PP^{n-1}\hookrightarrow N$.

Fortunately, as we will show in Theorem \ref{bingo} below,
this  desirable picture is   applies 
to every compact Zoll surface. 
A key step in this direction is the following: 

\begin{prop} \label{meow} 
Any   Zoll projective structure $[\nabla ]$ 
on a  compact  {\em orientable} surface $M^2$  is tame. 
\end{prop}
\begin{proof}
Because $M$ is assumed to be a compact surface, 
 $\PP TM$ is a compact $3$-manifold, and 
the Zoll projective structure $[\nabla ]$ 
gives us a foliation $\mathcal F$ of  $\PP TM$ by circles. 
However, a theorem of
Epstein \cite{epstein} asserts that any
foliation of a compact $3$-manifold by circles is a 
Seifert fibration. Thus any leaf of  $\mathcal F$ has 
a basis of neighborhoods modelled on 
$$(\CC \times S^1)/\ZZ_m,$$
where the $\ZZ_m$ action on $\CC \times S^1\subset \CC^2$
is generated by $(z_1, z_2)\mapsto (e^{2\pi i\ell / m}z_1, e^{2\pi i /m}z_2)$,
for some integer $\ell$. All we therefore need to show is that
no leaf is non-trivially covered by nearby leaves.

%

Now, because we have assumed that $M$ is orientable, any geodesic circle 
${\mathfrak C}$ has a tubular neighborhood diffeomorphic to the cylinder
$S^1 \times \RR$.  Moreover, by Epstein's result, 
the lift of  ${\mathfrak C}$ to $\PP TM$  has a standard neighborhood
whose projection to $M$  is contained in the given cylindrical neighborhood. 
Thus, 
any geodesic circle ${\mathfrak C}'$  with initial point
and tangent 
sufficiently close to those of ${\mathfrak C}$ will remain within our 
cylindrical neighborhood, and indeed will do so in such a manner that  the projection ${\mathfrak C}'\to {\mathfrak C}$
induced by $S^1 \times \RR\to S^1$ has non-zero derivative everywhere, and so will be 
 a covering map. 
However,  our tubular neighborhood 
 $S^1 \times \RR$ can be identified with  
  $\RR^2 -0$ in such a manner that 
${\mathfrak C}$ becomes the unit circle, and the degree of the covering 
becomes the winding number of ${\mathfrak C}'$ around the origin.
But since ${\mathfrak C}'$ has been transformed into an embedded curve
in the plane, the Jordan curve theorem tells us that 
its winding number around the origin has absolute value 
$\leq 1$. Thus  the covering map in question must have degree $1$.
The associated foliation $\mathcal F$ of  $\PP TM$  is therefore  trivial in a neighborhood of 
 the lift of ${\mathfrak C}$. 
%
\end{proof}

Next, we wish to determine precisely which compact surfaces
admit  Zoll projective structures. Our solution to this problem begins with 
 the following   simple observation:

\begin{lem}\label{cover}
Let $[\nabla]$ be a tame Zoll projective structure on an $n$-manifold
$M$. Let $\varpi : \tilde{M}\to M$ be the universal cover of $M$.
Then $[\varpi^*\nabla ]$ is a tame Zoll projective structure on
$\tilde{M}$. 
\end{lem}

\begin{proof}
If $(M,[\nabla])$ is a tame Zoll manifold, all the 
lifted geodesics are freely homotopic embedded circles
in $\PP TM$; this is true because $\PP TM$
is connected, and  is the union 
of `trivializing' open sets for the foliation $\mathcal F$, in which all the  circular leaves are 
 freely homotopic. 
Hence all the 
geodesic circles in $M$ are freely homotopic. Moreover, 
by considering the geodesic circles through a given point $p\in M$,
one obtains a base-point homotopy between 
any geodesic circle ${\mathfrak C}\subset M$ and its reverse-parameterized
version $\overline{{\mathfrak C}}$. Hence ${\mathfrak C}$ either represents  an   element of order 
$1$ or $2$ in $\pi_1 (M, p)$. Thus either ${\mathfrak C}$ or a 2-fold cover $\hat{{\mathfrak C}}\to {\mathfrak C}$
lifts to the universal cover $\tilde{M}$ as  an embedded circle,
and this circle is geodesic with respect to the pull-back
connection $\varpi^*\nabla$. 
Acting on each such lift by the action of
$\pi_1(M)$, we thus see that every geodesic of $(\tilde{M}, [\varpi^* \nabla ])$
is an embedded circle., and  
$[\varpi^*\nabla ]$ is therefore a Zoll projective structure on 
$\tilde{M}$.

It remains to show that $[\varpi^*\nabla ]$ is tame. To see this, 
first observe that $\mathcal F$ of $\PP TM$  
pulls back to the foliation $\hat{\mathcal F}$ of $\PP T\tilde{M}$ given 
by lifted geodesics of $[\varpi^*\nabla ]$. Moreover, the induced map
$\hat{\varpi}: \PP T\tilde{M}\to \PP TM$
is a covering map. If $U\subset \PP TM$ is any connected  open set, 
and if $\hat{U}\subset  \PP T\tilde{M}$ is any connected
component of $\hat{\varpi}^{-1}(U)$, then $\hat{\varpi}|_{\hat{U}}: \hat{U}\to U$
is also a covering map. But if $U$ is  a trivializing neighborhood for  $\mathcal F$,  then the finite cover $\hat{U}$ of $U\approx S^1\times \RR^{2n-2}$
will therefore provide a local trivialization of $\hat{\mathcal F}$.
 Since $\PP T\tilde{M}$ is covered by 
such neighborhoods, this shows that $(\tilde{M}, [\varpi^*\nabla ])$
is tame, as claimed. 
\end{proof}

This leads to  constraints on the topology of $M$. 

\begin{lem}\label{top}
Suppose that the $n$-manifold  $M$ 
admits a tame Zoll projective structure $[\nabla]$.
Then $M$ is compact, and has finite fundamental group.
Moreover, every two points $x$ and $x'$ of $M$ are joined by a 
geodesic of $\nabla$. 
\end{lem}

\begin{proof}
Choose an arbitrary point  $x\in M$. 
In $\PP TM$, consider the union 
$$\hat{X}=\nu^{-1}\left(\nu \left[\mu^{-1}(x)\right] \right)$$ of the lifts of geodesics through 
$x$. Then $\hat{X}$ is a compact
differentiable   $n$-manifold. 
But since  $\mu^{-1}(x)\subset \hat{X}$
is an $\RR\PP^{n-1}$ whose normal bundle is the universal
 line bundle, $\hat{X}$ may be    blown down
along $\mu^{-1}(x)$
to produce a new compact differentiable  $n$-manifold\footnote{We remark in passing that 
it is not difficult to show that ${X}$ is always diffeomorphic to $\RR\PP^n$.}  
 ${X}$.
Moreover, 
$\mu$ induces a differentiable map 
${\wp}:{X}\to M$. Indeed, 
if $\check{x}\in X$ denotes the point 
obtained by blowing down $\mu^{-1}(x)$, 
 then, in a neighborhood of $\check{x}$, 
$\wp$ is modeled on the exponential map of $\nabla$
near $0\in T_xM$. In particular, $\check{x}$
is a regular point of  $\wp$.  But, because $[\nabla ]$ is Zoll,  
a  geodesic circle can pass through $x$ 
{\em only  once},  
so it follows that ${\wp}^{-1}(x)=\{ \check{x}\}$. Thus $x$ is a regular value of 
the proper map $\wp$ with $\# {\wp}^{-1}(x)=1$.  This shows that   
the $\bmod$-$2$ degree
of the proper map $\wp$ is $1\in \ZZ_2$. In particular, 
 ${\wp}$ is onto, and 
 $M= {\wp}({X})$ is therefore compact. The very definition
of the surjective map $\wp$ now tells us that 
any point $x'$ of $M$ is joined to $x$ by 
some geodesic of $\nabla$.

Since the universal cover $\tilde{M}$ also admits a tame Zoll projective 
structure by Lemma \ref{cover}, the above argument now also shows that
$\tilde{M}$ is  compact. Hence the universal covering map $\varpi : \tilde{M}\to M$
is  finite-to-one, and $\pi_1(M)$ is therefore finite, as claimed. 
\end{proof}

Applying this  to the two-dimensional case, we  obtain the following:

\begin{prop} \label{class} 
A compact surface $M^2$ admits a  Zoll projective structure iff 
$M$ is diffeomorphic to either $S^2$ or $\RR \PP^2$. 
\end{prop}
\begin{proof}
By pulling the projective structure back to a double cover $\tilde{M}$ of
$M$ if necessary, we obtain a  Zoll
projective structure  on a compact orientable surface $\tilde{M}$, 
and this pulled-back structure is then 
  tame by 
Proposition \ref{meow}. This  forces $\tilde{M}$, and hence $M$, to have 
 finite fundamental group 
by Lemma \ref{top}. The classification of compact surfaces then tells us that
$M$  must be diffeomorphic to 
either $S^2$ or $\RR \PP^2$. Conversely, the Levi-Civita connection $\triangledown$
of  the standard, homogeneous metric determines a Zoll projective structure 
$[\triangledown ]$ on either of these spaces. 
\end{proof}

The following information thus becomes pertinent to our discussion: 

\begin{lem}\label{order}
If $M=S^2$, $|\pi_1(\PP TM)|=4$. If $M=\RR\PP^2$, 
$|\pi_1(\PP TM)|=8$.
\end{lem}
\begin{proof}
The unit bundle of $S^2$ may be identified with  $SO(3)$
by thinking of the first column of an orthogonal matrix as a point of
$S^2\subset \RR^3$, and the second  column as a unit tangent vector
at that point. 
Thus  $\PP TS^2$ may  be identified with $SO(3)/\ZZ_2$, where the 
$\ZZ_2$ action is generated by left multiplication by
$$\left[
\begin{array}{ccc}
1&0&0\\0&-1&0\\0&0&-1
\end{array}\right] ~.
$$
Lifting to the universal cover $Sp(1)=S^3\subset \HH^{\times}$ of $SO(3)$, we thus have
$\PP TS^2=Sp(1)/\ZZ_4$, where the $\ZZ_4$ is generated by $i$. 
Hence $\pi_1(\PP TS^2)\cong \ZZ_4$ has order $4$, as claimed. 

The antipodal map on $S^2$ acts on the unit tangent bundle 
via
$$
\left[
\begin{array}{ccc}
-1&0&0\\0&-1&0\\0&0&1
\end{array}\right] \in SO(3) ~,
$$
and this lifts to $Sp(1)$ as $\pm k$. 
Thus $\PP T\RR\PP^2= Sp(1)/ \{ \pm 1 , \pm i , \pm j , \pm k\}$, and hence 
$\pi_1(\PP T\RR\PP^2)\cong \{ \pm 1 , \pm i , \pm j , \pm k\}$
 has order $8$, 
as claimed. 
\end{proof}

In particular, $\pi_1 (\PP TM^2 )$ must be finite.  Hence: 

\begin{prop}
Let $(M,[\nabla ])$ be a  compact surface with tame Zoll projective structure. Then its space $N$ 
of unoriented geodesics is diffeomorphic to $\RP^2$.
\end{prop}
\begin{proof}
The group homomorphism  
$$\nu_{\natural}: \pi_1 (\PP TM)\to \pi_1 (N)$$
induced by the fibration $\nu$ is surjective, since  
 each fiber of $\nu$ is path connected. But Proposition \ref{class} and 
Lemma \ref{order} together 
tell us  that 
$\PP TM$ has finite fundamental group. Hence $\pi_1(N)$ is finite, 
and the classification of $2$-manifolds therefore tells us that $N$ must be 
diffeomorphic to either 
$S^2$ or $\RR\PP^2$. 
But we also know that  $N$ is not simply connected, since it has a non-trivial cover
$\tilde{N}$, given by the space of  {\em directed} geodesics of $[\nabla]$. 
This shows that $N\approx \RP^2$, as claimed. 
\end{proof} 

Next, we would like to understand the  topological structure of  the $S^1$-bundle 
$$\nu : \PP TM \to N .$$
Our method will  simultaneously allow us to analyze the conjugate points of the 
projective structure $[\nabla ]$. Let us thus begin by recalling the notion of a Jacobi field.

If $\nabla$ is a connection on a manifold $M$, and if 
$c: (a,b) \to M$ is an affinely  parameterized 
geodesic of $\nabla$, then a 
 {\em Jacobi field} 
along   $c$
is by definition a vector field 
${\bf y}\in \Gamma (c^*TM)$ 
along $c$ 
which satisfies the linear differential equation
$$\nabla_{\bf v}\nabla_{\bf v}{\bf y}= R_{{\bf v}{\bf y}}{\bf v},$$
where $R$ denotes the  curvature tensor of $\nabla$, and where 
 the standard tangent vector
$${\bf v}=\frac{dc}{dt}$$
of our parameterized geodesic 
satisfies the {\em auto-parallel condition} 
\begin{equation}\nabla_{\bf v} {\bf v} =0.\label{parl}\end{equation}
It is not difficult to see that ${\bf y}$ is a Jacobi field iff 
it is  locally 
the joining vector field for a  $1$-parameter family of geodesics of $\nabla$.
More precisely,  for any $[a',b']\subset (a,b)$, there is an  $\varepsilon > 0$ and 
a differentiable map \begin{eqnarray*}
\hat{c}: [a',b']\times (-\varepsilon, \varepsilon )&\to& M\\
(t,u) &\mapsto& \hat{c}(t,u)
\end{eqnarray*}
with $\hat{c}(t,0)=c(t)$, such that, setting 
$$\tilde{{\bf v}}=  \frac{\partial \hat{c}}{\partial t}, ~~ \tilde{{\bf y}} = \frac{\partial \hat{c}}{\partial u},$$
one has 
$$\nabla_{\tilde{{\bf v}}}\tilde{{\bf v}}=0$$
and $$\tilde{{\bf y}}|_{u=0}= {\bf y}.$$

The notion of a Jacobi field is not actually projectively invariant, but there is a 
closely related concept which {\em is}.

\begin{defn} Let $[\nabla ]$ be a $C^1$ projective connection on $M$, and 
let ${\mathfrak C}\looparrowright M$  be any geodesic of $[\nabla ]$.
Then a section ${\mathfrak Y} $ of the normal bundle $TM/T{\mathfrak C}$ of ${\mathfrak C}$ will be called 
a {\em Jacobi class} on ${\mathfrak C}$ iff,
near any given point $p\in {\mathfrak C}$,  
$${\mathfrak Y} \equiv {\bf y} \bmod T{\mathfrak C}$$
for some locally defined Jacobi field ${\bf y}$.
\end{defn}
In other words, ${\mathfrak Y} $ is a Jacobi class iff it locally joins infinitesimally
separated {\em unparameterized} geodesics. Thought of this
way, it thus becomes immediately apparent that 
the notion of Jacobi class is projectively invariant. 

\begin{defn} Let $[\nabla ]$ be a $C^1$ projective connection on $M$, and 
let ${\mathfrak C}\looparrowright M$  be any geodesic of $[\nabla ]$.
We will say that two points $p,q\in {\mathfrak C}$ are {\em conjugate}
along ${\mathfrak C}$ iff there is a Jacobi class  ${\mathfrak Y} $ on ${\mathfrak C}$ 
with ${\mathfrak Y} (p)={\mathfrak Y} (q)=0$. 
 \end{defn}
Very roughly, conjugate points are thus the places where two 
infinitesimally separated geodesics of $[\nabla]$ meet.

Let us now make all  of this more explicit in the special case of $\dim M=2$. If 
${\mathfrak C}\looparrowright M$ is a  geodesic of an affine connection $\nabla$ on 
a surface $M$, the normal bundle $TM/T{\mathfrak C}$
is  a real line bundle $E\to {\mathfrak C}$. Since $T{\mathfrak C}\subset TM$ is
parallel, $\nabla$ defines a connection ${\mathfrak D}$ on $E$. 
Let us take an affine parameterization
$c: (a,b)\to {\mathfrak C}$, so that ${\bf v}=dc/dt$
satisfies (\ref{parl}). Let us then trivialize $c^*E\to (a,b)$  by means of $[{\bf e}]$,
where ${\bf e}\not\propto {\bf v}$ is a generic  parallel section of $c^*TM$, and 
where the brackets $[\cdot]$ indicate the equivalence class $\bmod  ~ T{\mathfrak C}$. 
Defining $\kappa : (a,b)\to \RR$ by 
$$\kappa =  r({\bf v},{\bf v}),$$
where $r_{ab}={R^c}_{acb}$
is the Ricci tensor of $\nabla$, we then have 
$$
R_{{\bf v}{\bf e}}{\bf v}\equiv -\kappa {\bf e} \bmod {\bf v}, 
$$
so that $y(t) {\bf e}\equiv {\bf y}\bmod T{\mathfrak C}$ for some Jacobi field 
${\bf y}$ iff $y : (a,b)\to \RR$ satisfies the second order linear differential equation 
\begin{equation}\label{jack}
\frac{d^2y}{dt^2}+ \kappa y=0.\end{equation}
More abstractly, (\ref{jack}) becomes
\begin{equation}\label{lack}
{\mathfrak D}_{\bf v}{\mathfrak D}_{\bf v}{\mathfrak Y} + r({\bf v},{\bf v}) {\mathfrak Y} =0\end{equation}
in terms of the connection ${\mathfrak D}$ induced on the normal bundle 
$E$, and this in turn generalizes to becomes 
\begin{equation}\label{quack}
{\mathfrak D}_{\bf v}{\mathfrak D}_{\bf v}{\mathfrak Y} -{\mathfrak D}_{\nabla_{\bf v}{\bf v}}{\mathfrak Y} + r({\bf v},{\bf v}) {\mathfrak Y} =0\end{equation}
if we drop the auto-parallel condition (\ref{parl}) on our tangent field ${\bf v}$.
Let us remark  that if 
 $\nabla$ is replaced by the projectively equivalent connection $\hat{\nabla}$
defined by 
$$\hat{\nabla}_{\bf u}{\bf v} = \nabla_{\bf u}{\bf v}+ \beta ({\bf u}) {\bf v} + \beta ({\bf v}) {\bf u} ,$$
one then has
\begin{eqnarray*}
\hat{\mathfrak D}_{\bf v}\hat{\mathfrak D}_{\bf v} {\mathfrak Y} &=& 
{\mathfrak D}_{\bf v} {\mathfrak D}_{\bf v} {\mathfrak Y}
+ 2 \beta ({\bf v}) {\mathfrak D}_{\bf v} {\mathfrak Y}
+ \left[{\bf v} \beta ({\bf v})+ \beta ({\bf v})^2 \right] {\mathfrak Y}  ,
\\
\hat{\mathfrak D}_{\hat{\nabla}_{\bf v}{\bf v}}{\mathfrak Y} &=& 
{\mathfrak D}_{{\nabla}_{\bf v}{\bf v}}{\mathfrak Y} + 2 \beta ({\bf v}) 
{\mathfrak D}_{\bf v}{\mathfrak Y} 
+ \left[\beta (\nabla_{\bf v}{\bf v})+ 2 \beta ({\bf v})^2\right] {\mathfrak Y} , \\
\hat{r}({\bf v},{\bf v}) &=& r ({\bf v},{\bf v}) + (n-1) 
\left[ \beta (\nabla_{\bf v}{\bf v}) - {\bf v} \beta ({\bf v}) + \beta({\bf v})^2\right] , 
\end{eqnarray*}
so that blind, brute-force calculation does indeed show  that 
(\ref{quack}) is   projectively invariant   in dimension $n=2$, as   previously 
deduced by pure thought.

Now the vector space of solutions of (\ref{jack}) is  two dimensional, corresponding to  choices of 
$y$ and $y'$ at an arbitrary base-point of the interval $(a,b)$. Let
$\{ y_1, y_2\}$ be an arbitrary basis for this solution space, and 
consider the Wronskian 
$$W(t)= \left|\begin{array}{cc}
y_1(t)&y_1'(t)\\y_2(t) & y_2'(t)
\end{array}
\right|=y_1y_2'-y_2y_1'.
$$
The differential equation (\ref{jack})
then tells us that 
\begin{eqnarray*}
\frac{dW}{dt}&=&y_1'y_2'+ y_1y_2''-y_2'y_1'-y_2y_1''\\
&=&y_1(-\kappa y_2)- y_2(-\kappa y_1) =0,
\end{eqnarray*}
so that $W(t)$ is constant. Moreover,  this
constant must be non-zero, since $y_1$ and $y_2$ have linearly 
independent initial values at the base-point. 
 The map 
 \begin{eqnarray*} 
 \phi : (a,b) & \longrightarrow &~~~~\RP^1\\ 
t~~~&\mapsto&  [y_1(t):y_2(t)]\end{eqnarray*} 
is therefore well defined for all $t$, since $y_1$ and $y_2$ cannot
 simultaneous vanish. Moreover, $\phi$ is an {\em immersion}, since
$$\frac{d}{dt}\left( \frac{y_1}{y_2}\right) = \frac{W}{y_2^2}
~~~ \mbox{  and  } ~~~\frac{d}{dt}\left( \frac{y_2}{y_1}\right) = -\frac{W}{y_1^2}
$$
are never zero.  Geometrically, $\phi$ may be interpreted as sending $x\in (a,b)$
to the set of Jacobi classes ${\mathfrak Y}$ with ${\mathfrak Y} (x)=0$, since a Jacobi class 
$$y(t) = \lambda_2 y_1(t) -\lambda_1 y_2 (t)\not\equiv 0$$
vanishes at $x$ iff $[ \lambda_1 : \lambda_2] =\phi (x) := [y_1(x) : y_2(x)]$. 
In particular, two points are conjugate along $c\left[(a , b )\right]$ 
iff they have the same image under $\phi$.

For a tame $C^k$ Zoll projective structure 
$[\nabla ]$ on a surface $M^2$, there are two linearly independent  
 Jacobi
classes defined along the entirety of any  closed geodesic ${\mathfrak C}$; indeed, if $y\in N$
represents ${\mathfrak C}$ in the space of geodesics, $T_yN$ is naturally
in one-to-one correspondence with the space of Jacobi classes
along ${\mathfrak C}$ via $\mu_*\circ (\nu_*)^{-1}$. 
The  above construction thus gives
us a $C^{k+1}$ covering map $\phi : {\mathfrak C}\to \RP^1$
for every geometrically closed geodesic ${\mathfrak C}$, and this map is uniquely defined
modulo the action of $SL(2, \RR )$ on $\RP^1$. 
The order  of the covering $\phi : {\mathfrak C} \to RP^1$ will be called the 
{\em conjugacy number}
of the geodesic, since it exactly counts how many points of ${\mathfrak C}$ are  
conjugate to $x\in {\mathfrak C}$, of course including  $x$ itself. We will 
now see that this number actually has a rather deeper meaning.

\begin{prop}\label{blimey} 
Let $[\nabla ]$ be a tame $C^k$ Zoll projective connection, $1\leq k \leq \infty$,
 on a compact $2$-manifold $M$, and 
consider the $C^{k-1}$ map
 \begin{eqnarray*} 
\varphi :\PP TM &\longrightarrow& ~~~~\PP TN\\
z~~~~&\mapsto& \nu_* \left( \ker \mu_{*z}\right) ,
 \end{eqnarray*} 
where $\mu_*$ and  $\nu_*$ denote the derivatives of $\mu$ and $\nu$,
respectively. Then $\varphi$ is 
a  covering map. Moreover, the order of the covering $\varphi$ exactly equals 
 the conjugacy number of 
any closed geodesic ${\mathfrak C}\subset M$. In particular,
all the geodesics of  $[\nabla ]$ have the same conjugacy number. 
\end{prop}
\begin{proof}
Let us first notice that  we have a commutative diagram 
\setlength{\unitlength}{1ex}
\begin{center}\begin{picture}(36,19)(0,3)
\put(10,17){\makebox(0,0){$\PP TM$}}
\put(18,19){\makebox(0,0){$\varphi$}}
\put(18,5){\makebox(0,0){$N$}}
\put(26,17){\makebox(0,0){$\PP TN$}}
\put(15,12){\makebox(0,0){$\nu$}}
\put(21,12){\makebox(0,0){$\pi$}}
\put(11,15.5){\vector(2,-3){6}}
\put(25,15.5){\vector(-2,-3){6}}
\put(14,17){\vector(1,0){8}}
\end{picture}\end{center}
where $\pi$ denotes the relevant canonical projection. Moreover, 
since $N$ is by definition the leaf space of the 
  foliation $\mathcal F$, we also know that $\varphi$ maps  
each leaf of  $\mathcal F$ to a different
 fiber of $\pi$.

Now 
the tangent space of $N$ at any point can be canonically
identified with the space of Jacobi classes on the corresponding
geodesic in $M$. With this identification,   $\varphi$ then sends a point
of a  geodesic ${\mathfrak C}$ (identified, by lifting, with a leaf of $\mathcal F$) to the 
 set of Jacobi classes which vanish at that point. In other words, 
on each leaf of $\mathcal F$, thought of as a geodesic ${\mathfrak C}\subset M$
of $[\nabla ]$, $\varphi$ precisely coincides with the map $\phi$ described above. 
This shows that $\varphi$ immerses each leaf in $\PP TN$ as a 
fiber of $\pi$.  Since $\nu$ is a submersion, it follows, for $k\geq 2$,
that $\varphi_*$ is injective, and hence that $\varphi$ is a 
local diffeomorphism; for $k=1$, one instead may observe that 
 $\varphi$ must be injective on some
 neighborhood of any point, and so must be 
a local homeomorphism by the  open mapping theorem.
But since  $\PP TM$ is compact,   this implies that  $\varphi$ is  a  covering map. 
Moreover, the order of this covering is precisely the number
of points on a leaf of $\mathcal F$ which are sent to the same point of a fiber 
$\pi$. This shows that the order of covering $\varphi$ is precisely the
conjugacy number of any geodesic of $[\nabla]$. 
\end{proof}

If $X$ is any manifold,  let us use ${\mathbb S}TX$ to denote the 
sphere bundle $(TX- 0_X)/\RR^+$.  In other words, 
${\mathbb S}TX$ may be thought of as the set of unit tangent vectors
 for an arbitrary Riemannian metric on $X$. 

\begin{thm} \label{wooster}
If $[\nabla ]$ is any $C^k$ Zoll projective structure, $k \geq 1$, on $M\approx S^2$,
its conjugacy number is two, and 
there is a $C^{k-1}$ diffeomorphism  $\PP TM \approx {\mathbb S}TN$ such that 
$\nu$ becomes the canonical projection ${\mathbb S}TN \to N$.
Moreover, the real line bundle $\ker \mu_*$ over $\PP TM$
is trivial.
\end{thm}

\begin{proof} Let us first recall that Proposition \ref{meow} tells us 
that $[\nabla ]$ is tame. But now, 
with a nod to  Lemma \ref{order}, we see that 
the covering map 
 $\varphi : \PP TM \to \PP TN$
 has order 
$$\frac{|\pi_1(\PP TN )| }{ |\pi_1 (\PP TM )|}=
\frac{|\pi_1(\PP T\RR\PP^2 )| }{|\pi_1 (\PP TS^2 )|} =\frac{8}{4} =2,$$
and the conjugacy number is therefore $2$, by Proposition \ref{blimey}. 

Now notice that the real line bundle $\ker \mu_*$ over $\PP TS^2$
is trivial. Indeed, after the choice of a metric and orientation, 
$\PP TS^2$ can be identified with the $SO(2)$ bundle of 
oriented orthonormal frames divided by $\langle -{\bf 1}\rangle \subset SO(2)$, 
and carries an induced ${\mathfrak so}(2)$ action which trivializes 
$\ker \mu_*$.. 

Imitating our  construction of $\varphi$, we now  obtain a diagram 
\setlength{\unitlength}{1ex}
\begin{center}\begin{picture}(36,19)(0,3)
\put(10,17){\makebox(0,0){$\PP TM$}}
\put(18,19){\makebox(0,0){$\hat{\varphi}$}}
\put(18,5){\makebox(0,0){$N$}}
\put(26,17){\makebox(0,0){${\mathbb S}TN$}}
\put(15,12){\makebox(0,0){$\nu$}}
\put(21,12){\makebox(0,0){$\wp$}}
\put(11,15.5){\vector(2,-3){6}}
\put(25,15.5){\vector(-2,-3){6}}
\put(14,17){\vector(1,0){8}}
\end{picture}\end{center}
by defining $\hat{\varphi}(z)= \RR^+\nu_*({\bf v}_z)$; here $\wp: {\mathbb S}TN \to N$
of course denotes the canonical projection. 
Now   $\hat{\varphi}$ is a  covering map, since it  lifts
$\varphi$. 
But 
$$\frac{|\pi_1(S TN )| }{ |\pi_1 (\PP TM )|}=
\frac{|\pi_1(S T\RR\PP^2 )| }{|\pi_1 (\PP TS^2 )|} =\frac{4}{4} =1,$$
so it now follows that $\hat{\varphi}$ is a homeomorphism if $k=1$, and 
a diffeomorphism if $k\geq 2$. 
\end{proof}

This finally allows us to definitively  dispense with the tame condition.

\begin{thm} \label{bingo}
Any $C^1$ Zoll projective structure  on a compact surface
$M^2$ is tame.  
\end{thm}
\begin{proof} 
Proposition \ref{meow} covers   the orientable case, so we  may assume 
henceforth that 
 $M$ is non-orientable.  Proposition \ref{class} then tells us that 
$M$ is diffeomorphic to $\RP^2$, so we have 
$M= \tilde{M}/\langle a\rangle$, where 
$\tilde{M}\approx S^2$, and where $a :  \tilde{M} \to  \tilde{M}$
corresponds to the  antipodal map on $S^2$. 
Any Zoll projective structure  on $M$ then pulls
back to a tame Zoll projective structure on $\tilde{M}$
each of whose geodesics is  sent to some geodesic by $a$.
If $\tilde{N}\approx \RP^2$ is the space of 
unoriented geodesics of $\tilde{M}$, then $a$ thus induces
a diffeomorphism $\hat{a}: \tilde{N}\to \tilde{N}$. We claim  
that $\hat{a}$ is in fact the identity. 

Suppose not. Then $\hat{a}$ generates a non-trivial $\ZZ_2$ action.
But any  action by a finite group of diffeomorphisms is isometric with respect to some 
 Riemannian metric, and so has fixed-point set
consisting of a disjoint union of closed submanifolds. In our case, 
the fixed-point set would be  a finite  union  of disjoint circles and points. 
Moreover, the quotient $\tilde{N}/\ZZ_2$ would have Euler characteristic 
$$\chi (\tilde{N}/\ZZ_2) = \frac{\chi (\tilde{N})+m}{2} =  \frac{1+m}{2} ,$$
where $m$ is the number of isolated fixed points. Since the Euler
characteristic is an integer, this shows that $\hat{a}$
has at least one  isolated fixed point. At such an isolated fixed point, the
derivative of $\hat{a}$ must be $-1$, as this is   the unique  order-$2$ element of 
$O(2)$  with trivial  $+1$-eigenspace. 

Back in $\tilde{M}$, this fixed point would correspond to a geodesic circle 
$\mathfrak C$ with $a ({\mathfrak C})={\mathfrak C}$
along which $a_*$ induced the action 
${\mathfrak Y} \mapsto -{\mathfrak Y}$ on the vector space of Jacobi classes. 
In particular, the zero locus of a Jacobi class ${\mathfrak Y}\not\equiv 0$
would  necessarily be sent to itself by $a_*$. But since $\tilde{M}$ has
conjugacy number $2$ by Theorem \ref{wooster}, and since $a$ has no fixed points, 
this means that $a$ acts 
on $\mathfrak C$ by sending each point to the  unique other point
to which it is conjugate. Now trivialize the normal bundle $E= T\tilde{M}/T{\mathfrak C}$,
so that we can talk about whether a non-zero element of $L$ is `positive'
or `negative'. Then, since any Jacobi class ${\mathfrak Y}\not\equiv 0$
meets the zero section of $E$ transversely in exactly $2$ points, 
the subsets of $\mathfrak C$ given by ${\mathfrak Y} > 0$ and
${\mathfrak Y} < 0$ are necessarily intervals, and are necessarily 
interchanged by the fixed-point-free map $a$. But since 
${\mathfrak Y} \mapsto -{\mathfrak Y}$
under $a_*$, this shows that $a_*$ acts on the normal bundle $E$
in an {\em orientation-preserving} manner. Moreover, $a_*$ is also
orientation-preserving on $T{\mathfrak C}$, since $a : {\mathfrak C}
\to {\mathfrak C}$ has no fixed point. 
 Hence $a$ acts on $\tilde{M}$ in an orientation-preserving 
manner --- contradicting the fact that, by construction,  $a$ is an orientation-reversing map!

This contradiction shows that $\hat{a}$ must be the identity on $\tilde{N}$. Hence 
$a_*$ induces an action on 
$\PP T\tilde{M}$ which sends each leaf to itself, and  holonomy around any leaf  
in $\PP TM$ is therefore trivial. Hence the given Zoll projective structure on
$M\approx \RP^2$ is tame, as claimed. 
 \end{proof}

In particular, it now makes sense to talk about the 
conjugacy number of any Zoll projective structure 
on $\RP^2$. 

\begin{thm} \label{jeeves} 
If $[\nabla ]$ is any $C^k$  Zoll projective structure, $k\geq 1$, 
 on $M\approx \RR\PP^2$,
its conjugacy number is $1$. Moreover, 
there is a $C^{k-1}$ diffeomorphism $\PP TM \approx \PP TN$ such that 
$\nu$ becomes the canonical projection $\PP TN \to N$, and such that 
 $\ker \mu_*\to\PP TM$
becomes  the `tautological' real line bundle
 $L \to \PP TN$, whose  frame bundle bundle is the principal $\RR^\times$-bundle 
$(TN - 0_N )\to \PP TN$. 
\end{thm}

\begin{proof}
Since $[\nabla ]$ is tame by Theorem \ref{bingo}, we are free to  
consider the covering map $\varphi : \PP TM \to \PP TN$
of Proposition \ref{blimey}. By construction,  the tautological
line bundle $L\to \PP TN$ then satisfies $\varphi^*L=\ker \mu_*$.
Since 
the order of this covering is $$\frac{|\pi_1(\PP TN )| }{ |\pi_1 (\PP TM )|}=
\frac{|\pi_1(\PP T\RR\PP^2 )| }{|\pi_1 (\PP T\RR\PP^2 )|} =1,$$
we conclude that  $\varphi$ is a homeomorphism, and  
 the conjugacy number is therefore $1$ by Proposition \ref{blimey}. 
Moreover,  the same argument also shows that $\varphi$ is actually a  diffeomorphism
if $k\geq 2$.
 \end{proof}

\begin{cor} \label{lcapl}
For  any  Zoll projective structure  $[\nabla ]$ on $M\approx \RR\PP^2$,
any two distinct points are joined by a unique geodesic circle $\mathfrak C$. 
\end{cor} 
\begin{proof}
As in the proof of Lemma \ref{top}, 
let 
$$\hat{X}=\nu^{-1}\left(\nu \left[\mu^{-1}(x)\right] \right)$$ be the 
union of the lifts of geodesics through 
$x$. Then $\hat{X}$ is a compact
differentiable   surface and may be    blown down
along $\mu^{-1}(x)$
to produce a new smooth compact  surface $X$.   
Since $\hat{X}$ is a circle bundle over the circle $\ell_x =  \nu \left[\mu^{-1}(x)\right]$,
and since a neighborhood of  $\mu^{-1}(x)$ is a M\"obius band $B$, 
it follows that 
 $X$  contains a  M\"obius band $B'=\hat{X}-B$, 
and hence  is not orientable. 

 On the other hand, Theorem \ref{jeeves}
tells us that each geodesic in $M$ has conjugacy number $1$, and
hence  no point $x'\neq x$ is conjugate to $x$ along any geodesic. 
Hence the canonical projection $\hat{X}\to M$ is an  immersion 
away from $\mu^{-1}(x)$, and the induced map $\wp :X\to M$ is therefore
an immersion everywhere.  Since $X$ is compact, $\wp$  is therefore
 a covering map. Since $X$ is not simply connected and 
$\pi_1 (M) = \ZZ_2$, it follows that $\wp$ is a one-to-one and onto. 
But, by the very definition of $\wp$, this means that there is one 
and only one geodesic between $x$ and any other point $x'\neq x$
in $M$. 
\end{proof}

\begin{cor} \label{lowercase}
Let $(M^2, [\nabla ])$ be a compact surface with Zoll projective structure.
Let ${\mathfrak C}\subset M$ be any geodesic circle.  
Then the following conditions are equivalent:
\begin{itemize}
\item $\langle w_1 (M) , [{\mathfrak C}]\rangle = 1 \in \ZZ_2$;
\item the conjugacy number of ${\mathfrak C}$ is odd;
\item $M$ is not orientable;
\item $M$ is diffeomorphic to $\RP^2$.
\end{itemize}
\end{cor}

\begin{proof}
At points where a Jacobi class ${\mathfrak Y}\not\equiv 0$ 
vanishes along ${\mathfrak C}$, the covariant derivative ${\mathfrak D}_{\bf v}{\mathfrak Y}$
must be nonzero, since ${\mathfrak Y}$ satisfies (\ref{jack}). Thus the 
mod-$2$ reduction of the conjugacy
number of ${\mathfrak C}$ calculates $\langle w_1 (E) , [{\mathfrak C}]\rangle$, where 
$E=TM/T{\mathfrak C}$ is the normal bundle, and this of course coincides with 
 $\langle w_1 (M) , [{\mathfrak C}]\rangle := \langle w_1 (TM) , [{\mathfrak C}]\rangle$, 
since $T{\mathfrak C}$ is trivial. But Theorems \ref{jeeves} and \ref{wooster} 
tell us that the only possible values of the conjugacy number are $1$ and $2$, and 
that the value of the conjugacy number determines whether $M$ is diffeomorphic
to $\RP^2$ or $S^2$. 
\end{proof}

The same argument also yields the following: 

\begin{cor}\label{uppercase}
Let $(M^2, [\nabla ])$ be a compact surface with Zoll projective structure.
Let ${\mathfrak C}\subset M$ be a geodesic circle.  
Then the following conditions are equivalent:
\begin{itemize}
\item $\langle w_1 (M) , [{\mathfrak C}]\rangle = 0\in \ZZ_2$;
\item the conjugacy number of ${\mathfrak C}$ is even;
\item $M$ is orientable;
\item $M$ is diffeomorphic to $S^2$.
\end{itemize}
\end{cor}

Let us now take a moment  to  compare our definitions with those 
previously used by others in the Riemannian context \cite{beszoll,guillzoll}. 

\begin{prop}\label{others}
Let $(M^2,g)$  be a compact surface with $C^k$ Riemannian metric, $2\leq k\leq\infty$. Let
$\triangledown$ be the Levi-Civita connection of $g$. Then $[\triangledown]$ is a 
$C^{k-1}$  Zoll projective
structure on $M$ iff the geodesics of $g$ are all simple  
closed curves of equal length. 
\end{prop} 
\begin{proof}
If $[\triangledown]$ is  a Zoll projective structure,  Theorem \ref{bingo} then
tells us it is tame,
and its geodesic circles are therefore freely homotopic to one
another {\em through  geodesic circles}. 
But the affinely parameterized closed geodesics of $g$ are
precisely those differentiable maps $c: S^1 \to M$ which are critical points of 
the energy functional
$$E( c ) = \int_{S^1} g(  c^\prime (t) ,  c^\prime (t) ) dt ~;$$
 thus the energy is necessarily constant for any $1$-parameter
family of closed geodesics. This shows that  the 
geodesic circles of $g$ must all have equal 
energy,  and hence equal 
length. 
\end{proof}

We conclude this section with an aside which  plays no r\^ole whatsoever 
in what follows, but which, in light of Proposition \ref{others},
 has a certain intrinsic interest. 
Given a  Zoll projective structure $[\nabla ]$ on a compact 
surface $M$, it is natural to  ask whether there is a connection $\nabla$ 
representing $[\nabla ]$ such that every affinely parameterized geodesic
 is {\em
periodic}. 
The answer is affirmative.

\begin{prop}
If $[\nabla ]$ is any Zoll projective structure on a compact surface $M^2$,
then there is a symmetric affine connection $\nabla\in [\nabla ]$ 
for which each affinely parameterized geodesic extends as a 
periodic function $c: \RR \to M$. 
\end{prop}

\begin{proof} If
$M=S^2$, let 
$\omega$ be an arbitrary area form on $M$, and let $\nabla$ be \cite{schouten} the 
unique connection in the equivalence class such that $\nabla \omega =0$.
If $c: [a,b]\to M$ is an affine parameterization of a geodesic of $\nabla$,
with $c(b)=c(a)$ and $c^{\prime}(b)=\lambda c^{\prime}(a)$, then
 any   parallel vector field ${\bf e}$ along $c$
 must satisfy ${\bf e}(b)=  \lambda^{-1} {\bf e}(a)\bmod c^{\prime}$. 
Now the Zoll  condition guarantees the 
existence of  a two-parameter family of  solutions of  (\ref{lack}) which  satisfy the
"periodicity" condition 
$${\mathfrak Y}|_{c (b)} = {\mathfrak Y} |_{c(a)}, 
~~ {\mathfrak D}{\mathfrak Y} |_{c (b)}  = {\mathfrak D} {\mathfrak Y}  |_{c(a)} .$$
 Every solution of (\ref{jack})  
must therefore 
satisfy $$y(b)=\lambda y(a), ~~~ y^\prime (b) = \lambda^2 y^\prime (a).$$ 
Hence the Wronskian $W= y_1y_2^\prime - y_2y_1^\prime$ of two linearly
independent solutions of (\ref{jack}) must  satisfy 
$W(b)=\lambda^3 W (a)$. But $W$ is constant! Thus $\lambda = 1$, and 
the given geodesic is therefore periodic. But this argument applies to {\em any}
geodesic on $M$. 
Hence every geodesic of 
the chosen connection $\nabla$  is   periodic, and the claim follows if $M=S^2$. 

The case of $\RR\PP^2$ now follows easily; one simply  takes the area form 
$\omega$ on $S^2$  to be anti-invariant under the antipodal map  
$a : S^2\to S^2$, and then notices that the corresponding connection 
$\nabla$ then descends to 
$\RR\PP^2$. 
\end{proof}


\section{The Blaschke  Conjecture Revisited}
\label{rip}

If $[\nabla ]$ is a Zoll projective structure on a compact surface $M$, 
we saw in \S \ref{prelim}  that its space of unoriented geodesics $N$
is diffeomorphic to $\RP^2$.  Now notice that $N$ also comes
equipped with a family 
$$\ell_x = \nu [ \mu^{-1} (x) ]$$
of embedded circles $\ell_x\subset N$, $x\in M$. (For any given $x\in M$, 
this is to say that 
 $\ell_x$  consists  precisely of the geodesics  passing through $x$.)
If we were simply given $N$ and this family of curves, we could then 
completely reconstruct the given projective structure on $M$. 
Indeed, $M$ could be redefined as the parameter space or `moduli space' 
of these curves $\ell_x$, and the geodesics ${\mathfrak C}_y \subset M$  would then become the
set of curves $\ell_x$ passing through some given point $y\in N$. 
The utility of this point of view might seem to be rather questionable, however,
as there is no obvious  geometric structure one might impose on $N$ 
in order to keep track of which embedded circles $\ell \subset N$ are to be the 
elements of the family $\{ \ell_x \}_{x\in M}$.  However, our main observation,
extrapolated from  a twistor correspondence due to 
Hitchin \cite{hitproj} and the first author \cite{lebthes}, 
is that one {\em can} naturally keep track of these curves by `complexifying' the
picture, and embedding $N$ in a complex $2$-manifold ${\mathcal N}$.

Let us suppose we are given a $C^2$ Zoll projective structure $[\nabla ]$
on $M=\RP^2$.
 Consider the $\CP_1$-bundle 
$$\PP T_\CC M = \left( \CC \otimes TM - 0_M
\right)
/\CC^\times,$$
and observe that the circle bundle
$$\PP T M = \left(  TM - 0_M
\right)
/\RR^\times $$
is a hypersurface in the $4$-manifold $\PP T_\CC M$. 
For brevity, we introduce the notation
$${\mathcal Z}= \PP T_\CC M, ~~ Z = \PP TM.$$
Because each fiber of $\PP T_\CC M$ has a canonical 
complex structure $J^\parallel$, the normal bundle of $\PP T M\subset \PP T_\CC M$
is just $J^\parallel (\ker \mu_*)$, where $\mu : \PP T M \to M$
is the bundle projection. Now recall that  
our Zoll projective structure gives us a foliation 
$\mathcal F$ of $\PP T M$ by circles,  and the leaves of 
$\mathcal F$ are precisely the fibers of a $C^2$ submersion 
$\nu : \PP T M \to N\approx \RP^2$. Moreover, 
Theorem \ref{jeeves} tells us that there is a $C^1$ 
diffeomorphism $\varphi : \PP T M \to \PP T N$ 
such that the real line bundle $\ker \mu_*$ 
becomes the pull-back $\varphi^* L$ of the 
tautological line bundle $L\to  \PP T N$.
The latter line bundle is by definition a sub-bundle of 
$\pi^* TN$, where $\pi : \PP TN \to N$ is the canonical 
projection; namely, for any non-zero vector ${\bf v}\in T_yN$, 
the fiber over $[{\bf v} ] \in \PP TN$ is 
$L_{[{\bf v}]}= \mbox{span} ({\bf v}) \subset T_yN$. In particular, 
there is a tautological $C^1$ `blowing down' map 
$\psi : L\to TN$ which is a diffeomorphism 
away from  the zero section $\PP TN$ 
of $L$, but collapses this zero section to the zero section
$N$ of $TN$ via $\pi : \PP TN \to N$. 
On the other hand,  the tubular neighborhood theorem tells
us that $Z=\PP T M$ has a neighborhood $\hat{\mathcal V}$ in ${\mathcal Z}=\PP T_\CC M$
which is $C^\infty$ diffeomorphic to the total space of 
$J^{\parallel}\ker \mu_*$, in such a manner that the derivative along 
$Z$ is the identity. Letting $\mathcal V$ denote the total space of 
$TN=T\RP^2$, we then have a $C^1$ map $\tilde{\psi} : \hat{\mathcal V}\to {\mathcal V}$
which corresponds to $\psi$ via our $C^1$ diffeomorphism $J^{\parallel}\ker \mu_*\to L$. 
We may now define a 
 new $C^1$ compact $4$-manifold 
$${\mathcal N} = {\mathcal U} \cup_{\tilde{\psi}} {\mathcal V}$$
by gluing together ${\mathcal U} := {\mathcal Z}-Z$
and ${\mathcal V}=TN$ via $\tilde{\psi}$. 
 By construction, 
we also  have a $C^1$ 
`blowing down' map 
$$\Psi : {\mathcal Z} \to {\mathcal N},$$
given by the identity on ${\mathcal U}$ and by $\tilde{\psi}$ on $\hat{\mathcal V}$. 





If we suppose that $[\nabla ]$ is
$C^k$ for $k > 2$, the above construction allows us to  impose a  
$C^{k-1}$ structure on ${\mathcal N}$  
in such a manner that $\Psi$ becomes a $C^{k-1}$ map. 
While this will actually turn out to be  technically useful, 
the reader should be warned, however, that such a 
$C^{k-1}$ structure is in no sense be
{\em natural} or {\em canonical},  because 
it depends on the $(k-1)$-jet of 
 our identification of the tubular neighborhood $\hat{\mathcal V}$ 
with $L\to \PP TN$, and such a choice 
  is  uniquely specified  by the geometry only when $k=2$;
for this reason, we will refer to such  a choice  as a {\em provisional}  $C^{k-1}$ structure. 
Fortunately, however, this apparent shortcoming will 
soon  be remedied. 
Indeed, the thrust of our argument is that  that $[\nabla ]$ induces  a 
certain complex structure $J$ on $\mathcal N$, 
and so endows 
$\mathcal N$  with a canonical $C^{\infty}$  structure.
In order to see this, we will proceed by first constructing a  certain 
 involutive complex distribution ${\bf D}$ on $\PP T_\CC M$,
and then analyzing its image under $\Psi$.

Since ${\mathcal Z}= \PP T_\CC M$, we have a bundle projection,
which we will denote by 
 $\hat{\mu}: {\mathcal Z}\to M$. The sub-bundle ${\bf V} = \ker \hat{\mu}_*\subset T{\mathcal Z}$
will be called the {\em vertical sub-bundle}.  Now
choose a connection 
 $\nabla$ representing
the given projective structure $[\nabla ]$,
and let ${\bf H}\subset T{\mathcal Z}$ be the {\em horizontal sub-bundle},
corresponding to parallel transport with respect to $\nabla$,
so that we have a 
 direct-sum decomposition
$$
T{\mathcal Z} = {\bf V}\oplus {\bf H}.
$$
Complexifying these bundles, we thus have
$$T_\CC {\mathcal Z} = {\bf V}_\CC\oplus {\bf H}_\CC ,$$
where $T_\CC  {\mathcal Z} = \CC \otimes T{\mathcal Z}$, etc. 
Notice that the derivative of the projection also gives us a 
 canonical isomorphism
$$\hat{\mu}_*: {\bf H}_\CC \stackrel{\cong}{\longrightarrow} \hat{\mu}^*T_\CC M.$$

Using this picture, we will now  define two 
line sub-bundles 
$${\bf L}_j \subset T_\CC {\mathcal Z} = \CC \otimes T{\mathcal Z}, ~~
j=1,2.$$ To this end, let us  
first  recall that each fiber of ${\mathcal Z}\to M$ 
is a $\CP_1$, so that we  have a fiber-wise complex structure tensor
$$J^\parallel : {\bf V}\to {\bf V}, ~~ (J^\parallel)^2= - {\bf 1},$$
and we  define ${\bf L}_1\subset {\bf V}_\CC$ to be the $(-i)$-eigenspace of 
$J^\parallel$:
$${\bf L}_1 = {\bf V}^{0,1}_{J^\parallel}.$$
On the other hand, each element of ${\mathcal Z} = \PP T_\CC M$
may be identified with a $1$-dimensional complex-linear subspace 
of $T_\CC M$, and 
 this picture gives us  
a tautological line sub-bundle  ${\bf L}_2$ of ${\bf H}_\CC\cong \hat{\mu}^* T_\CC M$:
\begin{equation}\label{dumbo}
\left. {\bf L}_2\right|_{[{\bf w}]} = (\hat{\mu}_{*[{\bf w}]})^{-1}( \mbox{span } {\bf w}).
\end{equation}
Set \begin{equation}
\label{recipe}
{\bf D}= {\bf L}_1\oplus {\bf L}_2 \subset T_\CC {\mathcal Z} .
\end{equation}
Then ${\bf D}$ is a $C^2$ distribution of complex $2$-planes on ${\mathcal Z}$.
We will now see that ${\bf D}$ is involutive, in the sense that 
$$[ C^1 ({\bf D}) , C^1 ({\bf D})] \subset C^{0} ({\bf D}).$$
Moreover,  ${\bf D}$ will turn out to be  unchanged if
we replace $\nabla$ with a projectively equivalent 
connection $\hat{\nabla}$. 

Indeed, let  $(x^1, x^2): {\Omega} \to \RR^2$ be a local coordinate system 
on ${\Omega}\subset M$, and let
$$
\Gamma^j_{k\ell }= \left\langle ~dx^j ~,~ \nabla_{\frac{\partial}{\partial x^k}}{\frac{\partial}{\partial x^\ell}}~\right\rangle
$$
be the corresponding Christoffel symbols of the connection 
$\nabla$. We can then introduce local coordinates
$(x^1, x^2, \zeta ) : \hat{\mu}^{-1} ({\Omega}) \to \RR^2 \times \CC$
on $\hat{\mu}^{-1} ({\Omega})\subset {\mathcal Z}$ by 
$$
\left[
\left. 
\left(
\frac{\partial}{\partial x^1} + \zeta \frac{\partial}{\partial x^2}
\right)\right|_{(x^1 , x^2)}
\right]   \longleftrightarrow (x^1,x^2,\zeta ) .$$ 
Then, in these coordinates, 
${\bf L}_1$ is spanned by $\partial/\partial \bar{\zeta}$, whereas 
${\bf L}_2$ is spanned  by 
$$
\Xi_0 = \frac{\partial}{\partial x^1}+ \zeta \frac{\partial}{\partial x^2}+ 
Q(x,\zeta , \zeta ) \frac{\partial}{\partial \zeta} + Q (x, \zeta , \overline{\zeta} )
 \frac{\partial}{\partial \overline{\zeta}}, 
$$
where 
$$Q(x, u, v) = 
 - \Gamma^2_{11} - \Gamma^2_{12} (u+v) -\Gamma^2_{22}uv 
+ \Gamma^1_{11} v+ \Gamma^1_{12}v(u+v)
+  \Gamma^1_{22}uv^2$$
encodes  the Christoffel symbols $\Gamma^j_{k\ell}$ of our
chart, which are of course functions of $x=(x^1,x^2)$. 
In particular, $\bf D$ is spanned by $\partial /\partial \overline{\zeta}$ and 
\begin{equation}
\Xi = \frac{\partial}{\partial x^1}+ \zeta \frac{\partial}{\partial x^2}+ 
P(x,\zeta ) \frac{\partial}{\partial \xi},\label{zuppa}
\end{equation}
where $\zeta = \xi + i \eta$ and where
\begin{equation}
P(x,\zeta)   = Q(x,\zeta, \zeta) = - \Gamma^2_{11}
+\left[\Gamma^1_{11}-2 \Gamma^2_{12}\right]\zeta
+ 
\left[ 2\Gamma^1_{12}-\Gamma^2_{22}\right]\zeta^2
+ \Gamma^1_{22}\zeta^3
\label{pesce}
\end{equation}
is evidently of the same differentiability class as $\nabla$. 
But 
$$\left[ ~\frac{\partial}{\partial \overline{\zeta}} ~,~\Xi~
\right ]
=
\left[ ~ \frac{\partial}{\partial \overline{\zeta}} ~ ,~ \frac{\partial}{\partial x^1}+ \zeta \frac{\partial}{\partial x^2}+ 
P(x,\zeta ) \frac{\partial}{\partial \xi} ~ \right ]=0 ~,$$ 
because 
$$\frac{\partial}{\partial \overline{\zeta}}\zeta= 0, \hspace{1cm}
 \frac{\partial}{\partial \overline{\zeta}}P(x,\zeta ) =0.$$
It therefore follows that 
${\bf D}= \mbox{span } \{ \Xi , \partial /\partial \bar{\zeta} \}$
is involutive, as claimed.

Notice that  the replacement 
$$\Gamma^i_{jk} \rightsquigarrow \Gamma^i_{jk} + \delta^i_j\beta_k + \beta_j \delta^i_k$$
leaves  $P(x,\zeta)$ unaltered. Thus replacing $\nabla$ with 
a projectively equivalent  connection $\hat{\nabla}$ leaves  $\Xi$ unchanged, and 
 ${\bf D}=\mbox{span} \{ \Xi , \partial/\partial \overline{\zeta}\}$ 
is   therefore 
projectively invariant.

The distribution ${\bf D}$ does not quite define a complex structure on ${\mathcal Z}$, 
because  certain real tangent vectors  are elements of 
${\bf D}$. Indeed, notice that, because ${\bf D}$ is the direct sum of ${\bf L}_1\subset {\bf V}_\CC$ and 
${\bf L}_2 \subset {\bf H}_\CC$, and because the projections $T_\CC M \to  {\bf V}_\CC$
and $T_\CC M \to  {\bf H}_\CC$ commute with complex conjugation, any  real
element of ${\bf D}$ must have real components in ${\bf L}_1$ and ${\bf L}_2$. 
But since ${\bf L}_1$ contains no non-zero real element, we therefore have
$$
{\bf D}\cap \overline{{\bf D}} = 
({\bf L}_1\cap \overline{\bf L}_1)+ ({\bf L}_2\cap \overline{\bf L}_2) = 
{\bf L}_2\cap \overline{\bf L}_2.$$
On the other hand, equation (\ref{dumbo}) tells us that 
${\bf L}_2$ contains a non-zero real element precisely at the
hypersurface $Z=\PP TM$ in ${\mathcal Z}=\PP T_\CC M$: 
\begin{equation}
\dim ({\bf D}_z\cap \overline{{\bf D}}_z ) = 
\left\{ \begin{array}{cl}
0, & z\not\in Z \\
1, & z\in Z .
\end{array}
\right.
\label{acs}
\end{equation}
Indeed,  ${\bf L}_2|_{Z}$ is simply 
the complexification 
$\CC \otimes \ker \nu_*$ of the  tangent space 
of the foliation $\mathcal F$ of $\PP TM$ by lifted geodesics. 
This  observation gives a somewhat more
geometric explanation for the previously noted projective 
invariance of ${\bf D}$. Indeed,  in equation 
(\ref{pesce})
 we   carefully chose our  complex vector field  $\Xi$  so that  at the locus $Z$, 
given by $\eta =0$, 
$\Xi$  is real and tangent to $\mathcal F$, with coefficients that 
are holomorphic in $\zeta=\xi + i \eta$, and so determined by the behavior of $\Xi$ along
$\eta =0$.

\begin{prop}\label{key} Let $[\nabla]$ be a Zoll projective structure which is 
represented by a $C^3$ connection $\nabla$ on $M\approx \RP^2$. Then 
there is a unique integrable almost-complex structure $J$ on $\mathcal N$ such that 
$$\Psi_* [{\bf D}] \subset T^{0,1}({\mathcal N}, J).$$
The unique $C^\infty$  structure on ${\cal N}$
associated with its maximal  atlas of  $J$-compatible complex charts
is compatible with the previously-constructed $C^1$ structure on 
$\mathcal N$, so that  $\Psi : {\mathcal Z} \to {\mathcal N}$  remains a $C^1$ map relative to 
this smooth structure; moreover, 
$\Psi$ actually becomes  $C^3$ on the open dense set ${\mathcal Z}-Z$.
Moreover, if $[\nabla ]$ is  represented by a 
$C^{k,\alpha}$ connection $\nabla$ on $M$,  $3\leq k\leq \infty$,
  $0< \alpha < 1$, 
and if ${\cal N}$ is again given the natural $C^\infty$  structure
associated with  $J$,
then  $\Psi : {\mathcal Z} \to {\mathcal N}$ is actually a $C^{k+1,\alpha}$ map
on  ${\mathcal Z}-Z$. 
\end{prop}

\begin{remark}
With the same hypotheses, 
we will later also show (remark, page \pageref{tricky}) that
$\Psi$ is actually   $C^{k+1,\alpha}$ on {\em all} of $\mathcal Z$.
\end{remark}

\begin{proof}
We begin by defining $J$ point-wise. 
On the open set ${\mathcal N} -N = \Psi ({\mathcal Z}-Z)$, we may do this by 
first observing that 
$$
T_\CC ({\mathcal N} -N) = \Psi_*{\bf D} \oplus \overline{\Psi_*{\bf D}} 
$$
by (\ref{acs}) and the fact that $\Psi|_{{\mathcal Z}-Z}$ is a diffeomorphism;
on ${\mathcal N}-N$, we now set 
$$J=
\left[
\begin{array}{cc}
-i& 0\\
0& +i
\end{array}\right]
$$
with respect to this direct sum decomposition. 
On the other hand, since ${\mathcal V}\subset {\mathcal N}$
is, by definition,  a copy of the total space of $TN\to N$, 
we have 
 a canonical identification 
$$T{\mathcal N}|_N = TN\oplus TN,$$ where the first factor is tangent to
$N$, and where the second factor is transverse to it; 
and along $N\subset {\mathcal N}$ we  can therefore  set 
$$J=
\left[
\begin{array}{cc}
0& -1\\
1& 0
\end{array}\right]
$$
with respect to this second direct sum decomposition. 
This defines the almost complex structure $J$ at all points of $\mathcal N$.

While it is not yet even yet clear that this $J$ is  
continuous, 
 it is at least  easy to see that $\Psi_*{\bf D} \subset T^{0,1}({\mathcal N}, J)$. Indeed, 
by construction, $\Psi_*{\bf D} = T^{0,1}({\mathcal N}, J)$ away from $N$. 
 On the other hand, $\Psi_*{\bf D} = \Psi_*{\bf V}^{0,1}$ along $Z$, and since 
we used $J^\parallel$ to pick out the normal factor of $T{\mathcal Z}|_{Z}= TZ \oplus L$
before blowing down, 
 $\Psi_*\circ J^\parallel= J \circ \Psi_*$ on ${\bf V}|_{Z}$,  and it
follows that $\Psi_*{\bf D} \subset  T^{0,1}({\mathcal N}, J)$ along $Z$, too. 
Moreover, $J$ is certainly the {\em only} almost-complex structure 
with this property, since, for any $y\in N$,   $$T_y{\mathcal N}= 
\Psi_* {\bf V}_x\oplus \Psi_*{\bf V}_{x'}$$
whenever $x\neq x'$ are 
 distinct  points  of the geodesic ${\mathfrak C}_y\subset M$ represented by 
$y$. 

Now since $[\nabla]$ has been assumed to be $C^3$, we can 
can give  ${\mathcal N}$  a `provisional' $C^2$ structure, compatible with 
its fixed $C^1$ structure,  relative to which 
  $\Psi$ 
becomes a $C^2$ map. We now claim that $J$ is  actually 
Lipschitz continuous in the associated charts on ${\mathcal N}$. 
Of course, this is is only a non-trivial statement near a point $y\in N$, since 
the restriction of 
$J$ to ${\mathcal N}-N$ corresponds, via $\Psi$, to a 
 $C^3$ almost-complex structure on ${\mathcal Z}-Z$.

Now  let us recall that  we have written down an 
explicit local framing $(\Xi, \partial/\partial \overline{\zeta})$ of  
${\bf D}$ such that $[\Xi, \partial/\partial \overline{\zeta}]=0$, and such that 
$\Xi$ is real along
$Z=\PP TM$, and spans the tangent space of the foliation $\mathcal F$
there. 
Giving an arbitrary  leaf $\hat{{\mathfrak C}}_y$ a parameter $t$ such that 
$\Xi = d/dt$ along the leaf, then, for any $C^2$ function $f$ on ${\mathcal N}$
 we  have
\begin{eqnarray*}
\frac{d}{dt}\left[ \Psi_* (\frac{\partial}{\partial \overline{\zeta}}) f\right]&=& 
\frac{d}{dt}\frac{\partial}{\partial \overline{\zeta}} \Psi^* f\\
&=&\Xi\frac{\partial}{\partial \overline{\zeta}}\Psi^* f\\
&=& \frac{\partial}{\partial \overline{\zeta}}\Xi\Psi^* f\\
&=&  \frac{\partial}{\partial \overline{\zeta}}\left[ \Psi_*(\Xi)f\right],
\end{eqnarray*}
so that, setting $\zeta = \xi + i\eta$,  
$$\frac{d}{dt}\left[ \Psi_* (\frac{\partial}{\partial \overline{\zeta}})\right] 
=  \frac{\partial}{\partial \overline{\zeta}}\left[ \Psi_*(\Xi)\right] =
  \frac{i}{2} \frac{\partial}{\partial \eta}
\left[ \Psi_*(\Xi)\right]
$$
at $y\in N$, 
since $\Psi_*(\Xi)\equiv 0$ along $Z$, where $\eta=0$.
Here the 
 right-hand side should   be  interpreted as the
  invariant  derivative 
{\em at a zero} of a section of a  vector bundle on 
 $\Sigma_x := \Psi [\hat{\mu}^{-1} (x)]\cong \CP_1$. 
On the other hand, 
$$\Psi_* \left(\frac{\partial}{\partial \overline{\zeta}}\right)\in T_y^{0,1}({\mathcal N},J)$$
for all $t$, by our previous discussion, so it follows that 
$$\left.
\frac{\partial}{\partial \eta}
\left[ \Psi_*(\Xi)\right]\right|_{\eta =0} \in T_y^{0,1}({\mathcal N},J),$$
too. Along $\Sigma_x$, we therefore have, 
near an arbitrary point $y\in N$,    two 
continuous sections of $T^{1,0}$
given by ${\bf e}_1=\Psi_*(\partial/\partial \overline{\zeta})$ and 
$$
{\bf e}_2 = \left\{
\begin{array}{cc}
\left[ \Psi_*(\Xi)\right]/\eta& \eta \neq 0\\
\frac{\partial}{\partial \eta}
\left[ \Psi_*(\Xi)\right]&\eta =0. 
\end{array}
\right.
$$
These sections are linearly independent at every point, and so 
span $T^{1,0}_y$, because $\det (\Psi_*)$ only vanishes to first
order at $Z$. Moreover, since $\Psi$ appears to be  $C^2$ in our coordinates, 
 these sections are both appear to be continuously differentiable  in our chart,
with derivatives that may be expressed in any coordinate
system in terms of partial derivatives of $\Psi$ of 
order $\leq 2$. Hence $J$ is also differentiable,  and in particular is Lipschitz, along
$\Sigma_x$, with Lipschitz constant controlled by the 
partial derivatives of $\Psi$ of order $\leq 2$.
Since the family $\{ \Sigma_x\}$ sweeps out all the 
radial lines  in our tubular neighborhood  $TN$  of $N\subset {\mathcal N}$,
if follows that the tensor field $J$ on ${\mathcal N}$ is  Lipschitz. 

Now recall that Rademacher's theorem asserts that all  the distributional 
first partial derivatives of a Lipschitz function are 
 locally  bounded measurable functions. 
The Nijenhuis tensor 
$$\tau ({\bf v},{\bf w}) = [{\bf v}, {\bf w}] - [J{\bf v}, J{\bf w}] +J[{\bf v}, J{\bf w}] +J[J{\bf v}, {\bf w}]$$
of our almost-complex structure $J$ is therefore well-defined in the distributional sense, 
and has components that are
locally  bounded and  measurable. 
On the other hand, since ${\bf D}$ is involutive and 
 $\Psi |_{{\mathcal Z}-Z}$ is a diffeomorphism, 
$\tau \equiv 0$ on a set of full measure, and therefore 
vanishes in the distributional sense. But Hill and Taylor  \cite{hiltay}
 have
recently shown that the Newlander-Nirenberg theorem 
holds for Lipschitz almost complex structures for which  $\tau =0$
in just this  distributional sense. Thus every point of $\mathcal N$
has a neighborhood on which we can find a pair $(z^1, z^2)$ of differentiable 
complex-valued functions with $dz^k\in \Lambda^{1,0}({\mathcal N}, J)$
and $dz^1\wedge dz^2 \neq 0$. Taking these to be the complex
coordinate systems gives ${\mathcal N}$ the structure of a 
compact complex surface. In particular, this gives ${\mathcal N}$ a
specific real-analytic structure, and hence  a specific $C^\infty$ structure. 

Finally,  we  address  
 the smoothness of  $\Psi : {\mathcal Z}\to {\mathcal N}$. 
Suppose that $\nabla$ is of differentiability class $C^{k,\alpha}$, and 
suppose that $f$ is a holomorphic function on some open subset
of ${\mathcal N}$; we then consider the function  $\Psi^*f$
on ${\mathcal Z}=\PP T_\CC M$. Now,  by \cite{hiltay},
$f$ is a $C^1$ function with respect to our (original, unchanged) $C^1$ structure on 
${\mathcal N}$, 
and $\Psi^* f$ is therefore a $C^1$ function, since $\Psi$ was
$C^1$ by construction.  Moreover, since $\Psi_* {\bf D} \subset T^{0,1}{\mathcal N}$, 
$\Psi^* f =0$ solves the Cauchy-Riemann equations 
$\bar{\partial}_{\bf D} (\Psi^* f)=0$  with respect to the $C^{k,\alpha}$ almost-complex
structure which ${\bf D}$ determines on  ${\mathcal Z}-Z$.
But since $\bar{\partial}_{\bf D}+ \bar{\partial}_{\bf D}^*$, defined with 
respect to an arbitrary $C^{k, \alpha}$ Hermitian metric on ${\mathcal Z}-Z$,
 is a
first-order elliptic system with $C^{k, \alpha}$ coefficients, 
elliptic regularity  \cite{morrey} tells us 
that $\Psi^* f$ is $C^{k+1,\alpha}$ on ${\mathcal Z}-Z$. 
Applying these observation when $f$ is any local complex coordinate 
$z^j$ on ${\mathcal N}$ then shows that $\Psi$
belongs to  the claimed differentiability class. 
\end{proof}

\begin{remark} \label{rocky}
The above proof uses a powerful recent analytic theorem  in order to obtain 
the result without too much hard work. Most readers will find it reassuring, however, that 
older technology may instead be used to prove a workable version of the proposition 
at  the price of a half-dozen derivatives and  a certain amount of careful calculation.
Moreover, this approach has the added benefit of providing some immediate
added information concerning the regularity of
$\Psi$ along $Z\subset {\mathcal Z}$.  
In particular, 
those primarily interested in the $C^\infty$ case  might
 well prefer the following  elementary argument.

Suppose that $\nabla$ is a $C^k$ connection, where $k=2\ell +2$. 
Choose $C^{2\ell+2}$ local real 
coordinates $(\check{y}^1,\check{y}^2)$ on $U\subset N$, and pull 
them back to $Z=\PP TM$ so as to obtain
 $C^{2\ell+2}$ functions $y^\jmath = \nu^*\check{y}^\jmath$
 on $\nu^{-1}U\subset Z$. By construction, these solve   the
equation $\Xi y^\jmath=0$.  We now extend the $y^\jmath$ as 
$C^{\ell+2}$ complex-valued  functions 
${\mathfrak z}^\jmath$ defined on an open set in 
${\mathcal Z}$ 
 by requiring  that $\partial y^\jmath/\partial
\bar{\zeta}$ vanish to order $\ell-1$ along $Z$.
This  completely specifies the $\ell$-jet of the 
function, and we must have 
$$
{\mathfrak z}^\jmath(x^1,x^2,\xi,\eta) =\sum_{r=0}^{\ell} \frac{i^r}{r!}
\eta^r \left.
 \frac{\partial^r y^\jmath}{\partial \xi^r}\right|_{(x^1,x^2,\xi)}  ~~+ O(\eta^{\ell +1}). 
$$
Indeed, this recipe does indeed give us  
\begin{eqnarray*}
 \frac{\partial {\mathfrak z}^\jmath}{\partial \overline{\zeta}}
&=&\frac{1}{2}
 \left( \frac{\partial}{\partial \xi}
+ i\frac{\partial}{\partial \eta}\right)  \left( \sum_{r=0}^{\ell} \frac{i^r}{r!} \eta^r \frac{\partial^r y^\jmath}{\partial \xi^r}
 ~~+ O(\eta^{\ell +1}) \right)
\\&=& 
\frac{1}{2}
\sum_{r=0}^{\ell} \frac{i^r}{r!} \eta^r\frac{\partial^{r+1} y^\jmath}{\partial \xi^{r+1}}
-\frac{1}{2}
\sum_{r=1}^{\ell} \frac{i^{r-1}}{(r-1)!} \eta^{r-1}\frac{\partial^r y^\jmath}{\partial \xi^r} ~~+ O(\eta^{\ell})
\\&=& 
\frac{1}{2}
\sum_{r=0}^{\ell} \frac{i^r}{r!} \eta^r\frac{\partial^{r+1} y^\jmath}{\partial \xi^{r+1}}
-\frac{1}{2}
\sum_{r=0}^{\ell-1} \frac{i^{r}}{r!} \eta^{r}\frac{\partial^{r+1} y^\jmath}{\partial \xi^{r+1}} ~~+ O(\eta^{\ell})
\\&=& 
\frac{1}{2}
\frac{i^{\ell}}{\ell!} \eta^{\ell}\frac{\partial^{\ell+1} y^\jmath}{\partial \xi^{\ell+1}} ~~+ O(\eta^{\ell})
\\&=&   O(\eta^{\ell}) , 
\end{eqnarray*}
and  since the cancellation is a  term-by-term matter, 
uniqueness of the $\ell$-jet follows. But since our condition on the 
$\ell$-jet is obviously independent of the choice of 
coordinates 
$(x^1,x^2)$ on $M$, global 
existence now follows by patching together any such
local choices via a partition of unity.

The uniqueness argument also has another useful consequence. Notice that there 
certainly are $C^2$ coordinates $(\tilde{\mathfrak z}^1, \tilde{\mathfrak z}^2)$ 
for our provisional $C^2$ structure on
${\mathcal N}$ whose restrictions to $N$ are the $\check{y}^\jmath$, and 
which are satisfy $\overline{\partial}_J \tilde{\mathfrak z}^\jmath =0$
to $0^{th}$ order along $N$, since the restriction of $J$ to 
$T{\mathcal N}|_N$ is $C^{k-1}$. 
But pulling these back to ${\mathcal Z}$ would gives us $C^2$ functions 
killed by  $\partial/\partial\overline{\zeta}$ to $0^{th}$ order along $Z$, and the
$\ell=1$ version of the above calculation therefore gives 
$$
\Psi^*\tilde{\mathfrak z}^\jmath =  {\mathfrak z}^\jmath + O(\eta^2).
$$
It follows that $({\mathfrak z}^1, {\mathfrak z}^2)$ is
actually a $C^1$ complex-valued coordinate system on 
${\mathcal N}$. Our strategy will now be to analyze the 
the almost-complex structure $J$ 
by thinking of $(x^1,x^2,\xi,\eta)\mapsto ({\mathfrak z}^1, {\mathfrak z}^2)$
as a representation of 
$\Psi$ in special coordinates

To this end, we next
 observe that, since  $[\Xi , \frac{\partial}{\partial \overline{\zeta}} ]=0$, the 
$C^{\ell+1}$ function $\Xi {\mathfrak z}^\jmath$ satisfies 
$$
\frac{\partial^m }{\partial \overline{\zeta}^m} \Xi {\mathfrak z}^\jmath = \Xi \frac{\partial^m }{\partial \overline{\zeta}^m}{\mathfrak z}^\jmath =
\Xi O(\eta^{\ell-m+1}) = O(\eta^{\ell-m+1}),
$$
so that  
$$\left.
\left(\frac{\partial}{\partial {\xi}} +i\frac{\partial}{\partial {\eta}}
\right)^m (\Xi {\mathfrak z}^\jmath)\right|_{\eta=0} \equiv 0,
$$
for $m=0, \ldots, \ell$.
But since $\Xi {\mathfrak z}^\jmath \equiv 0$ along $\eta =0$, this tells us that 
$$\left.
\frac{\partial^m }{\partial {\eta}^m} (\Xi {\mathfrak z}^\jmath)\right|_{\eta=0} \equiv 0
$$
for $m = 0 , \ldots , \ell$, and hence  that
$$
\Xi {\mathfrak z}^\jmath = O (\eta^{\ell+1}).
$$

Now we 
have already shown, by an 
elementary argument, that the almost-complex structure $J$ is characterized,
in a point-wise manner, 
by the fact that $\Psi_*\partial/\partial \overline{\zeta}$ and $\Psi_* \Xi$ are 
always elements  of $T^{0,1}({\mathcal N},J)$. 
Since $\mbox{span} \{ \partial/\partial {\mathfrak z}^1 , \partial/\partial {\mathfrak z}^2\}$ contains 
the image of $\partial/\partial \overline{\zeta}$ and (trivially) $\Xi$ along 
the locus $N$ given by $\Im  m {\mathfrak z}^\jmath =0$, we must therefore have 
$$T^{0,1}({\mathcal N},J)|_{\Im m {\mathfrak z}^\jmath =0}= \mbox{span}\left\{ \frac{\partial}{\partial \overline{{\mathfrak z}}^\jmath} 
\right\}_{\jmath=1,2} ~~,$$
and 
$$
T^{*1,0}({\mathcal N},J)|_{\Im m {\mathfrak z}^\jmath =0}= \mbox{span}\left\{ d{\mathfrak z}^\jmath\right\}_{\jmath=1,2}~~.$$
Elsewhere, 
$$T^{*1,0}({\mathcal N},J)= 
\mbox{span}\left\{ d{\mathfrak z}^\jmath -\sum_\imath a^\jmath_\imath d\overline{{\mathfrak z}}^\imath
\right\}_{\jmath=1,2}$$
and
$$T^{0,1}({\mathcal N},J)= \mbox{span} \left\{ \frac{\partial}{\partial \overline{{\mathfrak z}}^\jmath}
+ \sum_\imath a^\imath_\jmath \frac{\partial}{\partial {{\mathfrak z}}^\imath}
\right\}_{\jmath=1,2} ~~,$$
where the $a^\jmath_\imath$ are to be found by solving the equation 
$$
\left[
\begin{array}{cc}
a_1^1&a^1_2
\\a^2_1&a^2_2
\end{array}
\right] 
\left[
\begin{array}{cc}
\Xi \overline{{\mathfrak z}}^1& \frac{\partial \overline{{\mathfrak z}}^1}{\partial \overline{\zeta}}
\\
\Xi \overline{{\mathfrak z}}^2& \frac{\partial \overline{{\mathfrak z}}^2}{\partial \overline{\zeta}}
\end{array}
\right]
=
\left[
\begin{array}{cc}
\Xi {{\mathfrak z}}^1& \frac{\partial {{\mathfrak z}}^1}{\partial \overline{\zeta}}
\\
\Xi {{\mathfrak z}}^2& \frac{\partial {{\mathfrak z}}^2}{\partial \overline{\zeta}}
\end{array}
\right].
$$
But 
$$\frac{\partial}{\partial \overline{\zeta}} \overline{{\mathfrak z}}^\jmath = \frac{\partial y^\jmath}{\partial \xi}  + O (\eta),$$
and 
$$
\Xi \overline{{\mathfrak z}}^\jmath = \Xi \left ( -i\eta \frac{\partial y^\jmath}{\partial \overline{\xi}} + O(\eta^2)\right)
= -i \eta [\Xi ,  \frac{\partial }{\partial \xi} ] y^\jmath
= i\eta  \frac{\partial y^\jmath}{\partial x^2} + i\eta P^\prime (\xi) \frac{\partial y^\jmath}{\partial {\xi}} + O (\eta^2),
$$
so that 
$$
 \left|
\begin{array}{cc}
\Xi \overline{{\mathfrak z}}^1& \frac{\partial \overline{{\mathfrak z}}^1}{\partial \overline{\zeta}}
\\
\Xi \overline{{\mathfrak z}}^2& \frac{\partial \overline{{\mathfrak z}}^2}{\partial \overline{\zeta}}
\end{array}
\right|
= i\eta 
 \left|
\begin{array}{cc}
\frac{\partial y^1}{\partial x^2} + P^\prime (\xi) \frac{\partial y^1}{\partial {\xi}} 
& \frac{\partial {y}^1}{\partial \xi}
\\
\frac{\partial y^2}{\partial x^2} + P^\prime (\xi) \frac{\partial y^2}{\partial {\xi}} 
& \frac{\partial {y}^2}{\partial \xi}
\end{array}
\right| + O (\eta^2) = i\eta  \frac{\partial(y^1,y^2)}{\partial(x^2,\xi)}
  + O(\eta^2).
$$
But ${\partial(y^1,y^2)}/{\partial(x^2,\xi)}\neq 0$ everywhere, since $\Xi$ is always 
linearly independent from $\partial/\partial x^2$ and $\partial/\partial\xi$. Thus 
\begin{eqnarray*}
\left[
\begin{array}{cc}
a_1^1&a^1_2
\\a^2_1&a^2_2
\end{array}
\right] &=& 
\left[
\begin{array}{cc}
\Xi {{\mathfrak z}}^1& \frac{\partial {{\mathfrak z}}^1}{\partial \overline{\zeta}}
\\
\Xi {{\mathfrak z}}^2& \frac{\partial {{\mathfrak z}}^2}{\partial \overline{\zeta}}
\end{array}
\right]
\left[
\begin{array}{cc}
\Xi \overline{{\mathfrak z}}^1& \frac{\partial \overline{{\mathfrak z}}^1}{\partial \overline{\zeta}}
\\
\Xi \overline{{\mathfrak z}}^2& \frac{\partial \overline{{\mathfrak z}}^2}{\partial \overline{\zeta}}
\end{array}
\right]^{-1}
\\ &=&
\left[
\begin{array}{cc}
O(\eta^{\ell+1})& O(\eta^{\ell})\\
O(\eta^{\ell+1})&O(\eta^{\ell})
\end{array}
\right]
\frac{1}{i\eta} \left(
 \frac{\partial (x^2,\xi)}{\partial (y^1,y^2)} + O (\eta)
\right) 
\left[
\begin{array}{cc}
 \frac{\partial \overline{{\mathfrak z}}^2}{\partial \overline{\zeta}}& -\frac{\partial \overline{{\mathfrak z}}^1}{\partial \overline{\zeta}}
\\
-\Xi \overline{{\mathfrak z}}^2&\Xi \overline{{\mathfrak z}}^1
\end{array}
\right]
\\ &=&
\left(
 \frac{\partial (\xi,x^2)}{\partial (y^1,y^2)} + O (\eta)
\right) 
\left[
\begin{array}{cc}
O(\eta^{\ell})& O(\eta^{\ell-1})\\
O(\eta^{\ell})& O(\eta^{\ell-1})
\end{array}
\right]
\left[
\begin{array}{cc}
O(\eta^0)& O(\eta^0)
\\
O(\eta)& O(\eta)
\end{array}
\right]
\\
&=& O(\eta^{\ell-1}) 
\end{eqnarray*}
More precisely, for $(x^1,x^2,\xi, \eta)$ in any fixed compact set, there is a constant 
$C$ such that 
$$|a^\jmath_\imath | < C~ |~\eta ~|^{\ell -1}.$$
For the corresponding set in ${\mathcal N}$, this becomes
the statement that 
$$|a^\jmath_\imath | < C_1  ~|~\Im m ~\vec{\mathfrak{z}} ~|^{\ell -1}.$$
But since $\Psi$ is a proper map, it only takes a finite 
number of closed coordinate balls to cover the inverse image of 
any compact set in $\mathcal N$, and hence we have 
$$|a^\jmath_\imath | < C_2 ~|~\Im m ~\vec{\mathfrak{z}} ~|^{\ell -1}$$
as long as $\vec{\mathfrak z}=({\mathfrak z}^1, {\mathfrak z}^2)$ 
is constrained to lie in any fixed compact set.

Since $(x^1,x^2,\xi,\eta)\mapsto ({\mathfrak z}^1,{\mathfrak z}^2)$ is 
a $C^{\ell+1}$ diffeomorphism away from $\eta =0$, the 
$a^\jmath_\imath$ are $C^{\ell+1}$ functions of the 
$({\mathfrak z}^1,{\mathfrak z}^2)$ away from $\Im m ~{\mathfrak z}^1 = \Im m ~{\mathfrak z}^2
=0$, and on the other hand we have seen that they vanish to
order $\ell -2$ along this bad locus. Thus the 
$a^\jmath_\imath$ are $C^{\ell-2}$ functions of the 
${\mathfrak z}^\jmath$, and the  complex structure
$J$ on ${\mathcal N}$  is $C^{\ell-2}$
in these coordinates. If $\ell -2 \geq 1$, the Nijenhuis tensor 
therefore vanishes identically by continuity, since  it is already known to
vanish on an open dense set. 
If $\ell-2\geq 4$, or in other words if $[\nabla ]$  is at least $C^{14}$,
we may therefore apply the original Newlander-Nirenberg theorem
\cite{newnir} to get $C^{\ell-2}$ functions
$(z^1,z^2)$ of $({\mathfrak z}^1, {\mathfrak z}^2)$ which are holomorphic
with respect to $J$. The Malgrange refinement \cite{malgrange} of Newlander-Nirenberg
may similarly be applied if $\ell-2\geq 2$, 
or in other words if $[\nabla ]$ is at least $C^{10}$. 
The rest of the proof  then proceeds as before. Notice, however, 
that this second argument also  directly verifies  that  $\Psi : {\mathcal Z}\to {\mathcal N}$
is at least $C^{[k/2]-3}$ along $Z\subset {\mathcal Z}$.
\label{road}
\end{remark}

Having constructed our compact complex surface ${\mathcal N}$,
we will now try to unmask its identity. To this end, 
recall that we originally assembled ${\mathcal N}$ from two open sets, 
${\mathcal U}={\mathcal Z}-Z$ and ${\mathcal V} \approx T\RP^2$. 
However,  ${\mathcal U}$
may be identified with the space of all almost-complex structures\footnote{Indeed, 
the fact that $\bf D$ is a complex structure on
$\mathcal U$  thus naturally arises in the context of the 
O'Brian-Rawnsley 
generalization \cite{obiwon} of the    Atiyah-Hitchin-Singer 
approach \cite{AHS} to  twistor theory.} on 
$M$, since an almost-complex structure is completely characterized 
by its $(0,1)$-tangent space, and in dimension $2$ this may be taken
to be any $1$-dimensional subspace of $T_\CC M$ which is not spanned
by a real vector. 
 Thus ${\mathcal U}\to M$
may be identified with the space  of pairs $([h], \circlearrowleft)$, 
where $h$ is a Riemannian metric on some tangent space $T_xM$, 
$[h]$ is its conformal class, and $\circlearrowleft$ denotes a choice of
orientation of $T_xM$. Since  the space of Riemannian metrics
is a convex cone,  ${\mathcal U}$ therefore canonically 
 deform retracts to the set of point-wise orientations
on  $M$, once we choose a
single `background' Riemannian metric $h_0$ on $M$.
But the $2$-fold cover $\tilde{M}$ of $M$ by its set of local 
orientations is evidently just $S^2$, since $M=\RP^2$
by assumption. This shows that ${\mathcal U}$ is homotopy equivalent to
$S^2$. 

With this observation in hand, we are now in a position to 
list   some identifying traits of our complex surface $({\mathcal N}, J)$. 
 
 \begin{prop} \label{complexn} 
Let $[\nabla ]$ be a Zoll projective structure on $M=\RP^2$,
and let $N\approx \RP^2$ denote the corresponding space of
unoriented geodesics. Then there is a compact complex surface
$\mathcal N$ and an embedding $N\hookrightarrow {\mathcal N}$
such that 
\begin{itemize}
\item $\pi_1 ({\mathcal N}) =0$; 
\item there is an anti-holomorphic involution $\sigma : {\mathcal N}
\to {\mathcal N}$ with fixed-point set N; 
\item for all $x \in M$, there is a  
complex curve $\Sigma_x \subset 
{\mathcal N}$, 
$\Sigma_x \cong \CP_1$, such that $$\ell_x = \Sigma_x\cap  N;$$
\item the $\Sigma_x$ all represent the same element of $\pi_2 ({\mathcal N})$; and
\item if $x$ and $x'$ are distinct points of $M$, then 
$\Sigma_x$ and $\Sigma_{x'}$ are transverse, and meet in exactly one point. 
\end{itemize}
\end{prop}

\begin{proof}
By construction, 
${\mathcal N}= {\mathcal U}\cup {\mathcal V}$, where ${\mathcal U}={\mathcal Z} - Z$
and ${\mathcal V} =TN \approx T\RP^2$. But we have just seen that 
 ${\mathcal U}$ deform retracts to $S^2$. Moreover, ${\mathcal V}$ deform 
retracts to $N\approx \RP^2$,  and the inclusion map 
$\jmath : {\mathcal U}\cap {\mathcal V} \hookrightarrow  {\mathcal V}$
is homotopic to the bundle projection $\wp : (TN-0_{N}) \to N$. 
Because ${\mathcal U}$ is simply connected and ${\mathcal U}\cap {\mathcal V}$ is connected, 
the Seifert-van Kampen theorem 
tells us that 
$$\pi_1 ({\mathcal N} ) = \frac{\pi_1 ({\mathcal V})}{ \jmath_{\natural} [\pi_1 ({\mathcal U} \cap {\mathcal V})]} = 
\frac{\pi_1 (N ) }{\wp_{\natural}
 [\pi_1 (TN-0_N)]}.$$
But $\wp_{\natural} : \pi_1 (TN-0_N) \to \pi_1 (N)$ is surjective, since the 
fibers of $\wp$ are path connected. Hence ${\mathcal N}$
is simply connected. 

Complex conjugation $\PP T_\CC M \to \PP T_\CC M$
sends the distribution ${\bf D}$ to its conjugate $\overline{{\bf D}}$. The 
induced involution $\sigma : {\mathcal N}\to {\mathcal N}$ is therefore 
anti-holomorphic, and obviously has fixed point set precisely consisting of $N$. 

 For each $x\in M$, set $\Sigma_x = \wp (\PP T_{x\CC} M)$. Then 
$\Sigma_x$ is an embedded genus $0$ complex curve in ${\mathcal N}$. 
Since  the fibers of $\PP T_\CC M$ are all homotopic, 
so are their images in $\mathcal N$. Moreover, since the fibers of 
$\PP T_\CC M$ are all disjoint, we must have $\Sigma_x \cap \Sigma_{x'} \subset N$.
But, by construction, $\Sigma_x \cap N = \ell_x$, and so 
$$\Sigma_x \cap \Sigma_{x'} = (\Sigma_x \cap N)  \cap (\Sigma_{x'} \cap N) = \ell_x
\cap \ell_{x'}, $$
and if $x\neq x'$ this 
 consists of precisely one point $y$, representing  the unique geodesic joining 
$x$ to $x'$; cf.  Corollary \ref{lcapl}. Now $\Sigma_x$ and $\Sigma_{x'}$ are both 
$\sigma$-invariant, so $T_y\Sigma_x \cap T_y \Sigma_{x'}$ is invariant under
the complex anti-linear involution 
$\sigma_*$ of $T_y{\mathcal N}$, which we may identify with 
complex conjugation  on $\CC\otimes T_yN$. 
But since $T_y\ell_x\cap T_y\ell_{x'}=0$, 
its complexification $T_y\Sigma_x\cap T_y\Sigma_{x'}$
is also zero, and $\Sigma_x$ and $\Sigma_{x'}$ therefore intersect
transversely, at the unique point $y$, exactly as claimed. 
\end{proof}

We now come to the key step in our proof, which is to 
observe that $\mathcal N$ must be biholomorphic to $\CP_2$. 
It is a deep and remarkable  fact \cite{yau} that, 
up to biholomorphism, 
 $\CP_2$ is the only simply connected complex surface of Euler characteristic $3$,
and it might therefore be tempting to now  invoke this powerful result, much as we will 
later do in  \S \ref{zoe} below.  However, we will 
 actually need to know a great deal about the  biholomorphism 
$F: {\mathcal N}\to \CP_2$, and for this reason it is 
in every sense more satisfactory   to instead make use of 
the following low-tech lemma, based on the classical ideas of Castelnuovo,
Enriques and 
Kodaira;  cf. 
  \cite[Proposition V.4.3]{bpv}.
As a courtesy to the reader, as well as  to emphasize the 
elementary nature   of  the result, we  include a short, complete proof.

\begin{lem}\label{castle}
Let ${\mathcal S}$ be a simply connected compact complex surface, 
equipped with 
a fixed homology class
${\bf a} \in  H_2 ({\mathcal S}, \Z)$ such that  ${\bf a}\cdot {\bf a} =1$. 
For every 
$p\in {\mathcal S}$,  suppose that there exists a non-singular, embedded
 complex curve $\Sigma\subset {\mathcal S}$ of  genus $0$
passing through  $p$, with homology class $[\Sigma]={\bf a}$. Then
${\mathcal S}$ is biholomorphic to $\CP_2$,  in such  a manner  that 
all of the given curves  become projective lines. 
\end{lem} 

\begin{proof}
Since  the 
Fr\"olicher spectral sequence
of any complex surface degenerates at the $E_1$ level \cite[Theorem IV.2.7]{bpv}, 
we have 
$$H^1({\mathcal S},\CC) \cong H^1({\mathcal S},{\mathcal O}) \oplus H^0({\mathcal S},\Omega^1),$$
so  the assumption that $\pi_1 ({\mathcal S})=0$ 
immediately implies that $H^1({\mathcal S}, {\mathcal O}) =0$. 
But  the divisor line bundle  ${\mathcal O}(\Sigma)$ 
of  any of the curves $\Sigma\subset {\mathcal S}$ fits into an exact sequence 
\begin{equation}
\label{rest}
0\to {\mathcal O} \stackrel{{f}\cdot}{\to}  {\mathcal O}(\Sigma) \to n_{\Sigma} \to 0
\end{equation}
of sheaves on ${\mathcal S}$, where $n_{\Sigma}$ 
is the normal sheaf of $\Sigma$, extended 
to ${\mathcal S}$ by $0$,
and where ${f}\cdot$ denotes multiplication by a holomorphic section 
${f}$ of ${\mathcal O}(\Sigma)$ which vanishes only at $\Sigma$, 
with $d{f}\neq 0$ along $\Sigma$. 
 Now the normal bundle of $\Sigma$ has degree 
${\bf a}\cdot {\bf a} =1$,  and thus  $n_{\Sigma}$ can 
be identified with the unique degree-$1$ holomorphic line bundle ${\mathcal O}(1)$ on $\CP_1$. 
Since $H^1({\mathcal S}, {\mathcal O}) =0$, the 
long exact sequence in cohomology induced by (\ref{rest}) 
therefore gives us the short exact sequence 
\begin{equation}\label{shorty}
0\to \CC \stackrel{{f}\cdot}{\to} 
 \Gamma ({\mathcal S}, {\mathcal O}(\Sigma)){\to} \Gamma (\CP_1 , {\mathcal O}(1))\to 0. 
\end{equation} 
In particular,  $H^0( {\mathcal S}, {\mathcal O}(\Sigma)) \cong \CC^3$; moreover,  
 there is a holomorphic section of ${\mathcal O}(\Sigma)$ which is 
non-zero at any given point of ${\mathcal S}$. The  associated  map 
$$F: {\mathcal S} \to {\mathbb P}[H^0( {\mathcal S}, {\mathcal O}(\Sigma))^*]\cong \CP_2,$$
is  thus everywhere defined. Also notice that 
  $F(\Sigma)$ is a 
projective line  ${\mathcal P}  \subset \CP_2$,  and that  
the derivative of $F$ is  of maximal rank at any point $p$ of 
$\Sigma$, since (\ref{shorty}) 
 allows us to produce two  sections of ${\mathcal O}(\Sigma)$, $f$ and another one, 
 which vanish at $p$, but have linearly 
independent
derivatives there.

Since 
$H^1({\mathcal S}, \O)=0$, 
the exact sequence 
$$\cdots \to H^1({\mathcal S}, {\mathcal O}) \to H^1 ({\mathcal S}, {\mathcal O}^* ) \stackrel{c_1}{\to} H^2({\mathcal S}, \Z ) \to \cdots  $$
 tells us that holomorphic line bundles on ${\mathcal S}$ are classified by their 
first Chern classes.
But if 
$\Sigma$ and $\Sigma'$ are two complex curves in the homology class ${\bf a}$,
their  divisor line bundles ${\mathcal O}(\Sigma)$ and ${\mathcal O}(\Sigma')$ both have 
 Chern class equal to the Poincar\'e dual of ${\bf a}$.  
Thus ${\mathcal O}(\Sigma)\cong {\mathcal O}(\Sigma')$, and 
$\Gamma ({\mathcal S}, {\mathcal O}(\Sigma))=  \Gamma ({\mathcal S}, {\mathcal O}(\Sigma'))$. The 
 holomorphic map 
$F: {\mathcal S}\to \CP_2$ determined by $\Sigma$ therefore also maps $\Sigma'$ biholomorphically
to a projective line ${\mathcal P}'$, and the derivative of $F$ has maximal rank at every
point of $\Sigma'$. 
Since, by hypothesis,  we may find such a  curve through 
any point, $F$ is a local biholomorphism. But  since ${\mathcal S}$ is compact,
$F$ is therefore a covering map; and 
since $\CP_2$ is simply connected, we conclude  that 
 $F$ is a biholomorphism. \end{proof}

\begin{thm} \label{rigid}
Let $(M,[\nabla])$ be a compact 2-manifold with  Zoll projective structure 
of odd conjugacy number. Assume that 
 $\nabla$ is  of differentiability class $C^{k,\alpha}$,  for some $k\geq 3$,
and some 
$\alpha \in (0,1)$.
Then there is a $C^{k+2, \alpha}$  diffeomorphism 
$\Phi : M\stackrel{\approx}{\longrightarrow} \RP^2$ such that   
$[\nabla]=  [ \Phi^* \triangledown]$,
where $\triangledown$ is 
  the Levi-Civita connection 
$\triangledown$ of the standard, constant curvature 
Riemannian metric 
$g$ on $\RP^2$. 
\end{thm}
\begin{proof}
By  Proposition \ref{complexn}, the entire complex surface ${\mathcal N}$ is swept out 
by the genus zero curves $\Sigma_x$, $x\in M$, and the 
homology class 
$[\Sigma_x ]\in H_2({\mathcal N}, \ZZ)$
is independent of $x$. Moreover, this homology class has self-intersection
$$[\Sigma_x ]\cdot [\Sigma_x ] = [\Sigma_x ]\cdot [\Sigma_{x'} ] = 
1 ,$$
since $\Sigma_x$ and $\Sigma_{x'}$ intersect transversely  in one point
whenever $x\neq x'$. 
 Lemma \ref{castle} therefore tells us that 
there is a biholomorphism $F: {\mathcal N}\to \CP_2$ which sends each 
of the 
 complex curves $\Sigma_x$  to a corresponding projective line $\CP_1\subset \CP_2$.

Now the anti-holomorphic involution $\sigma : {\mathcal N} \to {\mathcal N}$
induces an anti-holomorphic involution $\tilde{\sigma}= F\circ\sigma\circ F^{-1}:
\CP_2\to \CP_2$. By taking the 
 Jacobian determinant of this map, we then obtain 
 to an anti-holomorphic involution $\tilde{\sigma}^* : K \to K$ 
of the canonical line bundle 
$K=\Lambda ^{2,0}$ of $\CP_2$.  But $K$ has a unique holomorphic 
cube-root $K^{1/3}$, the frame bundle of which is the universal cover of 
the 
frame bundle of $K$; 
and covering space theory now  tells us that 
$\tilde{\sigma}^*$ has three possible anti-holomorphic lifts $\varrho : K^{1/3} \to 
K^{1/3}$, differing by  multiplicative factors of a cube-root of unity. Choose
any such lift, and observe that $\varrho^2$ is the identity on any fiber over 
the fixed-point locus $F (N)$ of $\tilde{\sigma}$; since $F(N)$ is totally real
and 
of maximal dimension, 
the principle of analytic 
continuation therefore implies that 
the holomorphic map 
$\varrho^2$ must therefore be the identity. The anti-linear map
$$\varrho^* : \Gamma (\CP_2 , {\mathcal O}(K^{-1/3}))\to  \Gamma (\CP_2 , {\mathcal O}(K^{-1/3}))
$$
therefore satisfies $(\varrho^*)^2={\bf 1}$. It is therefore diagonalizable over $\RR$,
with eigenvalues $\pm 1$, and, because it is anti-linear, it  can be put in the form 
$$(z_1, z_2 , z_3 ) \mapsto (\bar{z}_1, \bar{z}_2 , \bar{z}_3 )$$
by choosing a suitable basis for 
$\Gamma (\CP_2 , {\mathcal O}(K^{-1/3}))\cong \Gamma (\CP_2 , {\mathcal O}(1))\cong \CC^3$. 
But $[z_1 : z_2 : z_3 ]$  gives us a set of homogeneous coordinates on 
$\CP_2$, so we have succeeded in identifying  $\sigma : {\mathcal N} \to {\mathcal N}$ 
with the
standard 
complex conjugation on $\CP_2$. In the process, we have thereby 
identified $N$ with $\RP^2\subset \CP_2$, and each complex curves
$\Sigma_x$ with a complex projective line $\CP_1$ which is invariant
under complex conjugation.

Now let  $ \CP_2^*=\PP (\CC^{3*})$ denote  the dual projective
plane of $\CP_2= \PP (\CC^{3})$, and consider the map 
\begin{eqnarray*}
\Phi_0 : M &\to&  \CP_2^*\\
x &\mapsto& F(\Sigma_x)^\perp, 
\end{eqnarray*}
where $\perp$ denotes the usual correspondence between lines in 
$\CP^2$ and points in $\CP^{2*}$. 
We claim that $\Phi_0$ is of differentiability class $C^{k+2,\alpha}$.
Indeed, let ${\mathcal C}\subset {\mathcal U}$ be a (non-compact) holomorphic
curve which is transverse to the fibers of $\hat{\mu}$, obtained by setting
some local complex coordinate ${\mathfrak z}^1$ equal to zero. Since the almost
 complex structure on $\mathcal U$ is of class $C^{k,\alpha}$, elliptic regularity 
tells us that the 
local complex coordinates $({\mathfrak z}^1,{\mathfrak z}^2)$ are of class $C^{k+1,\alpha}$,
and   $\mathcal C$ is therefore   representable as the image of a 
$C^{k+1,\alpha}$ map  from  an open set in $\CC$ to $\mathcal U$.
But the projection from  $\mathcal C$ to $M$ is a local diffeomorphism,
and so $\mathcal C$ may locally be thought of as the graph of a $C^{k+1,\alpha}$
local section $\varsigma$ of ${\mathcal U}\to M$. But such a section 
is precisely a local almost-complex structure on $M$ of 
differentiablity class $C^{k+1,\alpha}$. Since the map 
$F\circ \Psi \circ \varsigma$ is holomorphic with respect to this $C^{k+1,\alpha}$
almost-complex structure, and it is  therefore of class $C^{k+2,\alpha}$
by elliptic regularity. But on the domain of this function, 
$F(\Sigma_x)^\perp$ is the unique line joining $F(\Psi ({\varsigma}(x))$
to its complex conjugate, and so can be expressed in homogeneous
coordinates as $$\Phi_0 (x) = F(\Psi ({\varsigma}(x))\times \overline{F(\Psi ({\varsigma}(x))},$$
where $\times : \CC^3 \times \CC^3 \to \CC^{3*}$ is the 
vector cross-product. Since $M$ is covered by the domains of 
such local almost-complex structures $\varsigma$, this shows that 
 $\Phi_0$ is $C^{k+2,\alpha}$
on all of $M$.

Now notice that 
$\Phi_0$
is also an immersion, because $\Psi$ is a diffeomorphism on 
${\mathcal U}$, and the section of the normal bundle 
of $\Sigma_x\subset {\mathcal N}$ corresponding to a non-zero
element of $T_xM$ is therefore never identically zero. 
Moreover,  because each $F(\Sigma_x)$ is invariant under complex conjugation, 
  $\Phi_0 (M)$ actually lies
in the real dual projective plane $\RP^{2*}\subset \CP_2^*$. 
Thus, $\Phi_0$ actually gives us a $C^{k+2,\alpha}$  immersion  
$$
\Phi : M \to   \RP^{2*}$$
which can be described as 
$$
x \mapsto  F(\ell_x)^\perp .
$$
But since $M$ is a compact $2$-manifold, this immersion must be a covering map, 
and since $\pi_1(M)\cong \pi_1 (\RP^{2*})=\ZZ_2$, it follows that 
$\Phi$ is a diffeomorphism. Moreover,$\Phi$ sends the geodesic ${\mathfrak C}_y$
to the set of projective lines through the point $F(y)\in \RP^2$, or
in other words to the projective line $F(y)^\perp$ in $\RP^{2*}$.
This shows that $\Phi_* \nabla$ has the same geodesics as 
the Levi-Civita connection $\triangledown$ of the standard  metric
$g$ on $\RP^{2*}$, so that $\Phi^*\triangledown$ is projectively equivalent
to $\nabla$. Identifying $\RP^2$ with $\RP^{2*}$ via  any isometry
   now proves the claim. 
 \end{proof}

\begin{remark}\label{tricky}
Much the same trick used to check the regularity of $\Phi_0$
also allows one to show that  $\Psi : {\mathcal Z}\to {\mathcal N}$
is actually $C^{k+1,\alpha}$ along $Z$.
Indeed, let $\varsigma_0$, $\varsigma_1$ and $\varsigma_\infty$ be three smooth  sections 
of ${\mathcal U}\to M$ over a coordinate domain $U\subset M$ whose values are
all distinct at each point. In terms of our local coordinates $(x^1,x^2,\zeta)$, 
these correspond to three complex-valued functions $\zeta_\ell (x) = \zeta (\varsigma_\ell)$,
$\ell=0,1,\infty$,
whose values are all distinct, and never real. Set
$$\tilde{\zeta} (x,\zeta)= \frac{[\zeta - \zeta_0(x)][\zeta_\infty(x)-\zeta_1(x)]}{[\zeta_1(x) -\zeta_0(x)][\zeta_\infty(x)
 -\zeta]},$$ 
so that $\tilde{\zeta}(x,\zeta_\ell(x))=\ell$  for each $x=(x^1,x^2)$ and 
$\ell=0,1,\infty$. Choose  an inhomogeneous coordinate system on $\CP_2$
such  that $z^1(F(\Psi(\varsigma_\ell(0,0)))$, $\ell=0,1,\infty$, are 
all finite and distinct, 
and, for $x$ in a neighborhood of $0$,  
set
$$ (z^1_\ell (x), z^2_\ell(x)) =F \circ \Psi (\varsigma_\ell (x^1,x^2)), ~~ \ell = 0,1,\infty. $$
Then, in these coordinates, $F\circ\Psi$ must explicitly be given by 
$$(x,\zeta) \mapsto \left( 
\frac{\lambda z^1_0(x)+\tilde{\zeta}(x,\zeta)z^1_\infty(x)}{\lambda(x)+\tilde{\zeta}(x,\zeta)},
\frac{\lambda z^2_0(x)+\tilde{\zeta}(x,\zeta)z^2_\infty(x)}{\lambda(x)+\tilde{\zeta}(x,\zeta)}
\right),$$
where
$$
\lambda (x^1,x^2)= \frac{z^1_\infty(x)-z^1_1(x)}{z^1_1(x)-z^1_0(x)},
$$
since each $\CP_1$ fiber of ${\mathcal Z}\to M$ is sent to holomorphically to a projective line 
in $\CP_2$ by $F\circ \Psi$. If $\nabla$ is $C^{k,\alpha}$, 
this shows, albeit quite indirectly, that $\Psi$ is 
$C^{k+1,\alpha}$ on all of  ${\mathcal Z}$, and not just on 
${\mathcal U}= {\mathcal Z}-Z$. Needless to say, however, a direct analytic proof of this
fact, perhaps along the lines of \cite{eastgrah}, would be highly desirable. 
\end{remark}

If we start with a Zoll metric $h$ on $M=\RP^2$, rather than just a 
Zoll projective structure, the complex surface $\mathcal N$ comes
equipped with  a \label{conehead} 
certain additional  complex curve
${\mathcal Q} \subset {\mathcal N}$. Indeed, let  us consider the locus 
$${\mathcal C} = \{ [v] \in \PP T_\CC M ~|~ h(v,v) = 0\}, $$
where $h$ has been extended from $TM$ to $T_\CC M$ as a 
{\em complex bilinear} form, and set
$${\mathcal Q} = \Psi [{\mathcal C}].$$
 In any inhomogeneous coordinate $\zeta$
on the fiber $T_{x\CC}M$, $h(v,v)$ becomes a quadratic polynomial of
degree  $2$, and the corresponding locus in $\PP T_{x\CC}M$
thus consists of  two points, perhaps counted with multiplicity. However, since $h$ is real,
${\mathcal C}$ is invariant under complex conjugation, so a root of multiplicity 
two would have to lie in  the real slice $\PP T_x M$; but the latter is impossible, 
since $h$ is a positive-definite inner product on $T_xM$. Thus
${\mathcal C}$ intersects each fiber of $\PP T_\CC M$ in precisely 
two points, neither of which is  in $\PP TM$. 
 Indeed, if we choose to think of ${\mathcal U} = {\mathcal Z}-Z$
as the bundle of all point-wise almost-complex structures on $M$, 
$\mathcal C$ is  consists precisely of those almost-complex structures
 which are orthogonal transformations of $T_xM$ with respect to $h$; and there
are exactly two of these for each $x$, corresponding to the two possible 
orientations of $T_xM$.

Now ${\mathcal C}$
is horizontal with respect to the Levi-Civita connection $\triangledown$,
since parallel transport preserves $h$. This not only implies that 
${\mathcal C}$ meets each fiber of $\PP T_\CC M$ transversely, but also, 
more importantly, that  
there is a non-zero element  $\Xi_0$
of ${\bf D}$ which is   tangent to ${\mathcal C}$ at each point.
Thus  ${\mathcal C}$
is a complex curve in $\PP T_\CC M - \PP TM$, and its diffeomorphic image 
${\mathcal Q} = \Psi [{\mathcal C}]$ is a complex submanifold of  $\mathcal N$. 
Since $\mathcal C$ is invariant under complex conjugation, the
corresponding curve ${\mathcal Q}\subset {\mathcal N}$ is therefore
invariant under the action of $\sigma : {\mathcal N}\to {\mathcal N}$. 
Moreover, since $\mathcal C$ meets each fiber of $\PP T_\CC M$
transversely, in two points $\not\in \PP TM$, it follows that
$\mathcal Q$ meets $\Sigma_x$ transversely, in two points,
for any $x\in M$. 

Also notice  that the bundle projection $\mu : 
\PP T_\CC M\to M$ induces a $2$-to-$1$ covering map
 $\varpi: {\mathcal C}\to M\approx \RP^2$,
so ${\mathcal C}$ is therefore 
compact --- and indeed, must be diffeomorphic to  $S^2$.
Moreover, this covering map  $\varpi$ is a {\em conformal map} from the 
Riemann surface ${\mathcal C}$ to the Riemannian manifold
$(M,h)$, since 
$$\varpi_* [T^{0,1}_{[v]}{\mathcal C}]=\mbox{span} (v)\subset T_\CC M , $$
and $h(v,v)=0$.  With this observation in hand, we may now  prove the following:

\begin{thm} \label{rumple}
Let $(M,h)$ be a 
 Riemannian $2$-manifold whose geodesics are
all embedded circles of length $\pi$. 
If $M$ is not simply connected, there is a   diffeomorphism 
$\Phi : M\stackrel{\approx}{\longrightarrow} \RP^2$ such that   
$h=  \Phi^* g$,
where $g$ is the standard curvature $1$ 
Riemannian metric 
 on $\RP^2$. \end{thm}
\begin{proof}
With these hypotheses, 
the Hopf-Rinow theorem tells us  that   $M$ is  necessarily compact, since,
for any $x\in M$, 
the  closed disk of radius $\pi/2$ in  
$T_xM$ will surject  
onto $M$ under the exponential map.  Proposition \ref{others},
therefore tells us that $[\triangledown ]$ is a Zoll projective structure
on the compact surface $M$. Now assume henceforth that $M$ is not simply connected, 
We then know that $M\approx \RP^2$ by Proposition \ref{class}, 
and that $[\triangledown ]$    has conjugacy number $1$ by 
Theorem \ref{jeeves}. 

Now the proof of Theorem \ref{rigid} tells us that 
that there is a biholomorphism $F: {\mathcal N}\to \CP_2$ such that 
the  $F(\Sigma_x)$ is a projective lines $\CP_1\subset \CP_2$
for each $x\in M$, and such that $F\circ \sigma \circ F^{-1}$ is the 
complex conjugation map 
$$
[ z^1 : z^2 : z^3 ]  \mapsto  [\bar{z}^1: \bar{z}^2: \bar{z}^3 ] .
$$
Thus, $F( {\mathcal Q})$ is a non-singular compact complex curve in 
$\CP_2$ which is invariant under complex conjugation, and which 
meets certain projective lines transversely, in two points. Hence 
 $F( {\mathcal Q})$ is a non-singular conic,  and  
so is the zero locus of a quadratic polynomial 
$$0= q(z)= \sum_{j,k=1}^3q_{jk}z^jz^k.$$
But since  $F( {\mathcal Q})$ is invariant under complex conjugation,
it is also the zero locus of $\overline{q(\overline{z})}$, so that both
$$\sum_{j,k=1}^3(\Re e ~q_{jk})z^jz^k
~~~\mbox{ and }~~~
\sum_{j,k=1}^3(\Im m ~q_{jk}) z^jz^k$$
vanish along $F( {\mathcal Q})$; and at least one of these
quadratic forms is non-trivial, since $q\not\equiv 0$. 
Thus $F( {\mathcal Q})$ is the zero locus of a real 
quadratic form, represented by a real symmetric $3\times 3$ matrix 
$A= [a_{jk}]$. But any such $A$  is similar, over $GL(3, \RR  )$, to
a diagonal matrix whose  entries are all in $\{  1, 0, - 1\}$.  
On the other hand, 
since $F ({\mathcal Q}) \cap \RP^2= \emptyset$,
  the quadratic form represented by  $A$  must be 
 definite. Thus, by a suitable real change of 
coordinates,  we may arrange for
our map $F: {\mathcal N}\to \CP_2$ to send
$\mathcal Q$ to the standard conic ${\mathcal Q}_0$ given by 
$$ (z^1)^2+ (z^2)^2+ (z^3)^2 =0$$
without sacrificing any of the previously used properties of $F$.

On the other hand, we can repeat the entire construction for the 
standard metric $g$ on $\RP^2$. The map $\Phi: M\to \RP^2$
constructed in Theorem \ref{rigid} is then characterized by 
$$\Phi (x) = \tilde{x} \Longleftrightarrow F(\Sigma_x)= \tilde{F}(\tilde{\Sigma}_{\tilde{x}})$$
where untilded letters pertain to $(M,h)$ and tilded ones pertain to 
$(\RP^2, g)$. But since we have arranged for both ${\mathcal C}$ and 
$\tilde{\mathcal C}$ to map biholomorphically to ${\mathcal Q}_0\subset \CP_2$,
it follows that 
\begin{eqnarray*}
  F\left[\Psi [\varpi^{-1}(x)]\right]&=& F(\Sigma_x)\cap {\mathcal Q}_0\\
\tilde{F}\left[\tilde{\Psi} [\tilde{\varpi}^{-1}(\tilde{x})]\right]
&=& \tilde{F}(\tilde{\Sigma}_{\tilde{x}})\cap {\mathcal Q}_0.
\end{eqnarray*}
The holomorphic map $$\hat{\Phi} = \left(\left. (\tilde{F}\circ \tilde{\Psi})\right|_{\tilde{\mathcal C}}
\right)^{-1}\circ (F\circ \Psi ) :
{\mathcal C}\to \tilde{\mathcal C}$$
therefore makes the diagram 
\setlength{\unitlength}{1ex}
\begin{center}\begin{picture}(20,17)(0,3)
\put(2,17){\makebox(0,0){$\mathcal C$}}
\put(18,17){\makebox(0,0){$\tilde{\mathcal C}$}}
\put(2,5){\makebox(0,0){$M$}}
\put(18,5){\makebox(0,0){$\RP^2$}}
\put(0,12){\makebox(0,0){$\varpi$}}
\put(20,12){\makebox(0,0){$\tilde{\varpi}$}}
\put(10,6.5){\makebox(0,0){$\Phi$}}
\put(10,19){\makebox(0,0){$\hat{\Phi}$}}
\put(18,15.5){\vector(0,-1){9}}
\put(2,15.5){\vector(0,-1){9}}
\put(3.5,17){\vector(1,0){13}}
\put(3.5,5){\vector(1,0){12}}
\end{picture}\end{center}
commute, and, since $\varpi$ and $\tilde{\varpi}$ are
both conformal maps, it follows that $\Phi$ is also conformal.
In other words, 
 $\Phi^*g = e^{2u}h$ for some smooth function $u : M\to \RR$.
But the Levi-Civita connection $\tilde{\nabla}$ of $\Phi^*g$
is then related to the Levi-Civita connection $\nabla$ 
of $h$ by
$$\tilde{\nabla}_{\bf v}{\bf w}- \nabla_{\bf v}{\bf w} = du ({\bf v}) {\bf w} + du ({\bf w}) {\bf v} + 
h({\bf v},{\bf w}) \mbox{ grad}_h u .$$
However, the proof of Theorem \ref{rigid} tells us that $\tilde{\nabla}$ and
$\nabla$ are also projectively equivalent; that is,
$$\tilde{\nabla}_{\bf v}{\bf w}- \nabla_{\bf v}{\bf w} = \beta ({\bf v}) {\bf w} + \beta ({\bf w}) {\bf v} $$
for some $1$-form $\beta$. 
Thus 
$$\beta ({\bf v}) {\bf w} + \beta ({\bf w}) {\bf v} = 
du ({\bf v}) {\bf w} + du ({\bf w}) {\bf v} + 
h({\bf v},{\bf w}) \mbox{ grad}_h u$$
for all vectors ${\bf v}$ and ${\bf w}$. 
But if, for example,   we take ${\bf v}$ and ${\bf w}$ to be 
orthonormal, with 
$\beta ({\bf w}) =0$, we  then have
$ \beta ({\bf v}) ~ {\bf w}=
du ({\bf v}) ~{\bf w} + du ({\bf w})~ {\bf v}$, so that 
$du ({\bf v}) = \beta ({\bf v})$,  $du ({\bf w})=0=\beta ({\bf w})$; thus 
$du$ and $\beta$ must have the same components in 
the basis $({\bf v}, {\bf w})$, and 
hence  $\beta = du$. But if  instead we take ${\bf w}={\bf v}\neq 0$,
we  instead obtain
$$
2 ~du ({\bf v}) ~{\bf v} + |{\bf v}|^2\mbox{ grad}_h u = 2~\beta ({\bf v}) ~{\bf v}, 
$$
and the substitution $\beta = du$ then tells us that $\mbox{grad}_h u=0$.
Hence $u$ is constant. But, by hypothesis, 
 $h$ is  normalized  so that its
geodesic circles all have the same length as those of $g$. The constant
$e^{2u}$ must therefore equal $1$, and $\Phi$ is therefore an isometry 
between $(M,h)$ and $(\RP^2, g)$. 
\end{proof} 

This is essentially equivalent \cite{beszoll} to the classical Blaschke conjecture
first proved  by Leon Green \cite{grezoll} in the early 1960s. 

\begin{cor}[Blaschke Conjecture] 
Let $(M,h)$ be a compact  
 Riemannian \linebreak 
$2$-manifold for  which  the cut locus of each 
point $x\in M$ is a  one-point set $\{ x'\} \subset  M$. Then 
there is a diffeomorphism 
$\Phi : M\stackrel{\approx}{\longrightarrow} S^2$ such that   
$h= c \Phi^* g$,
where  $g$ is the standard curvature $1$ 
Riemannian metric 
 on $S^2$, and $c$ is some positive constant. \end{cor}

\begin{proof}
On a compact Riemannian manifold, any minimizing geodesic 
segment necessarily has finite length, so every arc-length-parameterized geodesic 
emanating from $x$ must arrive at the cut locus $\{ x'\}$, and must first do so 
precisely at  time  $\mbox{dist}(x,x')$. But since $x'$ represents the first
conjugate point on each geodesic leaving $x$, we see, by following these
geodesics backwards, that $x$ is an element of the cut locus of $x'$, and our
hypothesis therefore implies that the cut locus of $x'$ is {\em exactly} 
$\{ x\}$. Thus $x\mapsto x'$ is an involution $\imath : M\to M$. Moreover, every geodesic
of $M$ is a simple closed curve, and $\imath$ maps every such
geodesic circle to itself, by a rotation of $180^\circ$. In particular,  
 $\imath$ is an isometry, and  is therefore smooth. Moreover, $\mbox{dist}(x,\imath (x))$
is independent of $x$ along any particular geodesic, and thus 
is constant on $M$. Thus 
 the geodesics of the 
the quotient Riemannian metric on 
$M/\langle \imath \rangle$ are all 
simple closed  curves of equal length.
After a suitable rescaling, Theorem \ref{rumple} therefore tells us that 
the non-simply-connected Zoll manifold 
$M/\langle \imath \rangle$ becomes isometric to the standard $\RP^2$, and
hence that $M$ becomes isometric to the standard $S^2$. 
\end{proof}

\begin{remark}
Since the Christoffel symbols of the Levi-Civita connection of $h$ are
expressed  in terms of the first derivatives of $h$, Theorem \ref{rigid} 
 constructs an isometry of class $C^{k+2,\alpha}$
when we assume that $h$ is itself of class $C^{k+1,\alpha}$, $3\leq k$,
$0< \alpha < 1$.   Thus the regularity of the map $\Phi$
in Theorem \ref{rumple} is actually  
optimal, as one  certainly has every right to expect. 

It is  more important, however,  to inquire as to
 the minimal level of differentiability  
needed for our proof of Theorem \ref{rumple}.  
If we  assume that $h$ is of class $C^4$, then the proof goes through,
although the  constructed map $\Phi$ would appear only  
to be $C^4$. Nonetheless, 
 $\Phi^*g$ is still $C^3$, and its Gauss curvature
is therefore the pull-back of the Gauss curvature of $g$. This
shows any $C^4$  Zoll metric 
$h$ on $\RP^2$ must have constant curvature. 
However,    Green's proof  \cite{grezoll} actually 
draws the same conclusion even if  $h$ is merely assumed to be $C^3$. 
It would thus be extremley gratifying if  there were 
 some way of improving the present arguments so as to 
make them  work when, for example,  $[\nabla]$ is  merely asssumed to be 
of class $C^2$!
\end{remark}

\pagebreak

\section{Zoll Structures on the $2$-Sphere}
\label{zoe}

In light of our success in understanding Zoll structures of odd conjugacy number, 
it now seems reasonable to ask what our  techniques can tell us about 
the even case.  Let us 
therefore suppose that we are given a $C^3$ Zoll projective structure 
$[\nabla ]$ of even conjugacy number on a compact $2$-manifold $M$. By 
Corollary \ref{uppercase},  $M$ is then diffeomorphic to $S^2$.
Let us fix some orientation of $M$, and observe that 
$${\mathcal U}={\mathcal Z}-Z=
\PP T_\CC M - \PP TM$$
 can once again be identified with the space of 
all point-wise almost-complex structures on $M$. 
Thus 
$${\mathcal U} = {\mathcal U}_+ \cup {\mathcal U}_- , $$
where 
${\mathcal U}_+$ (respectively, ${\mathcal U}_-$) consists
of those almost-complex structures which are compatible 
(respectively, incompatible)
with the given orientation 
of $M$. These are both connected sets; indeed, either can 
be identified with the space of all point-wise conformal structures
on $M$.  Let us now consider the compact $4$-manifold-with-boundary
$${\mathcal Z}_+:={\mathcal U}_+\cup Z,$$
with   $\partial {\mathcal Z}_+ =Z$.
We can identify 
${\mathcal Z}_+$  with the non-zero, semi-positive
elements of $\odot^2T^*M$, modulo rescaling.  
Relative to some chosen `background' metric 
$h_0$ on $M\approx S^2$, we can then identify  
${\mathcal Z}_+\to M$ as the  unit disk bundle in  the traceless,
symmetric bilinear forms $\odot^2_0T^*M$. 
From a topological view-point, this 
allows us to think of   ${\mathcal Z}_+$ as
 the unique oriented $2$-disk bundle of 
Euler class $4$ over $S^2$. 

Let us now give the normal bundle $J^\parallel \ker\mu_*$
of $Z=\partial {\mathcal Z}_+$ the 
`inward pointing' orientation, and then give $\ker \mu_*$
the corresponding orientation. Having made such a choice, 
Theorem \ref{wooster} then tells us that $\nu : Z\to N$ can be 
canonically identified with the circle bundle ${\mathbb S}TN\to N$,
in such a way that $J^\parallel \ker\mu_*$ is canonically
identified with the pull-back of the (trivial) tautological 
line bundle over ${\mathbb S}TN$, meaning the 
 sub-bundle  
$L\subset \pi^* TN$, where $\pi :{\mathbb S} TN \to N$ is the canonical 
projection, whose fiber at $[{\bf v} ] \in \PP TN$
is $\mbox{span} ({\bf v})$. Now, with respect to  the canonical `outward pointing' orientation 
of $L\to {\mathbb S}TN$, let 
$L^+$ be the $[0,\infty )$-bundle consisting of vectors
which are not inward pointing. By the tubular neighborhood
theorem, $Z= \partial {\mathcal Z}_+$ has a neighborhood $\hat{\mathcal V}$ in 
${\mathcal Z}_+$ which can be identified with $L^+$
via a $C^1$ diffeomorphism whose derivative along the zero section 
of $L$ is given by our previous identification of $J^\parallel \ker\mu_*$
and $L$. But we have an obvious  $C^1$ `blowing down' map 
$\psi : L^+\to TN$, and, letting $ {\mathcal V}$ denote the 
total space of $TN$, 
 this now corresponds to
a $C^1$ map $\tilde{\psi} : \hat{\mathcal V}\to {\mathcal V}$
which is a diffeomorphism on the complement of $Z$. 
We may now define a 
 differentiable  
 $4$-manifold  
$${\mathcal N} = {\mathcal U}_+ \cup_{\tilde{\psi}} {\mathcal V}$$
by gluing together ${\mathcal U}_+$
and ${\mathcal V}=TN$ via $\tilde{\psi}$. 
 By construction, 
we   have a surjective $C^1$ 
`blowing down' map 
$$\Psi : {\mathcal Z}_+ \to {\mathcal N},$$
given by the identity on ${\mathcal U}_+$ and by $\tilde{\psi}$ on $\hat{\mathcal V}$,
so in particular we know that ${\mathcal N}$ is compact. 
Moreover, if $[\nabla ]$ is $C^k$, we can once again 
impose a `provisional' $C^{k-1}$ structure on 
$\mathcal N$ so that $\Psi$ will become a $C^{k-1}$ map.

Now $\mathcal Z$ still carries an involutive complex distribution 
$\bf D$, and the proof of Proposition \ref{key}, supplemented by the
remark on pp. \pageref{rocky}---\pageref{road},  then proves the following:

\begin{prop}\label{cle} Let $[\nabla]$ be a Zoll projective structure which is 
represented by a $C^3$ connection $\nabla$ on $M\approx S^2$. Then 
there is a unique integrable almost-complex structure $J$ on $\mathcal N$ such that 
$$\Psi_* [{\bf D}] \subset T^{0,1}({\mathcal N}, J).$$
The unique $C^\infty$  structure on ${\cal N}$
associated with its maximal  atlas of  $J$-compatible complex charts
is compatible with the previously-constructed $C^1$ structure on 
$\mathcal N$, so that  $\Psi : {\mathcal Z} \to {\mathcal N}$  remains a $C^1$ map relative to 
this smooth structure.
Moreover, if  $\nabla$ is of class $C^{2k+6}$, then $\Psi$ is $C^k$. 
\end{prop}

In order to unmask the identity of the complex surface $({\mathcal N}, J)$, 
we will now call in the heavy artillery, in the form of the 
 following fundamental result, which  is   due to Yau \cite{yau}.
We include the synopsis of a complete proof,
both as a courtesy to the reader, and for our    own enjoyment.

\begin{lem}[Yau] \label{wao}
Let ${\mathcal S}$ be a simply connected compact complex surface
with $b_2({\mathcal S})=1$. 
 Then ${\mathcal S}$ is biholomorphic to
$\CP_2$.
\end{lem}
\begin{proof}
Any   compact, oriented,   simply connected
$4$-manifold ${\mathcal S}$ has   Euler characteristic 
$\chi ({\mathcal S})  = 2+ b_2({\mathcal S})$, so that 
$\chi ({\mathcal S})  = 3$ if $b_2({\mathcal S})=1$.
On the other hand,  if $b_2({\mathcal S})=1$, the
 signature $\tau ({\mathcal S})$ is evidently 
 $\pm 1$, where the  $\pm$  sign indicates whether   the intersection form of 
${\mathcal S}$ is positive or  negative definite. 
But our $\mathcal S$ is assumed to admit a  complex structure, so its first Chern class
 has self-intersection 
$$
c_1^2 ({\mathcal S}) = 2\chi ({\mathcal S}) + 3\tau ({\mathcal S}) = 6\pm 3 > 0,
$$
 and the intersection form 
$H^2({\mathcal S}, {\ZZ}) \times H^2({\mathcal S},
{\ZZ}) 
\to \ZZ$ therefore 
cannot be negative definite. Thus  $\tau ({\mathcal S})=1$, and 
$c_1^2 ({\mathcal S})= 6+ 3 = 9$. 
Since this same calculation also shows that there is a holomorphic line bundle
of positive self-intersection, Grauert's criterion implies \cite{bpv} that 
${\mathcal S}$ is projective algebraic.  But since  $H^2({\mathcal S}, \ZZ)\subset 
H^2({\mathcal S}, \RR)\cong \RR$, 
 and 
$c_1 ({\mathcal S})\neq 0$,  this can only happen if 
 $c_1  ({\mathcal S}) = \pm [\omega ]$ for some K\"ahler form $\omega$.

Now if we had
 $c_1  ({\mathcal S}) = - [\omega ]$, the Aubin/Yau  theorem \cite{aubin,yau} 
would tell us that 
$\mathcal S$ admitted a K\"ahler-Einstein metric  of negative Ricci curvature. 
However, one has the Gauss-Bonnet-like formula
$$
 \chi - 3\tau = \frac{1}{8\pi^2} \int_{\mathcal S} \left[
3 |W_-|^2 - \frac{|\stackrel{\circ}{r}|^2}{2} 
\right]d\mu$$
for any K\"ahler metric on any compact complex surface, 
where $\stackrel{\circ}{r}$ is the trace-free  Ricci-curvature, and 
where the anti-self-dual Weyl  curvature 
$W_-$ is  the only piece of the 
curvature tensor not determined by the Ricci tensor. 
For our manifold, $\chi = 3\tau$, whereas $\stackrel{\circ}{r}$ vanishes
for any Einstein metric, so we would conclude that $W_-\equiv 0$.
 Our K\"ahler-Einstein manifold would therefore necessarily 
have negative sectional curvature,
and so would have  contractible universal cover.  
But   ${\mathcal S}$ has been assumed to be  compact and simply connected,  so this is a
contradiction. 

We must therefore have  $c_1  ({\mathcal S}) =  [\omega ]$
for some K\"ahler metric. Set $L= K^{-1/3}$, where $K=\Lambda^{2,0}$
is once again the canonical bundle,  so that $L$ is  the unique positive line bundle on 
${\mathcal S}$ with $c_1^2 (L) =1$.  By the Kodaira vanishing theorem,
$H^p ({\mathcal S}, {\mathcal O}( L))=0$ for  $p > 0$, and the Hirzebruch-Riemann-Roch
theorem therefore tells us that 
$$h^0 ({\mathcal S}, {\mathcal O}( L)) =\left\langle \left(1+\frac{c_1}{2}+\frac{c_1^2+ 
c_2}{12}\right) \exp (\frac{c_1}{3}) , [{\mathcal S}]\right\rangle 
= \frac{11}{36}c_1^2 +
\frac{1}{12} c_2 = 3.$$
Moreover, if $\Sigma \subset {\mathcal S}$ is the curve cut out by 
the vanishing of any non-trivial  holomorphic section of $L$, then,
because $L$ is positive on every curve  and satisfies $L\cdot L =1$, 
 $\Sigma$ can have only 
one irreducible component, and the zero of the section can only 
have multiplicity $1$ at a generic point of $\Sigma$.  
If
$\hat{\Sigma}$ is the normalization of $\Sigma$, the pull-backs of these
 sections  therefore give us a $2$-dimensional 
space of sections of the degree-$1$ line bundle  $L|_{\hat{\Sigma}}$.
But since $\hat{\Sigma}$ is connected, Abel's theorem tells us that 
this gives us
  a biholomorphism $\hat{\Sigma}\to \CP_1$; and   since 
this map  is induced  by  pull-backs of  sections
from ${\mathcal S}$,  
$\hat{\Sigma}\to {\mathcal S}$ is an embedding, so that 
 $\Sigma= \hat{\Sigma}$ is  a non-singular embedded  curve. 
Moreover, there is no point of $\Sigma$ at which every section 
of $L$ vanishes. This shows that  the linear system   $|L|$
has empty base locus, and 
the sections of $L$ therefore give us a 
well-defined holomorphic map
$$
F : {\mathcal S}  \to \PP [H^0({\mathcal S}, {\mathcal O}( L))^*]\cong \CP_2 .
$$
But since the inverse image of any $\CP_1\subset \CP_2$ is a smooth complex 
curve $\Sigma$ which is carried biholomorphically onto its image, this
map is a degree-$1$ holomorphic submersion, and  is therefore a biholomorphism.  
\end{proof}

Let us next recall that a differentiable 
$n$-dimensional submanifold $X$ of a complex $n$-manifold 
$(Y^{2n}, J)$ is said to be {\em totally real} if 
$T_pX\cap J(T_pX) =0$ at each $p\in X$. 
When $n=2$, which is the case of interest to us here, 
this is equivalent to the statement that $T_pX$ is never
a $1$-dimensional complex  subspace of $(T_pY, J) \cong \CC^2$. 
This  is  of course an open condition on $T_pX$; indeed, 
for $n=2$, this simply amounts to the observation that 
since $Gr_{1}(\CC^2)=\CP_1$ is a closed submanifold of 
$Gr_2(\RR^4) \cong (S^2\times S^2 ) /\ZZ_2$. 
To that  any submanifold which is $C^1$ close to a totally 
real submanifold will itself  be totally real.  

It  will also be convenient to introduce some terminology specifically tailored 
to discussions of 
differentiable embeddings of $\RP^2$ into $\CP^2$. 

\begin{defn}
A differentiable  embedding  $\jmath: \RP^2\hookrightarrow \CP^2$ will be said to be 
{\em weakly unknotted} if there exists a diffeomorphism $\mbox{\cyr   f}: \CP_2 \to \CP_2$ such that $\jmath = \mbox{\cyr   f}\circ j$, where $j: \RP^2\hookrightarrow \CP^2$ 
is the standard embedding $[x:y:z]\mapsto [x:y:z]$.
\end{defn}

\begin{remark}
By composing  with   complex conjugation $\CP_2\to \CP_2$ 
 if necessary, we may always arrange for
 $\mbox{\cyr   f}$ to  induce the identity on homology.  But since  two self-homeomorphisms of a simply connected 
compact $4$-manifold are $C^0$-isotopic iff they induce the same
maps on homology \cite{frequin}, our diffeomorphism
  $\mbox{\cyr  f}$ would then be in the identity component of 
of the homeomorphism group of $\CP_2$. Thus any weakly unknotted embedding 
of $\RP^2$ in $\CP^2$, as defined above,  may  be moved through locally flat
{\sl topological}  embeddings 
so as to ``unknot'' it into the standard $\RP^2$. {\em A priori}, however, 
it might still be impossible to carry out this unknotting process by
a path of {\em smooth} embeddings.
\end{remark}

\begin{thm} Let $[\nabla]$ be a $C^3$ Zoll projective structure on
an oriented surface 
$M\approx S^2$. Then,
up to  a projective linear transformation, the projective
structure 
 $[\nabla]$ uniquely determines a differentiable, totally real,
weakly unknotted
 embedding
of the space of geodesics $N\approx \RP^2$ into $\CP_2$. 
If $[\nabla]$ is $C^\infty$, so is the embedding. 
Moreover, the image of each of the circles $\ell_x\subset N$, $x\in M$,
 bounds a holomorphic embedding  of the 
disk $D^2\hookrightarrow \CP_2$, and the interiors
of these disks foliate the  complement $\CP_2-N$.
\end{thm}

\begin{proof} By construction, the smooth $4$-manifold ${\mathcal N}$  can be
 obtained by gluing  the unit disk bundle in $T\RP^2$
to the  Euler-class-$4$ $D^2$ bundle over $S^2$ via an orientation-reversing diffeomorphism 
of their common boundary, which is the Lens space $X=S^3/\ZZ_4$. 
However, the diffeomorphism type of the pair $({\mathcal N},N)$ 
only depends on the isotopy class of the diffeomorphism $X\to X$.
But the group of orientation-preserving diffeomorphisms of $X$ is connected
\cite{diffm3}, so it follows that the diffeotype of  the pair 
$({\mathcal N},N)$  is independent of which  Zoll projective structure 
$[\triangledown ]$ on $S^2$ we use. However, the standard structure 
$[\nabla ]$ gives us the pair 
$(\CP_2 , \RP^2 )$.  Thus 
there is a diffeomorphism $\mbox{\cyr \em f}: \CP_2 \to {\mathcal N}$
with $\mbox{\cyr \em  f}(\RP^2 )=N$.

In particular, this argument 
says that ${\mathcal N}$ is diffeomorphic to $\CP_2$. 
 Lemma \ref{wao} therefore
tells us that there is a {\em biholomorphism} $F: {\mathcal N}\to  \CP_2$, 
and this $F$ is unique modulo composition with elements of    $PSL (3, \CC )$. 
The promised embedding $N\hookrightarrow \CP_2$ is then given by 
$F|_N$, 
whereas the promised disks are the images of the the fibers of 
${\mathcal Z}_+\to M$ under $F\circ \Psi$. 
Moreover,  since the diffeomorphism $\mbox{\cyr  f}= F\circ \mbox{\cyr \em  f}: \CP_2 \to \CP_2$  
sends 
$\RP_2$ to  $F(N)$, our embedding $F|_N$
is weakly unknotted,  and we are done. 
\end{proof}

Now, in order to invert the above onstruction, let us  instead suppose that we 
are given a totally real submanifold $N\approx \RP^2$ of $\CP_2$,
and attempt to construct a suitable family of holomorphic disks $D\hookrightarrow \CP_2$
with boundary $\partial D = S^1 \hookrightarrow N$;  these circles
in $N$ will then eventually become the curves $\ell_x$ corresponding to a
Zoll projective structure on $S^2$. Our method of  accomplishing this
will be to invoke  the inverse function theorem, and so 
will apply only when the given 
 embedding $N\hookrightarrow \CP_2$
is $C^1$ close to the standard  embedding 
$\RP^2 \hookrightarrow \CP_2$. Thus, relative to a choice of tubular neighborhood,
we will 
 henceforth assume that $N$  is represented by 
by a section of the normal bundle of $\RP^2$. This allows us a
further technical simplification, since,  if such a section
has sufficiently small $C^1$-norm. Also notice that the
normal bundle of $\RP^2\subset \CP_2$ can be canonically identified,
via the complex structure, with $T\RP^2$, and so may also be identified 
with $T^*\RP^2$ by means of the standard Riemannian metric.
Thus the freedom of choosing the submanifold $N\subset \CP^2$ 
can be conveniently parameterized by the space of $1$-forms on 
on $\RP^2$ of sufficiently small $C^1$-norm. 

For the standard projective structure on $S^2$, the 
disks in question are obtained by considering those
complex projective lines $\CP_1\subset \CP_2$ which 
are complexifications of some real projective line 
$\RP^1\subset \RP^2$, and then choosing one of
the hemispheres into which such a  $\CP_1$
is divided by the corresponding $\RP^1$. 
In order to understand these disks more explicitly, let us begin with the 
 standard homogeneous coordinates
$[z_1 : z_2 : z_3 ]$ on $\CP_2$, with the usual convention that  
$\RP^2$ is represented by $z_1, z_2 , z_3$ real, 
and consider the affine chart  $(\zz_1 , \zz_2 )$ on $\CP_2$ defined by 
$$
\zz_1  =  \frac{z_1-iz_2}{z_1+iz_2} ~, \hspace{0.5in}
\zz_2  =  \frac{z_3}{z_1+iz_2}~.
$$
This chart realizes  $\RP^2-[0:0:1]$ as the M\"obius band $B\subset \CC^2$ 
given by 
$$
\zz_1\overline{\zz}_1  =  1  ~, \hspace{0.5in}
 \zz_1\overline{\zz}_2 = \zz_2 .  
$$
Note that  we may also parameterize $B$ by 
\begin{eqnarray*}
\zz_1 & = & e^{i\theta} \\
\zz_2 & = & t e^{i\theta /2} ,
\end{eqnarray*}
where  the real coordinates  $(\theta , t)$ are  best  thought of  as really 
taking values in the {\em abstract} M\"obius band $\RR^2/\ZZ$
corresponding to    the $\ZZ$-action  generated by 
$$ (\theta , t ) \mapsto (\theta + 2\pi , -t ).$$

Now the projective line $z_3=0$ in $\CP_2$ corresponds, in this picture, to 
the complex affine line $\zz_2=0$; and one hemisphere of this $\CP_1$
is  the disk $|\zz_1| \leq 0$ in this affine complex 
line, the  boundary of which is the circle 
$\theta \mapsto (e^{i\theta}, 0)$  in $B$. 
How many other ways can one  holomorphically  the disk $D\subset \CC$ 
in $\CC^2$ in such a manner that its boundary $\partial D = S^1$ both lies on
$B$, and is homotopic in $B$ to $\theta \mapsto (e^{i\theta}, 0)$? Projecting
any such disk to the $\zz_1$ axis would give a degree-$1$ holomorphic map
$D\to \CC$ with boundary map a degree-$1$ map $S^1\to S^1$, and any
such map is of course given by a M\"obius transformation
\begin{equation}
\label{moby}
\zeta \mapsto \frac{a\zeta  + b}{\overline{a}+ \overline{b}\zeta} ~, |a|^2-|b|^2 =1. 
\end{equation}
Thus, after composition with a M\"obius transfomation, 
any such disk is the graph $\zz_2 = F (\zz_1)$ of a holomorphic function $F$ 
on the unit disk $|\zz_1| \leq 1$. 
However, the requirement that $F(\partial D)$ lie in $B$ says that 
$$
 F(e^{i\theta}) = e^{i\theta} \overline{F(e^{i\theta})}.
$$ 
If $F$ has power series expansion
$$F(\zz_1 ) = \sum_{\ell =0}^\infty a_\ell \zz_1^\ell , $$
our boundary condition becomes  
$$
\sum_{\ell=0}^\infty a_\ell e^{i\ell\theta} = \sum_{\ell=-\infty}^1 \overline{a}_{-\ell+1} e^{i\ell\theta}.
$$
Hence every such disk is   the graph of an affine linear function
$$\zz_2 = a + \bar{a}\zz_1$$
restricted to the unit disk $|\zz_1|\leq 1$, where $a=a_0$. 
Each of these disks exactly represents one hemisphere of the 
projective line $\CP_1\subset \CP_2$ given by 
$$z_3= (2~\Re e ~a) ~z_1 + (-2~\Im m ~a) ~z_2,$$ 
 and the boundaries of these disks are  thus 
precisely  the real projective lines $\RP^1\subset \RP^2$
which do not pass through the point $[0:0:1]$ which was excluded by our
choice of coordinates. By considering all possible permutations of the 
homogeneous coordinates $z_1,z_2,z_3$, one obtains the entire
family of disks corresponding to the points of $S^2$ equipped with
its standard projective structure.

We now consider the problem of constructing an analogous family of 
disks with boundaries on a submanifold $N\subset \CP_2$ which is $C^1$ 
near to $\RP^2\subset \CP_2$. To do this, it is enough to completely 
analyze the corresponding problem arising when 
 intersection of the M\"obius Band $B$ and a large ball
is replaced with  a section of its normal bundle, since $N$ is covered
by a finite number of pieces of this form. 

To this end, we will begin by considering maps of the circle $S^1$
to the abstract M\"obius band $\RR^2/\ZZ$
with winding number $1$. For reasons of technical transparency,
we will consider maps of Sobolev class $L^2_k$, where $k\geq 1$.
Let us recall  that the Cauchy-Schwarz inequality 
immediately implies the Sobolev embedding theorem in this case, since
any smooth, real valued  function $f$ on the line satisfies
\begin{equation}
\label{kosher}
| f (a) - f(b) | \leq \left( \int_a^b\left|\frac{ df}{dx}\right|^2 dx\right)^{1/2} |a-b|^{1/2},
\end{equation}
whence $L^2_k (S^1)\subset C^{k-1,\frac{1}{2}}(S^1)$.
In particular,  maps from the circle of class  $L^2_k$ are continuous,
and   it  thus makes sense to talk about  winding numbers of such  maps. 
Moreover, this shows that point-wise multiplication of functions
gives us  a continuous bilinear map 
$L^2_k (S^1)\times L^2_k (S^1)\to  L^2_k(S^1)$. Also note
that  the composition of 
of any $C^k$ function with an $L^2_k$ function is again an  $L^2_k$
function.

 We will freely identify  $L^2_k(S^1)$ with the real Hilbert space of 
real-valued $L^2_k$ functions of  $\theta \in [0,2\pi ]$ 
with 
$u(\theta +2\pi ) = u(\theta  )$, and we will also need to consider the
real Hilbert space $\tilde{L}^2_k(S^1)$ of $L^2_k$ sections of the M\"obius band, 
which we may think of as functions 
 of $\theta \in [0,2\pi ]$ 
with $u(\theta +2\pi ) = -u(\theta  )$. 
Since any continuous section of the M\"obius band must have a zero, 
(\ref{kosher}) tells us that 
any $u\in \tilde{L}^2_k$, $k\geq 1$,  satisfies 
$$\sup |u| \leq \sqrt{\pi}\left( \int_0^{2\pi}\left|\frac{du}{d\theta}\right|^2d\theta \right)^{1/2}\leq \sqrt{\pi} \|u\|_{L^2_k} ,$$
so the elements $u$ of the ball of radius $R/\sqrt{\pi}$ in $\tilde{L}^2_k$ may be thought of as
defining a section $\theta \mapsto (\theta, u(\theta ))$ 
of    the {\em finite} M\"obius strip 
$$
B^R = \left( \RR \times [-R ,  R ]\right)/\ZZ ,
$$
where the $\ZZ$ action is again generated by $(\theta , t ) \mapsto (\theta + 2\pi , -t )$. 
We will use $C^k (B^R)$ to denote the real Banach space of $C^k$ real-valued
functions on this strip, and 
$$\tilde{C}^k(B^R)= \{ h: \RR \times [-R ,  R ]\stackrel{C^k}{\to} \RR ~|~ h(\theta +2\pi , -t) = 
-h(\theta , t)
 \} $$ to denote the real  Banach space
of $C^k$ sections of the non-trivial real line bundle on $B^R$, the Banach-space norms 
being of course the suprema of the absolute values of all partial derivatives of 
order $\leq k$. 

Any pair $(h_1,h_2) \in C^{k+1}(B^R) \times \tilde{C}^{k+1}(B^R)$ defines 
an embedding  $B^R\hookrightarrow \CC^2$  by 
$$(\theta , t ) \mapsto \left(e^{h_1(\theta , t )+i\theta}, [t+ ih_2(\theta , t )]e^{i\theta /2}\right),$$
and any $C^{k+1}$ submanifold $N\subset \CP_2$ which is sufficiently close to 
the standard 
$\RP^2 \subset \CP_2$ can be written as a finite union of images of 
such embeddings of finite strips via suitable systems of inhomogeneous coordinates. 
The general $L^2_k$ embedding of $S^1$
inside this strip with winding number $1$ can then be written as
$$\theta \mapsto  
 \left(e^{h_1(\theta + u_1(\theta) , u_2(\theta) )+i[\theta + u_1(\theta)]},
 [u_2(\theta)+ ih_2(\theta + u_1(\theta) , u_2 (\theta )  )]e^{i(\theta + u_1(\theta)) /2}\right)$$
for $u_1\in L^2_k(S^1)$ and $u_2\in \tilde{L}^2_k (S^1)^R$, where 
 $\tilde{L}^2_k(S^1)^R$ denotes the open ball of radius $R/\sqrt{\pi}$ 
centered at the origin in $\tilde{L}^2_k(S^1)$.
This motivates us to consider the  maps of Banach manifolds   
$${\mathcal F}_1,  {\mathcal F}_2: 
L^2_k (S^1) \times \tilde{L}^2_k(S^1)^R \times 
C^{k+\ell}(B^R) \times \tilde{C}^{k+\ell}(B^R) \longrightarrow L^2_k (S^1, \CC) \times L^2_k (S^1, \CC), $$
given by 
$$
[{\mathcal F}_1 (u_1 , u_2 , h_1 , h_2) ] (\theta ) = 
\exp \left[ h_1\Big(\theta + u_1(\theta) , u_2(\theta) \Big)+i\Big( \theta + u_1(\theta)\Big)\right]
$$
and 
$$
[{\mathcal F}_2 (u_1 , u_2 , h_1 , h_2) ] (\theta ) = 
\left[ 
u_2(\theta)+ ih_2\Big(\theta + u_1(\theta) , u_2 (\theta )  \Big)\right]  \exp \left( i\frac{\theta + u_1(\theta)}{2} \right).
$$
These maps are both  $C^\ell$; in particular, for $\ell \geq 1$  they have bounded continuous
derivatives  given by 
$$
(\dot{u}_1 , \dot{u}_2 , \dot{h}_1 , \dot{h}_2) ] 
\stackrel{{\mathcal F}_{1*}}{\longmapsto}
\Big[ \dot{h}_1(\theta +u_1, u_2) 
+\Big( i+ \frac{\partial h_1}{\partial \theta}\Big) \dot{u}_1
 + \frac{\partial h_1}{\partial t}\dot{u}_2  \Big] e^{h_1+i(\theta + u_1)}
$$
and 
$$
(\dot{u}_1 , \dot{u}_2 , \dot{h}_1 , \dot{h}_2) ] 
\stackrel{{\mathcal F}_{2*}}{\longmapsto} 
\left[ \Big( \frac{iu_2-h_2}{2}+i \frac{\partial h_2}{\partial \theta}\Big)\dot{u}_1+ 
\Big( 1+ i\frac{\partial h_2}{\partial t}\Big)\dot{u}_2   + i\dot{h}_2(\theta +u_1, u_2) \right]  e^{i(\theta + u_1) /2},
$$
where $h_1$,  $h_2$,  and their first  partial derivatives  with respect to $\theta$ and $t$
are understood to be evaluated at 
$(\theta + u_1(\theta) , u_2 (\theta )  )$, and thus are  functions of class $L^2_k$ 
which depend continuously on $(u_1,u_2,h_1,h_2)$. 
In particular, notice 
 that the derivatives of these maps at the origin  are respectively given by 
$$
[{\mathcal F}_{1*{\bf 0}} (\dot{u}_1 , \dot{u}_2 , \dot{h}_1 , \dot{h}_2) ] (\theta ) = 
\Big[ \dot{h}_1(\theta , 0)  + i\dot{u}_1 (\theta) \Big] e^{i\theta} 
$$
and 
$$
[{\mathcal F}_{2*{\bf 0}} (\dot{u}_1 , \dot{u}_2 , \dot{h}_1 , \dot{h}_2) ] (\theta ) = 
\Big[ \dot{u}_2 (\theta)  + i\dot{h}_2(\theta , 0) \Big]  e^{i\theta/2}~.
$$

Next, we introduce the orthogonal projection  
$$\Pi: L^2 (S^1, \CC) \to L^2{\downarrow}$$ 
to the closed linear subspace 
$$L^2{\downarrow}=
\left\{ \sum_{\ell <0} a_\ell e^{i\ell\theta }~\Big|~a_\ell\in \CC , ~\sum_{\ell <0} |a_\ell |^2< \infty
\right\} \subset  L^2 (S^1, \CC)$$
 of  {\em negative frequency}
functions 
given by 
$$
\Pi(\sum_{\ell=-\infty}^\infty a_\ell e^{i\ell\theta } ) = \sum_{\ell=-\infty}^{-1} a_\ell e^{i\ell\theta } .$$
This a bounded linear operator, and indeed has operator norm $1$ .  
Notice that the kernel of $\Pi$ precisely consists of
those  $L^2$ function on the circle which arise as the  boundary values
of holomorphic functions on the disk. 
Set 
$$
L^2_{k}{\downarrow}= \left\{ \sum_{\ell <0} a_\ell e^{i\ell\theta }
~\Big|~a_\ell\in \CC , ~\sum_{\ell <0} \ell^{2k}|a_\ell |^2< \infty
\right\}
=  L^2_k (S^1, \CC) \cap L^2{\downarrow}. 
$$
and  notice that 
$$\Pi: L^2_k (S^1, \CC)\to L^2_{k}{\downarrow}$$
is also bounded, and  indeed again has operator norm $1$. 

Similarly, let us define 
$$\mbox{\cyr p}: L^2_k (S^1, \CC)\to \CC$$
by 
$$
\mbox{\cyr p} (\sum_{\ell=-\infty}^\infty a_\ell e^{i\ell\theta } ) = 
a_0 .$$

\begin{remark}
The linear map  $\Pi$ is closely related to the Hilbert transform on the 
circle, and can be explicitly be realized \cite{taylor} as the singular integral operator  
$$
[\Pi (u)](\theta)  = u(\theta)  - \frac{e^{-i\theta}}{2\pi}  ~ p.v. \int_0^{2\pi} \frac{u(\phi) d\phi
 }{e^{i(\phi -\theta ) }-1}. 
$$
This can be used \cite{hilltai} to show that $\Pi$ is also bounded with
respect to in $C^{k,\alpha}$ norms. However,  
we have chosen, in the spirit of \cite{bisdisk}, to emphasize
Sobolev norms here, as this has the advantage of keeping
the technical details to a minimum. 
\end{remark}

Now, for $k, \ell \geq 1$,  consider the $C^\ell$ map 
$$
L^2_k (S^1) \times \tilde{L}^2_k(S^1)^R \times 
C^{k+\ell}(B^R) \times \tilde{C}^{k+\ell}(B^R)  \stackrel{{\mathcal F}}{\longrightarrow}  
L^2_{k}{\downarrow} \times L^2_{k}{\downarrow}  \times   C^{k+\ell}(B^R) \times 
\tilde{C}^{k+\ell}(B^R)
\times \CC \times   \CC \times  \RR 
$$
of real Banach manifolds
defined by 
$$
{\mathcal F}=( \Pi \circ {\mathcal F}_1) \times ( \Pi \circ {\mathcal F}_2) \times  
\mbox{\cyr L} \times \tilde{\mbox{\cyr L}}
\times (\mbox{\cyr p} \circ {\mathcal F}_1) \times ( \mbox{\cyr p} \circ {\mathcal F}_2)
\times 
\mbox{\cyr sh}  , 
$$
where
$$\mbox{\cyr L} : L^2_k (S^1) \times \tilde{L}^2_k(S^1)^R \times 
C^{k+\ell}(B^R) \times \tilde{C}^{k+\ell}(B^R) \longrightarrow C^{k+\ell}(B^R)$$
and 
$$\tilde{\mbox{\cyr L}} : L^2_k (S^1) \times \tilde{L}^2_k(S^1)^R \times 
C^{k+\ell}(B^R) \times \tilde{C}^{k+\ell}(B^R) \longrightarrow \tilde{C}^{k+\ell}(B^R)$$
are the factor projections, while 
$$\mbox{\cyr sh} : L^2_k (S^1) \times \tilde{L}^2_k(S^1)^R \times 
C^{k+\ell}(B^R) \times \tilde{C}^{k+\ell}(B^R) \longrightarrow \RR$$
is  given by 
$$
\mbox{\cyr sh} (u_1,u_2,h_1,h_2)= 
\frac{1}{2\pi}\int_0^{2\pi} u_1(\theta ) 
~ d\theta .
$$
Since $\Pi$, $\mbox{\cyr L}$, $\tilde{\mbox{\cyr L}}$, $\mbox{\cyr p}$ and
$\mbox{\cyr sh}$ are all bounded linear operators, 
this map is $C^1$, with  derivative  given by 
$$
{\cal F}_* = (\Pi\circ {\mathcal F}_{1*}) \times (\Pi\circ {\mathcal F}_{2*})
\times \mbox{\cyr L}\times \tilde{\mbox{\cyr L}}\times
(\mbox{\cyr p}\circ {\mathcal F}_{1*}) \times (\mbox{\cyr p}\circ {\mathcal F}_{2*})
\times 
 \mbox{\cyr sh} .$$
In particular, for any
$$
\dot{u}_1= b_0+ \sum_{\ell=1}^\infty b_\ell\cos (\ell\theta ) + c_\ell \sin (\ell \theta ) 
$$
in $L^2_k(S^1)$,
and any 
$$
\dot{u}_2= \sum_{\ell=0}^\infty \
\tilde{b}_\ell\cos \left[(\ell + {\textstyle \frac{1}{2}}) \theta \right] + \tilde{c}_\ell \sin 
\left[(\ell +  {\textstyle \frac{1}{2}})  \theta \right] 
$$
in $\tilde{L}^2_k(S^1)$, 
we see that 
 the derivative of $\mathcal F$ at the origin is explicitly by 
$$ {\cal F}_{*{\bf 0}}\left[
\begin{array}{c}
\dot{u}_1\\ 
\dot{u}_2\\ 
\dot{h}_1\\ 
\dot{h}_2
\end{array}\right]=  
\left[\begin{array}{c}
 \Pi \left(\dot{h}_1 (\theta, 0)e^{i\theta}\right)+  
\sum_{\ell =2}^\infty{ \frac{-c_{\ell}+ib_{\ell}}{2}}e^{-i(\ell -1)\theta} \\   
 \Pi \left(i\dot{h}_2 (\theta, 0)e^{i\theta /2}\right) +  \sum_{\ell =1}^\infty{ \frac{\tilde{b}_{\ell}+
i\tilde{c}_{\ell}}{2}}e^{-i\ell\theta}
\\ \dot{h}_1\\ \dot{h}_2\\   
\mbox{\cyr p} \left(\dot{h}_1 (\theta, 0)e^{i\theta}\right) + 
\frac{-c_{1}+ib_{1}}{2}\\  
\mbox{\cyr p} \left(i\dot{h}_2 (\theta, 0)e^{i\theta /2}\right) +
{\frac{\tilde{b}_0+i\tilde{c}_0}{2}} \\ 
b_0 \end{array}\right] ~.
$$
Since ${\cal F}_{*{\bf 0}}$
manifestly has bounded inverse, the Banach-space inverse function theorem \cite{schwartz}
tells us that there is an open  neighborhood $\mathfrak U$ 
of ${\bf 0}\in  L^2_k (S^1) \times \tilde{L}^2_k(S^1)^R \times 
C^{k+\ell}(B^R) \times \tilde{C}^{k+\ell}(B^R)$ and an open  neighborhood $\mathfrak V$
of ${\bf 0}\in L^2_{k}{\downarrow} \times L^2_{k}{\downarrow}  \times   C^{k+\ell}(B^R) \times \tilde{C}^{k+\ell}(B^R)
\times \CC \times \CC \times \RR$ such that 
$${\mathcal F}|_{\mathfrak U} : {\mathfrak U} {\longrightarrow}
  {\mathfrak V} $$ is a 
diffeomorphism. For any $h_1, h_2$ of sufficiently small $C^{k+\ell}$ norm, we 
therefore  obtain a $5$-parameter family of holomorphic disks
$D\to \CP_2$ with boundaries on the graph of $(h_1, h_2)$ by considering 
the unique disks with boundary values specified by 
$({\mathcal F}|_{\mathfrak U})^{-1}[  {\mathfrak V}\cap (\{(0,0,h_1,h_2)\} \times 
\CC \times \CC \times \RR )
]$. On the other hand, not all of these disks correspond to   geometrically
distinct unparameterized disks, since any parameterized  disk gives rise to 
a $3$-parameter family 
of other parameterized disks by composition with M\"obius transformations
of the form (\ref{moby}). However, we can easily kill this ``gauge freedom'' 
by   instead 
considering the $2$-parameter family of disks whose boundary values are 
 by the circles 
$$({\mathcal F}|_{\mathfrak U})^{-1} (0,0,h_1,h_2,-w^2, w, 0) 
, ~~w\in \CC , |w| < \varepsilon .$$
The other disks in our original $5$-parameter family can then all be obtained by 
composing the disks in this $2$-parameter family with M\"obius transformations. 
Notice, however, that we have now carefully constructed our disks so that their centers 
are on the complex curve
$${\mathfrak z}_1 + {\mathfrak z}_2^2=0$$
in $\CC^2$, and that our parameter $w$ exactly sweeps out a neighborhood of 
the origin in this curve. However, this curve is just an affine chart on the conic ${\mathcal 
Q}\subset \CP_2$ given by 
$$
z_1^2+z_2^2 + z_3^2 =0 .
$$
Now the subgroup $SO (3) \subset PSL (3, \CC )$ preserves both
${\mathcal Q}$ and $\RP^2\subset \CP_2$, and acts transitively 
on both ${\mathcal Q}$ and the set of real projective lines $\RP^1\subset \RP^2$. 
Thus, by 
considering only affine charts $({\mathfrak z}_1 , {\mathfrak z}_2)$
related to our original choice by the action of $SO (3)$, we can construct 
a  collection  of families of disks so that their centers  run through 
 a finite open cover of ${\mathcal Q}\approx S^2$, in  a uniform
manner depending on the submanifold $N \subset \CP_2$, thought of as 
the graph of a  section of  the normal bundle of 
$\RP^2\subset \CP_2$  of sufficiently small $C^{k+\ell}$ norm, corresponding to 
$(h_1,h_2)$ in local coordinates. 
Since ${\mathcal F}|_{\mathfrak U}$ is a  diffeomorphism,
we can also arrange that these disks coincide up to M\"obius transformations on overlaps
by at worst restricting to a smaller open set of $N$'s in the $C^{k+\ell}$ 
topology. 
 This yields  the 
following result: 

\begin{prop}
\label{smoothdisks}
 If  $N\subset \CP_2$ is the image of any embedding
$\RP^2\hookrightarrow  \CP_2$ which is sufficiently close to  the standard one in the 
$C^{2k-1}$ topology, then $N$ contains a unique family 
of embedded oriented circles $\ell_x \subset N$, $x\in S^2$, each of  which 
bounds an  embedded holomorphic 
disk $D^2\subset \CP_2$, and each of which is 
 $L^2_k$ close (and hence $C^{k-1}$
close) to the image of an  oriented real  projective line $\RP^1\hookrightarrow \RP^2$.
Moreover, if $k \geq 2$, the corresponding family of holomorphic disks   
can be realized  by a fiber-wise holomorphic,  
$C^{k-1}$ map from the unit disk bundle in the ${\mathcal O}(4)$ 
complex line bundle over $S^2= \CP_1$. 
these disks are all embedded, and their interiors foliate $\CP^2 - N$. 
\end{prop}
\begin{proof}
 Locally, our family of disks has
been found by using ${\mathcal F}^{-1}$ to construct
 a $C^{k-1}$ map from an open set  
$W\subset \CC$ to the space of $L^2_k$ maps from the circle to $N$ which bound maps 
of the $2$-disk. But, provided that $k \geq 2$, 
 the inclusion $L^2_k\hookrightarrow C^{k-1}$ is a bounded linear map, and
the maximum principle tells us that we therefore have a 
$C^{k-1}$ map from $W$  into 
the $C^{k-1}$ maps of the disk to $\CP_2$. But any such  map is 
given by a $C^{k-1}$ map $W\times D^2\to \CP_2$. 
Since we have also arranged for the centers of our disks to 
land on the conic $\mathcal Q$, our various local families 
of disks are related by M\"obius transformations which fix the origin,
and so are elements of $U(1)$; moreover, these transformations
are $C^{k-1}$ functions of our parameters, and so determine 
a $C^{k-1}$ disk bundle over ${\mathcal Q}\approx S^2$.

Now our family of disks is a $C^{k-1}$ map {\cyr f} from this
disk bundle to $\CP_2$, and sends the zero section 
to ${\mathcal Q}$. In our $({\mathfrak z}^1, {\mathfrak z}^2)$
 coordinates, each our disks is  $C^{k-1}$ close
to a disk in a complex line ${\mathfrak z}^2=a +\bar{a}{\mathfrak z}^1$.
By possibly shrinking our neighborhood of  $N$'s, we can thus arrange  that 
each is  embedded, and transverse to ${\mathcal Q}$. 
Similarly, we can arrange for the derivative of {\cyr f}  
to be non-zero everywhere, since  locally the 
map is $C^1$ close to our model example. 
Moreover,
each of our $N$'s can be obtained from $\RP^2$ by applying 
a self-diffeomorphism of $\CP_2$ which is $C^{k-1}$ close to the 
identity, and the push-forward of the 
local functions $|{\mathfrak z}^1|^2$ by these diffeomorphisms
will result in functions which are sub-harmonic on each disk of the 
family, and  the maximum principle therefore  shows each
of the disks will meet $N$ only along its boundary. 
Thus {\cyr f} gives us a proper local diffeomorphism, and hence a covering map, from 
the interior of the disk bundle to $\CP_2 - N$; but $\CP_2 - N$
is simply connected, so {\cyr f} is a diffeomorphism on the interior
of our disk bundle. 
In particular, the zero section of of our disk bundle, which  is 
sent to $\mathcal Q$, has self-intersection ${\mathcal Q}\cdot {\mathcal Q}=2^2=4$,
so our disk bundle has first Chern class $4$, and so must be $C^{k-1}$ 
isomorphic to the unit disk bundle in ${\mathcal O}(4)$. \end{proof}

Thus, we  have constructed  a family of  curves $\ell_x \subset N$,
$x\in S^2$, which bound holomorphic disks. We now wish  to  consider the 
 of curves ${\mathfrak C}_y\subset S^2$, 
$y\in N$, obtained by considering the set of all $\ell_x$'s passing through $y$,  
and we would  like to assert that these must be  the geodesics 
of a unique Zoll projective connection $[\nabla ]$ on $M=S^2$.
Our proof of this assertion will hinge on  

\begin{lem} \label{rival}
 Let $M$ be a smooth connected 
$2$-manifold,  $\varpi: {\mathcal X}\to M$  a smooth 
$\CP_1$-bundle, and   $\rho : {\mathcal X}\to {\mathcal X}$  an involution,
commuting with the projection $\varpi$, whose  fixed-point  set 
${\mathcal X}_\rho$ is an $S^1$-bundle over $M$ which disconnects
${\mathcal X}$ into two closed $2$-disk bundles ${\mathcal X}_\pm$
with common boundary ${\mathcal X}_\rho$. 
Suppose that $\Dye \subset T_\CC {\mathcal X}$ is
a    distribution of complex $2$-planes on ${\mathcal X}$
such that 
\begin{itemize}
\item 
$\rho^* \Dye = \overline{\Dye}$; 
\item   the restriction of $\Dye$ to ${\mathcal X}_+$
is $C^{k}$, $k \geq 1$, 
and  involutive; 
\item $\Dye \cap \ker \varpi_*$  is the $(0,1)$ tangent space of the $\CP_1$ fibers
of $\varpi$;  and 
\item  the restriction of $\Dye$ to a fiber of ${\mathcal X}$  has 
$c_1= -3$ with respect to the complex orientation.
\end{itemize}
Then there is a unique $C^{k-1}$
projective structure $[\nabla ]$ on $M$ such that $\Dye$ is 
obtained from the  
associated involutive distribution $\bf D$ on $\PP T_\CC M$
given by the recipe (\ref{recipe}), pulled back by a uniquely
determined  $C^{k}$ diffeomorphism $\phi: 
{\mathcal X}\to \PP T_\CC M$ which makes the diagrams 
\setlength{\unitlength}{1ex}
\begin{center}\begin{picture}(80,17)(0,3)
\put(10,17){\makebox(0,0){${\mathcal X}$}}
\put(18,19){\makebox(0,0){$\phi$}}
\put(18,5){\makebox(0,0){$M$}}
\put(26,17){\makebox(0,0){$\PP T_\CC M$}}
\put(11,15.5){\vector(2,-3){6}}
\put(25,15.5){\vector(-2,-3){6}}
\put(12,17){\vector(1,0){10}}
\put(36,10){{and}}
\put(52,17){\makebox(0,0){$\mathcal X$}}
\put(68,17){\makebox(0,0){$\PP T_\CC M$}}
\put(52,5){\makebox(0,0){$\mathcal X$}}
\put(68,5){\makebox(0,0){$\PP T_\CC M$}}
\put(50,12){\makebox(0,0){$\rho$}}
\put(70,12){\makebox(0,0){$c$}}
\put(60,6.5){\makebox(0,0){$\phi$}}
\put(60,19){\makebox(0,0){${\phi}$}}
\put(68,15.5){\vector(0,-1){9}}
\put(52,15.5){\vector(0,-1){9}}
\put(53.5,17){\vector(1,0){11}}
\put(53.5,5){\vector(1,0){11}}
\end{picture}\end{center}
commute, where $c: \PP T_\CC M \to \PP T_\CC M$ denotes 
the usual complex conjugation map. 
\end{lem}
\begin{proof}
Let us begin by noticing that, since   $\Dye= \rho^*\overline{\Dye}$ is continuous on the
closed sets ${\mathcal X}_+$ and 
${\mathcal X}_-$, it is continuous on all of $\mathcal X$. 
Also notice that the  $\phi$ makes the above diagrams commute.

Now let $L_1$ be the $(0,1)$ tangent space of the fibers. By hypothesis, 
$L_1\subset \Dye$, so that  $L_2= \Dye /L_1$ is a well
defined complex line bundle. Also notice that, since $\Dye \cap \ker \varpi_* = L_1$,
 the fibers of 
$L_2$ are carried injectively into $T_\CC M$ by $\varpi_*$. 
We may therefore define a  continuous map 
$\phi :  {\mathcal X} \to \PP T_\CC M$ by  $z\mapsto  \varpi_* (L_{2}|_z) = 
 \varpi_* (D_z)$. Now let $\zeta$ be a smooth, fiber-wise holomorphic
coordinate on $\mathcal X$, and notice that the corresponding 
vertical vector field $\partial / \partial \overline{\zeta}$
is a smooth section of $D$.  Next, near any point of the interior of ${\mathcal X}_+$,
let $\mathfrak w$ be any other local section of $\Dye$ which is linearly independent
from  $\partial / \partial \overline{\zeta}$, and then notice that the involutivity hypothesis
$[C^1 (\Dye) , C^1 (\Dye) ]\subset
C^0 (\Dye)$
tells us that 
$$\frac{\partial}{\partial \overline{\zeta}} \left( \varpi_* ({\mathfrak w}) \right)
=  \varpi_* ({\mathcal L}_{\frac{\partial}{\partial \overline{\zeta}}}{\mathfrak w})
=  \varpi_* \left(\left[  \frac{\partial}{\partial \overline{\zeta}} ,{\mathfrak w} \right]\right) 
\equiv 0 \bmod  \varpi_* ({\mathfrak w}),$$ 
so that $\phi$ is fiber-wise holomorphic on the interior ${{\mathcal X}}_+$. 
But since $\phi = c\circ \phi \circ \rho$, it then follows that $\phi$ is also 
fiber-wise holomorphic on the interior ${{\mathcal X}}_-$.
But since $\phi$ is also continuous across ${\mathcal X}_\rho = {\mathcal X}_+
\cap {\mathcal X}_-$, this implies that $\phi$ is actually fiber-wise holomorphic on
all of ${\mathcal X}$.

Now the restriction of $L_2$ to $\varpi^{-1}(x)$
is the  pull-back, via $\phi$, of the tautological 
$\O(-1)$ line bundle over $\PP (\CC\otimes T_{x} M)\cong \CP_1$.
Since $L_1$ is the $(0,1)$ tangent space of $\varpi^{-1}(x)$, and  $\varpi^{-1}(x)
\cong \CP_1$, 
 $c_1(L_1)=-2$ on any fiber of $\varpi$.  On the other hand, 
$c_1(D)=-3$ on $\varpi^{-1}(x)$, by hypothesis. Adjunction therefore
tells us that $c_1(L) = -1$ on any fiber. However,  $c_1({\mathcal O}(-1))=-1$ on 
 $\CP_1$, and we have just observed that 
the $\phi^*c_1( {\mathcal O}(-1)) = c_1(L_2)$. This shows that the fiber-wise 
degree of $\phi$ is $(-1)/(-1)=+1$. But since $\phi$ is also fiber-wise holomorphic,
it follows that $\phi$ maps each fiber of $\mathcal X$   biholomorphically
to the corresponding  fiber of $\PP T_\CC M$. 
This in turn implies that    $\phi$ is $C^k$ on all of ${\mathcal X}$,
since it sends any three pointwise-distinct local $C^k$ sections of ${\mathcal X}_+$
to  three pointwise-distinct local $C^k$ sections of $\PP T_\CC M$, and 
$\phi$  is then algebraically determined by its value along these
sections.

Let us now try to analyze the   distribution of complex $2$-planes
$\phi_* D$ on ${\mathcal Z}= \PP T_\CC M$. To this end, let us begin
by choosing an arbitrary $C^{k-1}$ torsion-free affine connection 
$\nabla_0$ on $M$, and then considering the  corresponding $C^{k-1}$ integrable 
distribution of complex $2$-planes ${\bf D}_0$ on $\mathcal Z$
given by  (\ref{recipe}). By construction, $\phi_* D$ and
${\bf D}_0$ both intersect the vertical in the $(0,1)$ tangent spaces of 
the fibers. Moreover, letting ${\bf V}^{0,1}$ denote the $(0,1)$ vertical tangent bundle of 
$\PP T_\CC M$,  ${\bf D}_0/ {\bf V}^{0,1} = (\phi_* D )/ {\bf V}^{0,1} = {\mathcal O}(-1)$, where
${\mathcal O}(-1)$ of course denotes the tautological line bundle. Thus there is 
a unique  continuous  section $\gamma$ 
 of ${\bf V}^{1,0}\otimes {\mathcal O}(1)$ 
such that ${\mathfrak w} \in {\bf D}_0$ 
iff 
${\mathfrak w} + \gamma (\pi_*{\mathfrak w} ) \in \phi_* \Dye$;
here we  have used the notation ${\bf V}^{1,0}=\overline{{\bf V}^{0,1}}$
and ${\bf V}^{1,0}\otimes {\mathcal O}(1) = {\mathcal H}om ( {\mathcal O}(-1) , {\bf V}^{1,0} )$.
Moreover, the regularity of $\Dye$ guarantees that $\gamma$ is 
$C^{k-1}$ away from the real slice $\PP TM \subset \PP T_\CC M$. 
Now, let ${\mathfrak w}$ be a $C^{ k-1}$ local section of ${\bf D}_0$ 
for which $\pi_* {\mathfrak w}$
is a fiber-wise holomorphic section of ${\mathcal O}(-1)$;
such a section may always be constructed by 
multiplying a generic section by  a  suitable complex-valued function. 
Set $f \partial/\partial \zeta = \gamma (\pi_*{\mathfrak w} )$. 
Then, away from the real slice,  the involutivity 
of $\phi_* \Dye$ and ${\bf D}_0$ then tells us  that 
$$
\left[  \frac{\partial}{\partial \overline{\zeta}} ,{\mathfrak w} \right] \equiv  0 \bmod \frac{\partial}{\partial \overline{\zeta}}
$$
and
$$
\left[  \frac{\partial}{\partial \overline{\zeta}} ,{\mathfrak w} + f  \frac{\partial}{\partial {\zeta}}\right]
 \equiv  0 \bmod \frac{\partial}{\partial \overline{\zeta}} ~, 
$$
so that
$$
 \frac{\partial f }{\partial \overline{\zeta}} \frac{\partial}{\partial {\zeta}}  \equiv  0 \bmod \frac{\partial}{\partial \overline{\zeta}} ,
$$
and hence 
$
\partial f /\partial \overline{\zeta} = 0
$.
This shows that $\gamma$ is fiber-wise holomorphic
away from the real slice. But $\gamma$ is also continuous
across the real slice. It follows that $\gamma$  is fiber-wise holomorphic
on all of $\PP T_\CC M$. 

Now any holomorphic section of $(T^{1,0}\CP_1)\otimes {\mathcal O}(1)\cong {\mathcal O}(3)$ 
arises from a unique trace-free element of $\CC^2 \otimes \odot^2 (\CC^2)^*$. 
Thus $\gamma$
is uniquely expressible as a   trace-free symmetric tensor field
$$
\mbox{\cyr  g} \in T_\CC M \otimes \odot^2 T^*_\CC M .
$$
Since $\gamma$ is $C^{k-1}$ away from the real slice, 
it follows that $\mbox{\cyr  g}$ must be $C^{k-1}$. 
Moreover, because $\phi_* \Dye$ and ${\bf D}_0$ are both 
sent to their complex conjugates by $c$, so is $\gamma$, and 
$\mbox{\cyr  g}$ is therefore  real-valued. Setting
$\nabla = \nabla_0 +  \mbox{\cyr  g}$ now gives us  a
$C^{k-1}$ symmetric affine connection on $M$ such that 
$\phi_* \Dye$ coincides with the distribution $\bf D$ defined
by
(\ref{recipe}). Since this last requirement certainly also  determines  $\nabla$ up to
projective equivalence, we are therefore done. 
\end{proof}

This allows us to finally show that our constructed families of holomorphic 
disks actually give us Zoll projective structures.

\begin{thm}
Let $N$ be any embedding of $\RP^2$ into $\CP_2$ which is 
$C^{2k+5}$ close to the standard one. Let $\{\ell_x~|~ x\in S^2\}$ be the 
constructed family of circles which bound holomorphic disks.
For each $y\in N$, set 
$${\mathfrak C}_y =\{ x\in S^2 ~|~ y\in \ell_x \}.$$
Then there is a unique $C^k$ Zoll projective 
structure $[\nabla ]$ on $S^2$ for which every ${\mathfrak C}_y$
is a geodesic.  
\end{thm}
\begin{proof}
Let ${\mathcal X}_+$ be the unit disk bundle in ${\mathcal O}(4)$, and let 
$\mathcal X$ be its double, obtained by identifying two copies
of ${\mathcal X}_+$
along their boundaries. Let ${\mathcal X}_-$ be the second copy 
of  ${\mathcal X}_+$, and let $\rho : {\mathcal X}\to {\mathcal X}$
be the smooth map which interchanges  ${\mathcal X}_+$ and  ${\mathcal X}_-$.
Notice that one may think of   ${\mathcal X}\to S^2$ as the fourth Hirzebruch surface,
and that, while  ${\mathcal X}$ is itself diffeomorphic to $S^2\times S^2$, the 
`real slice' ${\mathcal X}_\rho \to S^2$ is the circle bundle of Euler class $4$. 

Next, we consider the constructed family of holomorphic disks $\mbox{\cyr f} : {\mathcal X}_+
\to \CP_2$ with boundary on $N$. Let $\mbox{\cyr f}_*^{1,0}: T_\CC {\mathcal X}_+
\to \mbox{\cyr f}^*T^{1,0} \CP_2$ be the $(1,0)$ component of its derivative. 
Since  $\det \mbox{\cyr f}_*^{1,0}$ is $C^{k+1}$ close to the corresponding,
non-zero   expression  arising in the  model case of the linear embedding 
$\RP^2\hookrightarrow \CP_2$,  it is also non-zero for every embedding
in an appropriate neighborhood with respect to the topology in question.
Thus we may arrange for $\Dye = \ker \mbox{\cyr f}_*^{1,0}$ to be a $C^{k+1}$ 
distribution of complex $2$-planes on ${\mathcal X}_+$ for each   of the 
embeddings in question. Moreover, $\Dye$ is involutive on the 
interior of ${\mathcal X}_+$, since $\mbox{\cyr f}$ is a diffeomorphism 
there, and sends $\Dye$  to the involutive distribution  $T^{0,1}\CP_2$.

Along ${\mathcal X}_\rho = \partial {\mathcal X}_+$, note that $\Dye$ is spanned by 
$\partial/\partial \overline{\zeta}$ and the  distribution of real lines
tangent to the fibers of 
$$\mbox{\cyr f}|_{\partial {\mathcal X}_+}:  {\mathcal X}_\rho\to N.$$
We may therefore extend $\Dye$  
to  ${\mathcal X}_-$ by declaring it  equal to $\rho^*\overline{\Dye}$ on this set. 
The resulting distribution is  $C^0$ close to the one corresponding to 
the model case, and so has $c_1 (\Dye ) =-3$ on every fiber of 
$\mathcal X$. Thus the hypotheses of Lemma \ref{rival}
are all fulfilled, and we therefore obtain  a unique $C^k$ projective structure $[\nabla ]$ on 
$M=S^2$ for which $\Dye$ corresponds to $\bf D$ via $\phi$. But
$\phi$ sends ${\mathcal X}_\rho$ diffeomorphically to $\PP TM$, and 
the fibers of $\mbox{\cyr f}|_{\partial {\mathcal X}_+}$ are thereby sent to 
a foliation $\mathcal F$ of $\PP TM$ by circles which is horizontal 
with respect to $[\nabla ]$, and must coincide with  the foliation by lifted
$[\nabla ]$-geodesics.   The projective structure $[\nabla ]$ is 
therefore Zoll, and so tame by Theorem \ref{bingo}. The space of geodesics
$\tilde{N}$ of $[\nabla]$ is then a compact manifold diffeomorphic to $\RP^2$,
and comes equipped with a tautological submersion to $N$; this 
map is necessarily a covering map, and hence is a diffeomorphism by 
comparison of fundamental groups. In particular, the
${\mathfrak C}_y$ are precisely the geodesics of the constructed  
projective structure.\end{proof}

We now address 
the issue of determining when  a given  projective
structure can be represented by the Levi-Civita connection of 
a Riemannian metric. 

Suppose that $g$ is a Zoll metric on $M \approx S^2$. 
Then, in analogy with the construction on page
\pageref{conehead}, we obtain a preferred holomorphic 
curve ${\mathcal C}\subset {\mathcal Z}_+$ of genus zero
and self-intersection $4$ 
by considering $T^{0,1}M$ for the unique complex structure 
compatible with $g$ and the fixed orientation of $M$. 
The image ${\mathcal Q}= \Psi [{\mathcal C}]$ of this
Riemann surface is then an embedded, non-singular rational curve 
of self-intersection $4$ in 
in 
${\mathcal N}\cong \CP_2$, and so must be a non-singular   conic\footnote{Note
 that one thus does not  need to invoke  
 Yau's deep contribution to Theorem \ref{wao} 
 in this  Zoll metric case, insofar as the 
existence of a rational curve of positive self-intersection
 in   ${\mathcal N}$ forces this compact  complex surface
 to be  rational,   
for strictly classical reasons \cite[Proposition V.4.3]{bpv}.}. 
After a projective linear
transformation, we may thus identify ${\mathcal Q}$ with 
the smooth 
conic  given by 
$$z_1^2+z_2^2+z_3^2=0 .$$
Henceforth, we will impose this choice as a matter of convention.

Now observe that 
the Riemann surface $\mathcal C$ is one of the two connected components of 
the locus $g ({\bf v},{\bf v}) =0$ in $\PP T_\CC M$. The complement of
this locus is doubly covered by
$$UT_\CC M = \{ {\bf v} \in T_\CC M ~|~ g( {\bf v},{\bf v}) = 1 \} ~, $$
 which we will think of as a fiber-wise complexification of
the unit tangent bundle of $(M,g)$. However, $UT_\CC M$
may be canonically identified, using $g$, with
 $$UT^*_\CC M = \{ \eta \in T^*_\CC M ~|~ g^{-1}( \eta , \eta ) = 1 \} ~, $$
and we may thus equip $UT_\CC M$ with a complex-valued $2$-form 
$\Upsilon$ obtained by restricting $d\Theta$ to $UT^*_\CC M$, where
$\Theta = \sum_{j=1}^2 y_j dx^j$ is the tautological complex-valued $1$-form 
on $T^*_\CC M$. Moreover, it is not hard  to see that  ${\bf D}= \ker \Upsilon$
on $UT_\CC M$, since, taking geodesic normal coordinates 
around an arbitrary point, we have
$$g = (dx^1)^2 + (dx^2)^2 + O (|x|^2),$$
and hence 
$$\Upsilon = d \left(
\frac{dx^1 + \zeta dx^2}{\sqrt{1+\zeta^2+ O (|x|^2)}}
\right) =
\frac{d\zeta \wedge (\zeta dx^1 - dx^2)}{(1+\zeta^2)^{3/2}}+ O (x^1,x^2)~,$$
and $\ker \Upsilon$ is therefore spanned by 
$\partial /\partial \overline{\zeta}$ and 
$\Xi = \partial /\partial {x^1} + \zeta ~\partial /\partial {x^2}+ O (x^1,x^2)$.  
 Away from the real slice, $\Upsilon$ is therefore a closed 
form of type $(2,0)$ with respect to $\bf D$, and hence is holomorphic;  and,
by the last calculation,  
 $\Upsilon \otimes \Upsilon$ descends to 
$\PP T_\CC M - {\mathcal C}$ so as to have a pole of order $3$
along $\mathcal C$. On the other hand, the restriction of $\Upsilon$ to the unit circle
bundle of $M$ is real-valued, and descends to the space of oriented geodesics by  symplectic reduction \cite{beszoll,weinstein}, and
$\Upsilon$ thus gives rise to a continuous $2$-form on the double cover of 
$\CP_2 -{\mathcal Q}$ which is holomorphic on the the complement of $N$, and hence 
holomorphic everywhere. Thus $\Upsilon \otimes \Upsilon$ is a well-defined
meromorphic section of $K^2$ on $\CP_2$ with  polar locus $3{\mathcal Q}$,
and it follows that 
$$
\Upsilon = \lambda \frac{z_1 ~dz_2\wedge dz_3 + z_2 ~dz_3 \wedge dz_1 + z_3~ dz_1\wedge d_2}{
(z_1^2 + z_2^2 + z_3^2)^{3/2}
}
$$
for some constant $\lambda \in \CC$. However, we also know that 
the restriction of $\Upsilon$ to the double cover $\tilde{N}\approx S^2$ of 
$N$ is real. Since $\tilde{N}$ is homotopic to the double cover $S^2$
of the standard $\RP^2\subset \CP_2$, and since $\Upsilon$ is 
closed, 
we have
$$
\int_{\tilde{N}} \Upsilon = \int_{S^2} \lambda \frac{z_1 ~dz_2\wedge 
dz_3 + z_2 ~dz_3 \wedge dz_1 + z_3~ dz_1\wedge d_2}{
(z_1^2 + z_2^2 + z_3^2)^{3/2}} = 4\pi \lambda ~,
$$
and it follows that $\lambda$ must be real. 
Since $\Upsilon$ is real along $\tilde{N}$, we thus conclude  that 
the Riemannian condition implies that 
 $N$  is Lagrangian with respect to the 
sign-ambiguous symplectic structure
\begin{equation}
\label{sympform}
\omega = \pm \Im m \left(
\frac{z_1 ~dz_2\wedge dz_3 + z_2 ~dz_3 \wedge dz_1 + z_3~ dz_1\wedge d_2}{
(z_1^2 + z_2^2 + z_3^2)^{3/2}
}
\right) 
\end{equation}
on  $\CP_2 -{\mathcal Q}$. 

Conversely, suppose the  surface $N$ corresponding to a given projective 
structure $[\nabla ]$ on $M=S^2$ does avoid the conic $\mathcal Q$ 
and is Lagrangian with respect to the sign-ambiguous symplectic structure
$\omega$ on $\CP_2 - {\mathcal Q}$. Assume, moreover, that $\mathcal Q$ 
generates $H_2 ( \CP_2 - N, \ZZ )\cong \ZZ$, as certainly happens whenever $N$
is sufficiently close to the standard $\RP^2 \subset \CP_2$. Then, since the 
family of holomorphic disks associated with $[\nabla ]$ generate
$H_2 ( \CP_2 , N ; \ZZ )$, their Poincar\'e duals generate $H^2 (  \CP_2 - N, \ZZ )$,
and each holomorphic disk therefore has intersection number $+1$ with 
${\mathcal Q}$, and hence geometrically intersects $\mathcal Q$
transversely in a unique point. This gives us a diffeomorphism 
${\mathcal Q} \to M$, and hence fixes a conformal structure on 
$M\approx S^2$. Since we also have an orientation, the structure group 
of the circle bundle ${\mathbb S}TM = (TM - 0_M) /\RR^+$ is thereby reduced to $SO (2)$;
let $\partial /\partial \theta$ be the vertical vector field on ${\mathbb S}TM$
which generates the corresponding $SO (2)$ action. 
Now make a particular choice of $\Upsilon$ by choosing the real constant $\lambda \neq 0$, 
and pull it back 
 to the double cover ${\mathbb S}TM$ of $\PP TM$;
 this is a real-valued $2$ form on ${\mathbb S}TM$, and 
we then have a  map ${\mathbb S}TM \hookrightarrow  T^*M$
given by  $\Upsilon ( \xi, \cdot) $, the image of which is the set of 
co-vectors for a unique Riemannian metric $g$ on $M$ in the given
conformal class. For this metric, the foliation $\mathcal F$ is given by symplectic reduction,
and the geodesics of our projective structure are then exactly those of $g$.   
Thus:

\begin{thm}
Let $N \hookrightarrow \CP_2$ be a totally real  embedding of $\RP^2$ 
which corresponds to a projective structure $[\nabla ]$ on $M\approx S^2$. 
Then there is a Riemannian metric $g$ on $M$ whose Levi-Civita
connection $\triangledown$ belongs to the projective class 
$[\nabla ]$ iff, after a $PSL(3,\CC )$ transformation
of $\CP_2$, the surface $N$ avoids the conic $\mathcal Q$, and 
is Lagrangian with respect to the signed symplectic structure 
$\omega$ on $\CP_2 - {\mathcal Q}$. Moreover, such a 
Lagrangian embedding ccompletely determines the metric $g$ up to 
an overall multiplicative constant. \end{thm}

\pagebreak

\section{Concluding Remarks}

A number of important technical issues remain to be 
resoved in connection to our treatment of Zoll structures
on $S^2$. For example, while we have shown that one can 
associate a totally real  embedding of $\RP^2$ in 
$\CP_2$ with each Zoll projective connection on $S^2$, 
and that such embedded surfaces can conversely  be 
used to determine a projective connection on $S^2$, one loses
a ridiculous number of derivatives in following the 
story full circle to ones starting point. Ideally, one might hope that 
 $C^{k,\alpha}$ projective structures on $S^2$ should  exactly correspond to 
$C^{k+1,\alpha}$ surfaces  $N\subset \CP_2$. However, we are at present quite
far from being able to make such an assertion in either direction. 

What is worse, we do not at present know  that our family of disks either exists or is 
unique when $N$ is very from the the standard $\RP^2$. Nonetheless, optimism
might well be appropriate in the present instance. Indeed, let us throw
caution to the wind and hazard the following: 

\begin{conjecture}
The moduli space of Zoll metrics on $S^2$ is connected. Moreover, 
once we mark our Zoll structures by choosing an orthonormal frame
at some base-point, 
the  
moduli space of marked Zoll structures is in natural $1$-$1$ correspondence with the set of
totally real Lagrangian embeddings of $\RP^2 \hookrightarrow (\CP_2 -{\mathcal C}, \omega )$
which are homotopic to the standard embedding.
\end{conjecture}

In fact, it does not seem not hard to show that the set of $N\subset \CP_2$ 
carrying suitable families of embedded holomorphic disks
is open, but there are a numerous 
 technical difficulties  involved in trying to show that it is closed, since
sequences of embedded disks may  have singular limits,
and one tends, in the limit,
 to lose regularity  of the  dependence of families on
parameters . Moreover,  
 one would need to know
that the relevant set of   Lagrangian $\RP^2$'s
in $(\CP_2 -{\mathcal C}, \omega )$ is actually connected
for this program to ultimately succeed. Fortunately, however, 
the latter is similar to problems already solved by 
Eliashberg \cite{elih,elipol} and his co-workers, so there is ample reason to 
hope for  such a  program to be viable. 

One might also want to hazard an analogous conjecture about Zoll projective 
structures. However, this would seem to be a considerably more difficult
problem, as there is as yet no good mechanism for trying to show that 
two weakly unknotted embeddings of $\RP^2$ in $\CP_2$ are
actually isotopic. On the other hand, Gromov's $h$-principle \cite{gromh,elih} at least 
provides a rather complete reduction of  questions concerning isotopy
through totally real submanifolds to questions of isotopy in the usual, 
elementary sense.

It seems improbable that the methods we have developed here
will shed  much light on higher-dimensional Zoll manifolds, at least in the 
near term. However, 
our techniques certainly have obvious extensions which could be brought to bear on   Zoll-like 
 Lorentzian $3$-manifolds \cite{mogul},  special classes of split-signature 
Einstein manifolds  \cite{massplit} and certain problems in  Yang-Mills fields 
\cite{maswood}.
 We 
look forward to watching the further development of the present circle of 
ideas in connection with these problems.

\pagebreak 

\appendix

\section{Appendix: The axisymmetric cases}
The main result of this appendix is a formula for the general
axisymmetric Zoll projective structure close to that for the round
sphere.  We first give the formulae for the connection, and show that
it gives rise to a Zoll projective structure.  We go on to show how
these examples are calculated from the twistor correspondence.
Although this latter step is not, strictly speaking, required in the
logical structure, this was how these examples were obtained and it is
difficult to see how the formulae would be obtained otherwise, or
indeed aspects of the proof of the Zoll property.  Furthermore, it
provides a family of detailed worked examples of the construction
given in the body of the paper.

\subsection{Axisymmetric Zoll examples}
In order to introduce the formulae, we first recall Zoll's original
family of axisymmetric metrics expressed here in spherical polar coordinates  
$$
g=(F -1)^2d\phi ^2 +\sin^2\phi ~
d\theta ^2\, .
$$
This metric is
the same as that given by (\ref{ansatz}),  after the
coordinate transformation $z=\cos\phi$ and the substitution
$F(\phi)=-f(\cos\phi)$.  

We will express the general Zoll projective sructure in terms of the
difference between a compatible affine connection and the metric
connection of the above metric.  
Consider the orthonormal frame 
$$
({\bf e}_1,{\bf e}_2)=(\frac{1}{F
-1}\p \phi ,\frac{1}{\sin\phi }\p \theta )
$$
and dual co-frame
$(\theta^1,\theta^2)=(({F -1})d\phi, \sin\phi d\theta)$.
In this frame, it is straightforward to calculate that the connection
1-form is 
$$
\omega=\frac{\cot\phi}{F -1}\theta^2\, .
$$ 
The associated Levi-Civita connection, $\nabla^g$, gives the most
general axisymmetric Zoll projective structure that is compatible with
a metric (at least close to the round metric).

In general, a compatible torsion-free affine connection for a
projective structure can be given by a connection $\nabla$ such that,
with $$
\gamma_{ij}^k=\left\langle\theta^k,(\nabla_i-\nabla^g_i){\bf e}_j\right\rangle \, ,
$$ $\gamma_{ij}^k$ is symmetric on the $ij$ indices (so that $\nabla$
is torsion-free) and trace free; this last condition corresponds to
fixing the connection in the projective equivalence class by requiring
that it preserve the metric volume form.  

The general axisymmetric Zoll projective structure close to the round
sphere is obtained with the choice 
$$
\gamma_{11}^i=0\, , \quad
\gamma_{22}^1=-\frac{h^2\cot\phi}{F-1}
\, , \quad
\gamma_{21}^1= \frac1{3(F-1)} \left(\frac{\partial h}{\partial \phi} -\frac{2h}{\sin\phi\cos\phi}
\right)
$$
where all the other components of $\gamma_{ij}^k$ are determined by
the trace and symmetry conditions and $h=h(\phi)$ is a smooth function
of $\phi$ vanishing in some small neighborhood of $0$ and $\pi$ and
odd under $\phi\rightarrow \pi-\phi$.

This information can be encapsulated in the geodesic spray.  
This is the vector field on the projective tangent bundle $PTS^2$ that at
$(v,x)\in PT_xS^2$ is the horizontal lift of the vector
$v$ at $x$.  We parametrize the fibre of $PTS^2$ by $\zeta$
corresponding to the vector ${\bf e}_1+\zeta {\bf e}_2$.  Then the geodesic spray
from the projective structure above is given by
\begin{equation}\label{zollspray}
\Xi
=\p \phi  + \left(\frac{F-1}{\sin\phi}\right)\zeta \p \theta 
-\zeta\left((1+\zeta^2(1+h^2))\cot\phi   -
\zeta \left(\frac{\partial h}{\partial \phi } -
\frac{2h}{\sin\phi \cos\phi }\right) \right) \p \zeta 
\end{equation}
on $PTS^2$ defines a Zoll projective structure if the smooth functions
$F(\phi)$ and $h=h(\phi)$ are respectively odd and even under $\phi
\rightarrow\pi-\phi$.  For regularity at $\phi=\pi/2$, we further
require that $h$ should vanish (and hence to second order) at $\phi
=\pi/2$.  For regularity at $\phi=0,\pi$, we assume that $F$ and $h$
vanish in some small neighborhood of these values. This is actually stronger
than necessary,  but makes the proof of the Zoll property more
straightforward; the minimal requirement would be to just stipulate  that they be smooth
functions of $\cos\phi$ that vanish at $\phi=0$ and $\pi$.

The metric case occurs when $h=0$, and in this case
there is the preferred overall scaling factor that gives the arc-length
parameterisation; this arises on dividing by $(F
-1)\sqrt{(1+\zeta^2)}$.  To see that the above gives a multiple of the
geodesic spray in this case, coordinatize the tangent bundle by
$(\mu_1,\mu_2)\rightarrow \mu_1 {\bf e}_1+\mu_2 {\bf e}_2$.  Then the horizontal
lift of ${\bf e}_1$ is just ${\bf e}_1$ since $\omega({\bf e}_1)=0$ and the horizontal
lift of ${\bf e}_2$ is ${\bf e}_2-\omega({\bf e}_2) (\mu_1\p {\mu_2} -
\mu_2\p {\mu_1})$.  Thus, using the affine coordinate
$\zeta=\mu_2/\mu_1$ on the projective tangent bundle, the geodesic
spray will be
$$
{\bf e}_1+\zeta \left({\bf e}_2 -(1+\zeta^2)\frac{\cot\theta}{F
-1}\p \zeta \right)
$$
and this can be seen to be proportional to the formula given above
when $h=0$ as required.  If we wish to normalize the horizontal part
to have unit length, then we must divide by $\sqrt{(1+\zeta^2)}$ and
this will give the overall factor required to give proper length
parameterisation.

We first give a direct proof of the Zoll property, and then in the
subsequent sections we show how the formula arises from the twistor
construction.  (The direct proof of the Zoll property below in fact
will use equations arising in the twistor derivation below, but it is
easily checked that these follow directly from the form of the
geodesic spray above.  It is difficult, however, to see how they might
have been anticipated without the twistor construction.)

\begin{thm}
Equation(\ref{zollspray}) defines a Zoll projective structure for all
smooth odd functions $F$ and even functions $h$ with $h(\pi/2)=0$
and both $h$ and $F$ vanishing in some neighborhood of
$\phi=0$.  
\end{thm}

The proof is divided into two parts.  We first 
analyze the flow of the projection of the geodesic spray under
$q:PTS^2\rightarrow \RP^1\times[0,\pi]$,
$q(\zeta,\phi,\theta)=(\zeta,\phi)$, to the space
of orbits of $\partial/\partial\theta$ in $PTS^2$.  We show first that the orbits
of the projected flow are circles, and secondly that the lifts of these
to orbits of the full geodesic spray are also circles in $PTS^2$.

\smallskip

\noindent
1)
We first study the integral curves of
$q_*\Xi$ for $(\zeta,\phi )\in \RP^1\times[0,\pi/2]$.
Introduce the angular coordinate $\psi\in
[0,\pi)$ on $\RP^1$ by $\zeta=\tan\psi$ so that $\psi$ is a
smooth coordinate near $\zeta=\infty$.  Then, the flow
becomes
\begin{eqnarray}
\dot \phi& =&\sin\phi \cos\psi\, ,\nonumber \\ 
\dot\psi &=&
-\sin\psi\left((1+h^2\sin^2\psi)\cos\phi -\cos\psi \sin\psi \left(\frac{\partial h}{ \partial \phi}
\sin\phi  -\frac{2h}{\cos\phi} \right)\right)\label{flow}
\end{eqnarray}
where $\dot\phi=d\phi/d t$ for the time parameter $t$ along the flow
defined by
$$
\frac d{dt}=\sin\phi\cos\psi q_*\Xi\, .
$$
These additional factors yield a smooth flow by inspection noting in
particular that our requirement that $h(\pi/2)=0$ implies that
$h/\cos\phi$ is smooth.  Note also that this flow is invariant under
the reflection in $\phi=\pi/2$: $(\psi,\phi,t)\rightarrow
(\psi,\pi-\phi,-t)$. 

[The perceptive reader might have noticed that the direction of the flow
changes sign across the identification of $\psi=\pi$ with $\psi=0$.
Although the flow defines a smooth distribution in the projective tangent
bundle away from the fixed points, to obtain a flow with continuous
direction, we would need to work on the double cover obtained by
factoring the tangent bundle by the positive scalings.  This will not
be a problem in the following as $\psi=0$ or $\pi$ is a flow line.]

\begin{lem}
The flow of $q_*\Xi$ has fixed points at $(\psi,\phi)=(0,0), (0,\pi)$
and $(\pi/2,\pi/2)$.  The integral
curves of $q_*\Xi$ are smoothly embedded curves in $(\zeta,\phi )\in
\RP^1\times[0,\pi]$ on which, for $\phi\in[0,\pi/2]$ (resp.\
$\phi\in[\pi/2,\pi]$) the coordinate $\phi$ decreases (resp.\
increases) from $\pi/2$ to a unique minimum (resp.\ maximum) value and
then increases (resp.\ decreases) again to $\pi/2$.  The extrema occur
when $\psi=\pi/2$.
\end{lem}

The fixed points are where both the right hand sides vanish, so that,
from $\dot\phi=0$ we obtain either $\phi=0,\pi$ or $\psi=\pi/2$.
At $\phi=0,\pi$, $h=0$ and so we find $\dot\psi=\mp\sin\psi$, i.e., a
fixed point at $\psi=0(=\zeta)$.  At $\psi=\pi/2$, we find that
$\dot\psi=0$ iff $\cos\phi=0$, i.e., $\phi=\pi/2$.

It is clear from the first of equations (\ref{flow}) that for
$\phi\in(0,\pi)$, $\dot\phi$ only vanishes when $\psi=\pi/2$.  The
second derivative at $\psi=\pi/2$ can be calculated to give
$$
\frac{\partial^2\phi}{\partial\psi^2}= (1+h^2)\cot\phi
$$
and it can be seen that this second derivative $\partial^2\phi/\partial\psi^2$ is positive
for $\phi \in(0,\pi/2)$ and so this
must be a minimum.  Similarly on $\phi\in(\pi/2,\pi)$, $\phi$ can only
be a maximum at a stationary point.  Thus, on an integral curve in
$\phi\in(0,\pi/2)$,
$\phi $ will descend to a unique minimum value, at which
$\zeta=\infty$ and then increase again.  
$\Box$

\smallskip

The key issue now is as to whether we can make these integral curves  
join up into a circle.  Firstly  note that $\psi=0$ and $\phi=0,\pi$
are all flow lines, and these are the only flow lines limiting onto
the fixed points $(0,0)$ and $(0,\pi)$ as we have assumed that $h=0$
in a neighborhood of $\phi=0$ and of $\pi$, and this means that the flow
lines in those neighborhoods are precisely those of the flat case,
and these are precisely the level curves of $\sin\phi\sin\psi$.

Let us suppose that a curve starts at some 
value of $\psi\in (0,\pi/2)$.  Then $\phi$ will descend to a minimum
and either (a) increase up to $\pi/2$ again, or (b) the minimum will be
$\phi=0$.  In case (a), the reflection of the orbit under the
involution $(\psi,\phi,t)\rightarrow (\psi,\pi-\phi,-t)$ will be an
orbit in $\phi\in (\pi/2,\pi)$ and this will join up to make a
circular orbit.  Case (b) will be the case $\psi=0$ since the orbit
must intersect $\phi=0$ at $\psi=0$, since the complement of that
point in $\phi=0$ is a regular orbit on its own, but the only orbit in
a neighborhood of $\phi=0$ that intersects this fixed point is
$\psi=0$ (or $\phi=0$).  

\begin{figure}[ht]
\caption{The flow diagram for the projected flow}\label{flowdiagram}
\begin{center}
\includegraphics[width=10cm,height=10cm,angle=0]{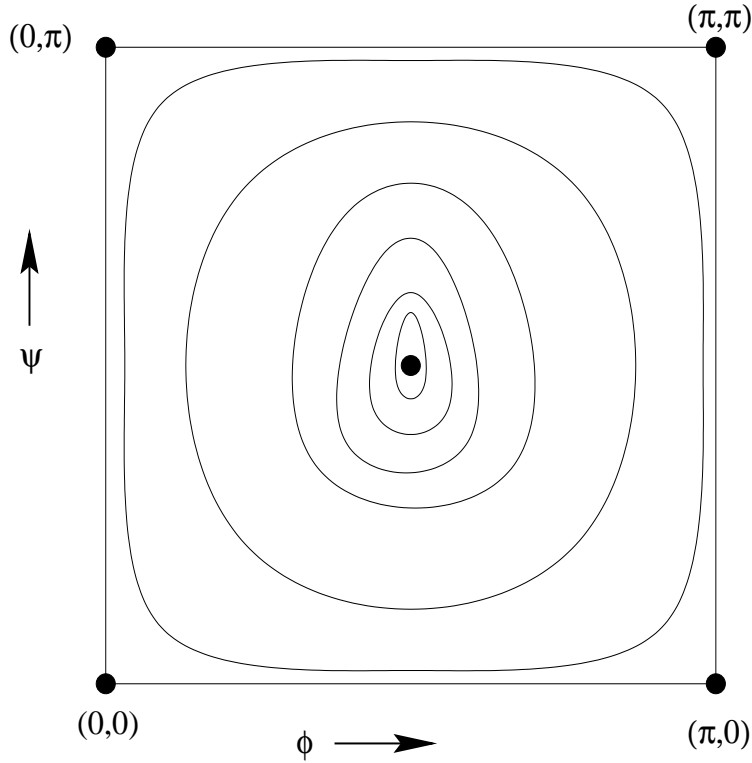}
\end{center}
\end{figure}

Thus, all the orbits of the flow are circles, except the above
mentioned fixed points and special orbits that limit onto the fixed
points; this gives the flow diagram
\ref{flowdiagram}. 

\medskip

2)  We now wish to show that these orbits in the $(\psi,\phi)$ plane
    only lift to give closed $S^1$ orbits in the full projective tangent
    bundle of the sphere.  In the above coordinates, the equation for
    $\theta$ will become
\begin{equation}\label{zolllift}
\dot\theta=(F -1)\sin\psi\, .
\end{equation}
In order for the geodesics to be circles, we need to prove that the
integral of the right hand side around an integral curve of $q_*\Xi$
is 0 modulo $2\pi$ for each integral curve.  The first and second
terms in the right hand side of equation (\ref{zolllift}) are
respectively odd and even under $\theta\rightarrow \pi-\theta$.  Since
the integral curves of $q_*\Xi$ are even, the first part will
automatically 
integrate to zero.  We need to show, then, that the second part
will in fact integrate to $0$ modulo $2\pi$ on all integral curves.

To integrate $\dot\theta=-\sin\psi$, from equations
(\ref{liftspray}) and (\ref{beta}) in the lifting part of the twistor
construction, we note that, with $1/a=(h-i)|\cos\phi|$
$$
\omega= \frac12 \arg {\frac{(1-\zeta/\bar a)}{1-\zeta/a}}
$$
satisfies
$$
\dot\omega=-\sin\psi\, .
$$
(We leave it to the assiduous reader to show that equations
(\ref{liftspray}) and (\ref{beta}) follow independently of the twistor
construction.)  Thus, $\theta=\omega$ is the solution to the even part
of the $\theta$ flow.  However, $\omega$ is the argument of a single
valued complex function on the $(\psi,\phi)$--plane, and so, when we
do a complete circuit around an integral curve of $q_*\Xi$ returning
to our original point, the argument must return to zero modulo
$2\pi$. $\Box$

\subsection{The twistor construction in the axisymmetric case}
In \S\ref{homcase} we study the structure of the action of
axisymmetry on the twistor space and the correspondence for the round
metric.  In \S\ref{defms} we give the axisymmetric deformations of
the real slice.  The subsequent subsection \S\ref{constrdisc} is
devoted to constructing the holomorphic disks, and then finally in
\S\ref{constrproj} the associated projective struture is constructed.

\subsubsection{The round sphere}\label{homcase}
We consider the action of the standard rotation on $\RR^3$, its
complexified action on $\C^3$ and induced action on $\CP^2$.  With
coordinates $(z,\tilde z, z_0)$, the $S^1$ action is generated by the real
part of the holomorphic vector field
$$
\p \theta =i(z\p z-\tilde z\p {\tilde z})\, ,
$$ 
where $\RR^3$ is taken to be $\tilde z=\bar z$ and $z_0=\bar z_0$.  If
we remove the the fixed points $(1,0,0)$, $(0,1,0)$ and $(0,0,1)$, the
generic orbits form the pencil of conics $(1-w)z_0^2=w z \tilde z$
that are tangent to the line $z=0$ at $(0,1,0)$ and also to the line
$\tilde z =0$ at $(1,0,0)$.  The degenerate orbits consist of the
double line $z_0=0$ at $w=0$ and the pair of lines $z=0$ and $\tilde
z=0$ at $w=1$.  They determine a fibration of
$\CP^2-\{(1,0,0),(0,1,0)\}$ over $\CP^1$ with affine coordinate $w$,
and, away from the exceptional fibres at $w=0,1$, we can coordinatise
$\CP^2$ with $(w,\xi)=(z_0^2/(z_0^2+z\tilde z), z/z_0)$.
In these coordinates $\p \theta =i\xi\p \xi$.

The real slice, $\RP^2$, is given by $w\in
[0,1]$ and $|\xi|^2=-1 + 1/w $. Note that the orbit $z_0=0$ intersects
$\RP^2$ in a real line, whereas the orbit $\{z=0\} \cup \{\tilde z=0\}$
intersects $\RP^2$ in a single point.  All the other real orbits are
contractible circles in $\RP^2$.

Introduce spherical polar coordinates $(\phi ,\theta )$ on $S^2$ so that
the symmetry is $\p \theta $.  We can coordinatize the fibres of the
projective tangent bundle by $\zeta$ so that $\zeta$ corresponds
to the vector $\p \phi  + \frac{\zeta}{ \sin \phi } \p \theta $.
(These coordinates will then only break down at the fixed points.)
The lines in $\CP^2$ corresponding to points of $S^2$ are $2
z_0=\tan\phi  (e^{i\theta }z+e^{-i\theta } \tilde z)$.  In terms of
$\zeta$, and the coordinates $(w,\xi)$ on $\CP^2$, the holomorphic
disks  are the images of the upper-half plane in $\zeta$ under
\begin{equation}\label{rounddiscs}
w=\frac{\zeta^2\sin^2\phi 
}{1+\zeta^2}
\, , \qquad \xi=e^{i\theta }
\frac{\zeta \cos\phi  +i}{\zeta\sin\phi }
\, ,
\end{equation}
and when $\zeta$ is real the image lies in $\RP^2$.\footnote{A
global and invariant formulation can be obtained in index notation by
letting $z_i$, $i=1,\cdots 3$ be homogeneous coordinates on $\CP^2$,
and $x^i$ coordinates on $\RR^3$, then the open disk in $\CP^2$
corresponding to $x^i$ on $S^2$ is given by the condition that $iz_i\bar
z_j\varepsilon^{ijk}$ be a positive multiple of $x^k$.}

It is worth noting for later use that, on these disks, $\zeta=
\sqrt{w/(\sin^2\phi -w)}$ defines the square root in the
upper-half plane.

In these coordinates, the geodesic spray takes the form:
$$
\Xi=\p \phi  +\frac{\zeta}{\sin\phi }\p \theta 
-\cot\phi 
(1+\zeta^2)\zeta\p \zeta\, .
$$

It should also be noted that the conserved quantity associated to the
axial symmetry $\p \theta $ and metric $g=d\phi ^2+\sin^2\phi  d\theta ^2$
is
$$
\frac{g(\p \theta ,\p \phi  +\frac{\zeta}{\sin\phi }\p \theta )}
{\sqrt{g(\p \phi  +\frac{\zeta}{\sin\phi }\p \theta ),\p \phi 
+\frac{\zeta}{\sin\phi }\p \theta )}} =\sqrt w\, .
$$
This formula can also be derived intrinsically on $\CP^2$; namely, $\sqrt w$
is the Hamiltonian for $\p \phi $ using the symplectic form
associated to the conic ${\mathcal Q}$  defined in
equation 
(\ref{sympform}).

\subsubsection{Deformation of the real slice}\label{defms}
We will represent a circle invariant deformed embedding of $\RP^2$
into $\CP^2$ as the set given by $$ w=\gamma(\phi )\, , \mbox{ and }
|\xi|^2= e^{g(\phi )} \left|\frac{1-\gamma(\phi )}{\gamma(\phi )}\right| $$ for
$\phi\in [0,\pi/2]$.  Here $g$ is a smooth real function with compact
support in $(0,\pi/2)$ and $\gamma:[0,\pi/2]\rightarrow\C$ is a smooth
embedded curve from $w=0$ to $w=1$ such that $\gamma(\phi )=\sin^2\phi
$ on the complement of some compact subset of $(0,\pi/2]$. 

In the homogeneous case, $\gamma(\phi)=\sin^2\phi$, and $g=0$.  The
compact support of the deviation from the homogeneous case will
guarantee smoothness of this deformation near the degenerate fibre
$z_0=0$.  In particular, the embedding of $\RP^2$ into $\CP^2$ near
the fixed line $z_0=0$, is the same as the canonical embedding, and so
the holomorphic disks  near those at $z_0=0$ will be those above in
equation (\ref{rounddiscs}) and so we will not need to concern
ourselves with singular behaviour there.

These assumptions amount to the assumption that our $S^1$--invariant
Zoll projective structure on $S^2$ will have two fixed points
corresponding to $\phi =0,\pi$ in a neighborhood of each of which
the projective structure will be that of the round sphere, and exactly 
one of the $S^1$ orbits will be a geodesic with $\phi =\pi/2$.  

In the metric case we will have that $\gamma(\phi )=\sin^2\phi $
since the square of the conserved quantity is determined by the
geometry of the action on $\CP^2$ relative to its fixed symplectic
structure.  It will necessarily be equal to $w$, and will be real on
the real slice.  The nontrivial information in this case is contained
only in the function $g(\phi )$.

For later convenience, we extend $\gamma$ and $g$ to $\phi \in
[0,\pi]$ by $\gamma(\phi )=\gamma(\pi-\phi )$ and $g(\phi )=g(\pi-\phi
)$. The data of the location of the deformation of $\RP^2$ could be
represented more economically by expressing the curve $\gamma$ as a
graph of the imaginary part over the real interval $[0,\pi/2]$.
However the formulation above will allow us to make a convenient
choice of the coordinate $\phi$ later.

\subsubsection{Construction of the holomorphic disks}\label{constrdisc}
The problem of finding the deformed disks  with boundary on the
deformed real slice decomposes into two parts: firstly that of
finding the projection of the disk to the $w$--Riemann sphere with
boundary on the projection of the real slice, the curve $\gamma$, and
secondly, the problem of lifting the disk to $\CP^2$.

\medskip

\noindent
1) The projected disks  must have their boundary on some subinterval of
the curve $\gamma$.  This subinterval must include the end at
$\phi=0$: this end corresponds to the line $z_0=0$ and each boundary
of a disk must be homologous to this line, but because these are all
generators of the homlogy of $\RP^2$, they must intersect each other
at least once.

Thus, the first task is to find, for each $\phi \in [0,\pi/2]$, a map
$\zeta \rightarrow w(\zeta,\phi )$ from the upper half plane into the
$w$-Riemann sphere such that the boundary of the disk is mapped to the
image of the interval $[0,\phi ]$ under $\gamma$.

To analyze this, first consider the conformal map $$ w\rightarrow
v(w,\phi )=\sqrt{\frac{w}{\gamma(\phi )-w}}\, , $$ where we fix the
branch of the square root by requiring that, near $w=0$,
$v\sqrt{\gamma(\phi )}$ lies in the upper-half-plane (there is no
obstruction to choosing $\sqrt{\gamma(\phi )}$ as $\phi $ varies so
that it is positive for small $\phi $).  In the $v$-Riemann sphere,
the image of $\gamma([0,\phi ])$ is a continuously differentiable
embedded circle tangent to $\sqrt{\gamma(\phi )}\times$ the real axis
at the origin and passing through the point $v=\infty$.  It will be
smooth except possibly at $0$ and $\infty$.  Thus the branch defined
above is well defined and determines a region $V_\phi $ in the
$v$-plane as the image of the complement of $\gamma([0,\phi ])$.

By the Riemann mapping theorem there will exist a conformal map from
the upper-half-plane in $\zeta$ to $V_\phi $ and hence to the
complement of $\gamma([0,\phi ])$ in the $w$-Riemann sphere.  It will
be smooth with non-vanishing derivative up to and including the
boundary on the $v$-Riemann sphere except possibly at $0$ and $\infty$
where it is nevertheless guaranteed to be continuous \cite[p. 340]{tayapp}.
  It is worth emphasizing that while Proposition \ref{smoothdisks}
guarantees that the disks will be smoothly embedded in $\CP^2$, but
they will be tangent to the fibres of the projection along the orbits
of the complexified axisymmetry at $w=0$ and $\gamma(\phi)$.  Hence,
the projection of the disks to the $w$-Riemann sphere will be smooth
up to $\gamma([0,\phi ])$ except at the points $0$ and $\gamma(\phi )$
which will be ramification points of order 2.  Using a Mobius
transformation of the upper-half plane to itself, this map
$w(\zeta,\phi )$ can be chosen so that $$ w(\zeta,\phi )=\zeta^2
\sin^2\phi + O(\zeta^3)\, , $$ at $\zeta=0$ and $w(\zeta,\phi
)=\gamma(\phi ) - k(\phi )\gamma'(\phi ) \zeta^{-2} + O(\zeta^{-3})$
at $\zeta=\infty$ for some real $k(\phi )>0$.

For later use we define the function $s(\zeta,\phi )$
for $\zeta\in\RR$, $\phi \in [0,\pi]$ by the condition that
$$
\gamma(s(\zeta,\phi )) = w(\zeta,\phi )\, .
$$ 
In the following we extend both $w(\zeta,\phi )$ and $s(\zeta,\phi
)$ to $\phi \in[0,\pi]$ so that they are even functions under $\phi
\rightarrow \pi-\phi $.

\medskip

\noindent
2) We now wish to find the lift of these conformal mappings to disks
in $\CP^2$ with boundary on the deformed real slice.  To do this we
need to obtain $\xi(\zeta,\phi ,\theta )$ holomorphic on the
upper-half-plane in $\zeta$ such that, for $\zeta\in\RR$, 
$$
|\xi(\zeta,\phi ,\theta )
|^2 = e^{g(s(\zeta,\phi)
)}\left|\frac{1-\gamma(s(\zeta,\phi ))} {\gamma(s(\zeta,\phi) )}\right|\, .  
$$
By symmetry we must have $\xi(\zeta,\phi ,\theta )=e^{i\theta
}\xi(\zeta,\phi ,0)$.

The orbits of the complexified axisymmetry corresponding to $w\neq
0,1$ are regular orbits.  Thus for $w(\zeta,\phi)\neq 0,1$, the lift
$\xi(\zeta,\phi ,\theta )$ cannot meet $\xi=0$ or $\infty$ since
$\xi=0$ is part of the orbit $w=1$ and $\xi=\infty$ is the orbit
$w=0$.  However, as $w\rightarrow 0$ we must have, by the above
condition on the real slice, $|\xi|^2\rightarrow
|(1-w)/w|\rightarrow\infty$.  Furthermore, if $\phi =\pi/2$, $w=1$ is
a real point on the boundary of the conformal mapping and must
therefore lift to the real point $\xi=\tilde\xi=0$.  Conversely, at
$w=1$, but $\phi \neq \pi/2$, the point $w=1$ is not a real point on
the disk and so we cannot have both $\xi=0$ and $\tilde\xi=0$.  Hence
either we will have $\xi=0$ and $\tilde \xi\neq 0$, or $\xi\neq 0$ and
$\tilde \xi= 0$.  We can therefore assume that, by continuity from the
round sphere case, $\xi\neq 0$ for $\phi \in [0,\pi/2)$, and
$\tilde\xi\neq 0$ for $\phi \in (\pi/2,\pi]$.

By taking logs, the problem of lifting the conformal maps to disks in
$\CP^2$, can be reduced to an abelian problem.  However, we cannot
proceed completely naively as we will still have $\xi\rightarrow
\infty$ as $w\rightarrow 0$, although we can guarantee that either
$xi$ or $\tilde\xi$ will be non-vanishing. We work first on $\phi \in
(0,\pi/2)$ so that $\xi\neq 0$, and divide that problem into a part
that is regular on taking logs, and one that can be handled
explicitly.  Set $$
\xi(\zeta,\phi ,\theta )=e^{i\theta  +
G(\zeta,\phi )}\Gamma(\zeta,\phi )
$$
then we wish to find $G(\zeta,\phi )$ that is
holomorphic for $\Im m\zeta >0$ such that for $\zeta$ real
$$
\Re eG(\zeta ,\phi )=g(s(\zeta,\phi ))
$$
and similarly we wish to find $\Gamma(\zeta,\phi )$ holomorphic on the
upper half plane  in $\zeta$, such that for $\zeta$ real
$$
|\Gamma(\zeta,\phi )|^2= 
\left|\frac{1-\gamma(s(\zeta,\phi )}{\gamma(s(\zeta,\phi ))}\right|\, .
$$
The first problem is solved in a standard way by a contour integral
along the real axis
$$
G(\zeta,\phi )=\frac{1}{2\pi i}
 \oint \frac{\Re eG(\mu,\phi )}{\mu -\zeta}
d \mu - \frac{1}{2\pi i}
 P.V. \int \frac{\Re eG(\mu,\phi )}{\mu }d \mu 
$$
where the purpose of the last term is to remove the ambiguity
associated with the addition of a constant (in $\zeta$ but perhaps
with $\phi $--dependence) to the
imaginary part of $G$.  This choice ensures $\Im mG(0,\phi )=0$.

The problem for $\Gamma$ cannot be solved so simply in the above way.
First we define the complex function $a(\phi )$ in the upper half
plane by the condition $w(a,\phi )=1$, i.e., the image in the
$\zeta$ plane of $w=1$.  Then the function
$$
\Gamma(\zeta,\phi )=i\sqrt{\frac{(1-\zeta/\bar a)}{(1-\zeta
/a)}\frac{(1-w)}{w}} 
$$
makes sense for $\zeta$ in the upper half plane  since the function whose root is
taken does not vanish on the upper half plane .  We choose the branch for the square
root that tends towards $i/\zeta\sin\phi $ as $\phi $ and
$\zeta$ tend to zero.  Then $\Gamma$ as defined is non-vanishing,
holomorphic in the upper half plane  and has the required modulus when
$\zeta\in\RR$ as then $|(1-\zeta/\bar a)/(1-\zeta/a)|=1$.

For $\phi \in [\pi/2,\pi]$ we work with $\tilde\xi$ as that will be
non-zero on this interval.  However, 
$$
|\tilde \xi(\zeta,\phi ,\theta )|^2 =\frac{|1-w|^2}{|w^2\xi^2|}=
 e^{-g(s(\zeta,\phi )}\frac{|1-\gamma(s,\zeta,\phi )|}
{|\gamma(s,\zeta,\phi )|}\, . 
$$
and so the solution will be
$$
\tilde\xi=e^{-i\theta -G(\zeta,\pi-\phi )}\Gamma(\zeta,\pi-\phi )\, ,
$$
where the $\Gamma$ and $G$ are the functions obtained above.

\subsubsection{Construction of the projective structure}\label{constrproj}
To reconstruct the corresponding projective connection on $S^2$, we
wish to construct the vector field determining the geodesic spray on
the correspondence space, $PTS^2$.  We use coordinates $(\phi ,\theta )$
on $S^2$, and $\zeta\in\RR$ on the fibres of $PTS^2$.  
We construct the geodesic spray $\Xi$ in two steps:

\medskip

\noindent
1) Under the projection $q:(\zeta, \phi ,\theta )\rightarrow
(\zeta, \phi )$, $\Xi$ projects to $q_*
\Xi=\p \phi -p(\zeta,\phi)\p \zeta$ for some $p(\zeta,\phi)$.  The
function $w$ is constant along the geodesic spray so that
$q_*\Xi w=0$ which gives $p=\partial_ \phi w/\partial_ \zeta w$.

When $\zeta\in\RR$, $w=\gamma(s(\zeta,\phi )$, so 
$$
p(\phi ,\zeta)=\frac{\gamma' \partial s/\partial \phi  }{\gamma' \partial s/\partial \zeta }=
\frac{\partial s/\partial \phi  }{ \partial s /\partial\zeta }
$$
is real.  Thus  $p$ can be extended meromorphically
over the $\zeta$ Riemann sphere by defining it in the lower-half
plane to be 
the complex conjugate of the pullback under $\zeta\rightarrow
\bar\zeta$.  The fact that it is real for $\zeta\in\RR$ ensures
continuity and hence holomorphy there.  It does, however, have simple
poles at $\zeta=0,\infty$ as $\p \zeta w$ has simple zeroes
there. However, the chosen form at $\zeta=0$ implies that in fact, $p$
vanishes at $\zeta=0$.  Thus, since $p\p \zeta$ is globally holomorphic
except a simple
pole at $\zeta=\infty$ (as a vector field on the Riemann sphere), zero at $\zeta=0$ and real for
$\zeta$ real, we can write 
$$
p=\zeta(\Gamma_2\zeta^2+\Gamma_1\zeta +\cot\phi )
$$
where $\Gamma_1$ and $\Gamma_2$ are real functions of $\phi $ and the
$\cot\phi $ follows from the expansion at $\zeta=0$.

Note here that since $w$ and $s$ are even functions under
$\phi \rightarrow \pi-\phi $, $\Gamma_2$, $\Gamma_1$ and
$\cot\phi $ are odd as they involve the $\phi $ derivatives of $s$.

We will need the fact later that $\Gamma_1$ and $\Gamma_2$ can be
expressed in terms of $a(\phi )$ and its first derivative by using
the condition
\begin{equation}\label{aeq}
q_*\Xi (\zeta -a)|_{\zeta=a}=0
\end{equation}
which follows from the fact that $\zeta=a$ corresponds to $w=1$
which is a holomorphic curve in $\CP^2$.  This yields the equation
$$
\p \phi  a+ a(\Gamma_2 a^2 +\Gamma_1 a +\cot\phi )=0
$$
and this together with its complex conjugate yields
\begin{equation}\label{connection}
\Gamma_2=\frac{\sin\phi }{\bar a -a}\p \phi \left(\frac{a-\bar
a}{|a|^2\sin\phi }\right) \, , 
\quad \Gamma_1=\frac{\sin\phi }{\bar a-a}\left(\bar
a\p \phi \left(\frac {1}{a\sin\phi }\right) - c.c.\right)\, .
\end{equation}
The number of free functions here is two: either the pair $\Gamma_1$
and $\Gamma_2$ or, equivalently, the real and imaginary parts of $a$.
This is to be compared to the one free function we have in the data of
the curve $\gamma(\phi )$ in the reduced twistor space and the second
free function we have in choosing the coordinate $\phi $, which, up to
now, has been arbitrary (at least away from $\phi =0,\pi/2$).  We will
fix this coordinate freedom subsequently.

\medskip

\noindent
2) The next step is to lift $q_*\Xi$ to the vector field $\Xi$ 
   on the full correspondence space $PTS^2$ that annihilates also $\xi$
   or equivalently $\tilde \xi$.   We will have
$$
\Xi= q_*\Xi - \frac{(q_*\Xi \xi)}{\partial_ \theta  \xi}\p \theta  
=q_*\Xi - \frac{(q_*\Xi \xi)}{i\xi}\p \theta  
=q_*\Xi +i(q_*\Xi \log\xi)\p \theta  
\, .
$$
In order to proceed further, note that the coefficient of $\p \theta $ is
$iq_*\Xi\log\xi$, and this is
   (a) holomorphic over upper-half-plane in $\zeta$, and 
(b) is real for $\zeta\in\RR$
since the imaginary part of
$$
i q_*\Xi \log \xi|_{\zeta=\bar \zeta }= i\left(\p \phi  -
\frac{\partial s/\partial \phi }{\partial s/\partial\zeta 
}\frac{\partial}{\partial\zeta}\right)\log \xi
$$
is just $q_*\Xi \log |\xi|$ but $\log|\xi|=\Re e\log \Gamma + \Re
G$ is a function of $\zeta$ and $\phi$ only through $s$, and such
functions of $s$ alone are annihilated by
$q_*\Xi$ by construction.  Thus, the imaginary part of the right
hand side of the above equation vanishes for $\zeta\in\RR$.  Hence,
we can extend it meromorphically over the $\zeta$--Riemann sphere by
setting it to be the complex conjugate of the pullback under
$\zeta\rightarrow \bar\zeta$ for $\Im\zeta<0$ and noting that
reality at $\zeta\in\RR$ implies continuity and hence holomorphy
across the real axis.

The function $iq_*\Xi \log\xi$ divides into two parts:
$$
iq_*\Xi \log\xi= iq_*\Xi G(\zeta,\phi )+ iq_*\Xi
\log\Gamma \, ,
$$
and since $w$ is constant along $q_*\Xi$, the second part reduces
to
$$
iq_*\Xi \log\Gamma = \frac{i}{2}q_*\Xi  \log \frac
{1-\zeta/\bar a}{1-\zeta/a} \, .
$$
They are both holomorphic on the full $\zeta$ sphere, except with
poles at $\zeta=\infty$ since $q_*\Xi$ has one there.  However,
they will also have a simple zero at $\zeta=0$ since the imaginary
parts of $G$ and the above expression for $iq_*\Xi\log \Gamma$
vanish there by construction.  (The possible apparent poles in
$iq_*\Xi\log\Gamma$ are removable as a consequence of equation
\ref{aeq}.) Therefore 
\begin{equation}\label{liftspray}
iq_*\Xi G=\frac{F(\phi)}{\sin\phi} \zeta \, , \quad\mbox{ and }\quad
iq_*\Xi\log\Gamma=\beta(\phi )\zeta
\end{equation}
for some real functions
$F$ and $\beta$ and the
geodesic spray is 
$$
\Xi=\p \phi  +\left(\frac{F}{\sin\phi} +\beta\right) \zeta\p \theta  -
\zeta(\Gamma_2\zeta^2+\Gamma_1\zeta +\cot\phi )\p \zeta\, .
$$
Using the above and equations (\ref{aeq}) and (\ref{connection}) we
calculate directly that
$$
\beta=-\Gamma_2 \Im ma\, .
$$

When $\phi \in[\pi/2,\pi]$ we should note first that $G$ and $\Gamma$
are even functions under $\phi \rightarrow\pi-\phi $.  Hence,
$F$ and $\beta$ are, as defined, odd functions.  However, there
is a further sign change on using $\tilde\xi$ instead of $\xi$ for
$\beta$ which yields an {\em even} contribution for $\beta$ and odd
for $F$ and $p$, i.e., for $\phi \in[\pi/2,\pi]$ 
$$
\Xi=\p \phi  +\left(-\frac{F(\pi-\phi )}{\sin\phi}
+ \beta(\pi-\phi )\right) \zeta\p \theta  
+p(\pi-\phi , \zeta) \p \zeta\, .
$$

\medskip

We now fix the choice of the coordinate $\phi $ which up to now has
been arbitrary except near $\phi =0$
and $\pi/2$.  We do this by imposing 
$$
\Im m\frac 1a=  - |\cos\phi |\, 
$$
(note that $a$ must always be in the upper half plane , and must be even
under $\phi \rightarrow \pi-\phi $).  This gives
\begin{equation}\label{beta}
\beta=-1/\sin\phi 
\end{equation}
Introduce the function $h(\phi)$ by 
$$
\Re e\frac 1a=h|\cos\phi |
$$ and this leads to the formulae
$$
\Gamma_1=-\p \phi  h
+\frac{2h}{\sin\phi \cos\phi } \, , \quad
\Gamma_2= \cot\phi  \left(1+h^2\right)\, .
$$
This leads to our final formula for the geodesic spray 
\begin{equation}\label{zollspray2}
\Xi
=\p \phi  + \frac{F-1}{\sin\phi}\zeta \p \theta 
-\left((1+ \zeta^2+\zeta^2h^2)\cot\phi   - 
\zeta \left(\frac{\partial h}{\partial \phi}   -
\frac{2h}{\sin\phi \cos\phi } \right)\right)\zeta \p \zeta 
\end{equation}
where $F$ must be
odd under $\phi \rightarrow\pi-\phi $ and $h$ must be even.  For
regularity, $h$
should vanish to second order at $\phi =\pi/2$.  From the assumption
that the twistor data was zero in some small neighborhood of the fixed line
$z_0=0$, we also deduce that the functions $h$ and $F$ should
vanish in some small neighborhood of $\phi =0, \pi$.
This is the formula that leads to the expressions given at the beginning of
this appendix.

\vfill

\noindent 
{\sc 
Department of Mathematics, SUNY, Stony Brook, NY 11794-3651 USA\\
The Mathematical Institute, 24-29 St Giles,
Oxford OX1 3LB,  England}

\bigskip

\bigskip

\noindent 
{\bf Acknowledgments.}
The first author would like to thank Denny Hill, 
Dusa McDuff, and Dennis Sullivan for helpful conversations, as well as  
Bob Gompf and Yasha Eliashberg for some helpful e-mail. 
The second author  would like to thank Mike Eastwood and Rafe Mazzeo for useful
discussions, and MSRI for its hospitality during the early stages of the writing of  this paper.

\pagebreak

\end{document}